\RequirePackage{fix-cm}
\documentclass[smallextended,dvipsnames]{svjour3}       
\smartqed  
\usepackage{graphicx}
\usepackage{amsmath,amsfonts}
\usepackage{tikz}
\usepackage{float}
\usepackage{subfig}
\usepackage{multirow}
\usepackage{xcolor}
\usepackage{footnote}

\usepackage{environ}         
\usepackage{etoolbox}        
\usepackage[pdftex]{hyperref}

\newlength{\myl}
\let\origequation=\equation
\let\origendequation=\endequation

\RenewEnviron{equation}{
  \settowidth{\myl}{$\BODY$}                       
  \origequation
  \ifdimcomp{\the\linewidth}{>}{\the\myl}
  {\ensuremath{\BODY}}                             
  {\resizebox{\linewidth}{!}{\ensuremath{\BODY}}}  
  \origendequation
}

\usepackage[sort&compress,square,comma,authoryear]{natbib}
\setcitestyle{numbers}
\begin{document}

\title{
An explicit exponential time integrator based on Faber polynomials and its application to seismic wave modelling
}
\titlerunning{Faber exponential integration and its application to seismic waves}

\author{Fernando V. Ravelo         \and
        Pedro S. Peixoto \and Martin Schreiber 
        }


\institute{
    Fernando V. Ravelo (ORCID: 0000-0003-1867-2123) \email{fernanvr@ime.usp.br} \and 
    Pedro S. Peixoto (ORCID: 0000-0003-2358-3221) \email{ppeixoto@usp.br}
    \at Instituto de Matem\'atica e Estat\'istica, Universidade de S\~ao Paulo, Brazil
    \and
    Martin Schreiber (ORCID: 0000-0002-4430-6779) \email{martin.schreiber@univ-grenoble-alpes.fr}
    \at Université Grenoble Alpes / Laboratoire Jean Kuntzmann, France
    \at Technical University of Munich, Germany
}

\date{Received: date / Accepted: date}

\maketitle

\begin{abstract}
Exponential time integrators have been applied successfully in several physics-related differential equations. However, their application in hyperbolic systems with absorbing boundaries, like the ones arising in seismic imaging, still lacks theoretical and experimental investigations.

The present work conducts an in-depth study of exponential integration using Faber polynomials, consisting of a generalization of a popular exponential method that uses Chebyshev polynomials.
This allows solving non-symmetric operators that emerge from classic seismic wave propagation problems with absorbing boundaries.

Theoretical as well as numerical results are presented for Faber approximations. One of the theoretical contributions is the proposal of a sharp bound for the approximation error of the exponential of a normal matrix. We also show the practical importance of determining an optimal ellipse encompassing the full spectrum of the discrete operator, in order to ensure and enhance convergence of the Faber exponential series. Furthermore, based on estimates of the spectrum of the discrete operator of the wave equations with a widely used absorbing boundary method, we numerically investigate stability, dispersion, convergence and computational efficiency of the Faber exponential scheme.

Overall, we conclude that the method is suitable for seismic wave problems and can provide accurate results with large time step sizes, with computational efficiency increasing with the increase of the approximation degree.\\

\textbf{Exponential integrator $\cdot$ Faber polynomial $\cdot$ Wave equation $\cdot$ Absorbing boundaries $\cdot$ Eigenvalue distribution }\\

\end{abstract}

\section{Introduction}\label{sec_intro}

The production of images from the subsurface using elastic waves is an important and challenging problem in geophysics \citep{ikelle2018introduction}. From a mathematical perspective, this process is divided into two sub-problems: the direct problem, consisting of solving the wave propagation equations; and the adjoint/reverse problem, which is part of the optimization problem that generally makes use of the adjoint method. Both tasks rely on the numerical solution of the propagating waves in the medium of interest and, due to the usual requirement of repeatedly having to solve the wave equations in the inverse problem, the overall computational cost and memory requirements tends to be very large. Additionally, this problem usually requires high-order numerical methods for accurate representations of wave dispersion, to ensure adequate assessment of subsurface media interfaces \citep{jing2019highly,wilcox2010high}. Naturally, in order to obtain more accurate representations of the subsurface, better algorithms are desired, posing a relevant challenge for numerical method development.

A particular class of methods, known as time exponential integrators, have been shown to outperform classical schemes in several partial differential equation (PDE) models from several applied areas in terms of accuracy and computational performance (e.g., \citep{loffeld2013comparative,pototschnig2009time,cohen2017exponential,brachet2022comparison,zhang2014nearly,schreiber2019exponential}). Such methods are known to accurately describe the solution of linear waves and, therefore, seem like natural candidates for seismic wave propagation problems. Overall, the study of such schemes with respect to solutions of the wave equations is not yet well established for seismic applications and only a few studies have been conducted in this direction \citep{kole2003solving,zhang2014nearly}. In  \citet{kole2003solving}, an explicit exponential integrator of arbitrary order in time using Chebyshev polynomials is applied to the elastic wave equations, but without the possibility to consider absorbing boundary conditions, hence rendering them inadequate for practical scenarios. In \citet{zhang2014nearly}, an implicit time exponential integrator of low order in time is developed to solve the wave equation with absorbing boundary conditions. Nonetheless, although implicit schemes have good stability properties allowing larger time steps, they are very costly and the numerical error is not granted to be small. 

Explicit exponential integrators based on Chebyshev polynomials have been recurrently mentioned by some authors as an efficient scheme of rapid convergence and easy implementation \cite{bergamaschi2000efficient,hochbruck2010exponential,ramadan2017exponential}, with further hardware-aware possible optimizations (c.f. \citet{Huber2021} for exploiting caches). They have been used to solve the wave equations \citep{kole2003solving}, obtaining not only high accuracy in the approximations, but also being able to compute the solution using remarkably large time-steps.
However, no analysis of the computational efficiency of the approximation was performed. Also, the proposed Chebyshev exponential was designed to be used only when the discrete system defines a symmetric or an antisymmetric operator. Nonetheless, in seismic wave propagation, an infinite domain is imitated using a limited area domain alongside with absorbing boundary conditions, to avoid spurious reflections. Such non-reflective boundary conditions usually break the anti-symmetry of the usual wave operator, therefore making Chebyshev exponential time integration of very limited applicability for seismic problems.
A generalization to non-symmetric (or non-anti-symmetric) operators has been proposed in \citet{bergamaschi2003efficient} based on Faber polynomials. They, in fact, showed that Faber polynomials are a generalization of Chebyshev polynomials, obtained by stretching and displacing Chebyshev polynomials which are able to handle systems with non-symmetric matrices. Moreover, the method was compared with a Krylov subspace algorithm, resulting in comparable accuracy and computational efficiency.
Withal, the proposed method is designed for and used to solve advection-diffusion equations but not wave equations.
Alongside, as in \citet{kole2003solving}, there was no analysis of the optimal polynomial degree looking towards computational efficiency and competitiveness with respect to traditional methods.

In this work, we target the investigation of the theoretical and numerical properties of Faber exponential methods and its application to seismic wave propagation problems. First, we develop the time exponential integrators based on Faber polynomials to solve acoustic and elastic wave equations with a commonly used absorbing boundary condition, the Perfectly Matching Layers (PML) \citep{assi2017compact}. Then, considering the optimization technique of scaling and squaring, we establish an error bound of the numerical approximation of the Faber exponential for operators given by a normal matrix, which improves the one reported in a previous study by \citet{bergamaschi2000efficient}. Furthermore, as a strategy to optimize the convergence of the Faber exponential series, we analyze the eigenvalues' distribution of the discrete and the continuous operators of the wave equations for different formulations, equation parameters, and spatial dimensions and propose sharp limits on the spectrum. A thorough numerical investigation of the stability and dispersion properties is also performed, calculating the Courant-Friedrichs-Lewy (CFL) number and the dispersion error for several polynomial degrees, and searching for optimal polynomials. We also study the performance of Faber approximations in practice, by solving the wave equation within several numerical examples with different levels of complexity. To ensure robustness in our results, we investigate stability, dispersion, convergence, and efficiency by varying all the model parameters in the continuum and the discrete models. Overall, this work contributes to be a first description and investigation of the viability of Faber exponential integration for classic seismic wave propagation problems, including absorbing layers, indicating the competitiveness of the scheme in terms of accuracy and computational cost.

This paper is organized as follows. In the next section, we introduce the basic formulations used to solve a general problem using time exponential integrators, with subsequent focus on the use of Faber polynomials on the exponential approximations. We also include an analysis of error bound estimates for system of normal matrices. In Section \ref{sec_application}, we present the wave equations with PML absorbing boundaries condition, with their variations, the numerical considerations, and the experiments we use through the work. The study of the spectrum of the discrete spatial operators of the wave equations is in Section \ref{sec_spectrum}, where we also derived formulas for the convex hull of the respective problem's eigenvalues.
Section \ref{sec_stability_dispersion} contains the theoretical stability and dispersion error analysis, together with an estimate of the method efficiency according to these criteria. Further on, we carry out several numerical experiments in Section \ref{sec_convergence_efficiency}, compare the accuracy and efficiency for different degrees of approximation across the seven numerical experiments designed in Section \ref{sec_application}.  Section \ref{sec_conclussion} concludes the work and summarizes the principal results.

\section{Time exponential integrators and Faber polynomials}\label{sec_faber}

In this section, we briefly present the theory of time exponential integrators, followed by a description of Faber polynomials. We also discuss some aspects of Faber polynomials, such as the error bounds of the approximation, and the dependence of the convergence on the region where the polynomials are defined.

\subsection{Exponential integrator perspective}

As described in \citet{hochbruck2010exponential}, and \citet{al2011computing}, exponential integrators are a class of time integrating methods used to solve ordinary differential equations of first order in time,
\begin{equation}\label{eq_1st_ord}
    \frac{d\boldsymbol{u}(t)}{dt}=\boldsymbol{H}\boldsymbol{u}(t)+\boldsymbol{f}(t,\boldsymbol{u}(t)),\quad \boldsymbol{u}(t_0)=\boldsymbol{u}_0,
\end{equation}
where $\boldsymbol{u}(t),\boldsymbol{u}_0\in\mathbb{C}^n$, $\boldsymbol{H}\in\mathbb{C}^{n\times n}$,  $\boldsymbol{f}$ may be a non-linear function, with \linebreak $\boldsymbol{f}:\mathbb{R}\times\mathbb{C}^{n}\rightarrow \mathbb{C}^n$, and the matrix $\boldsymbol{H}$ is usually a discretization of a continuous operator (originated from a partial differential equation, for instance).

Then, an exponential integrator is defined as an approximation of the semi-analytic solution of the constants variation formula,
\begin{equation}\label{eq_const_var}
     \boldsymbol{u}(t)=e^{(t-t_0)\boldsymbol{H}}\boldsymbol{u}_0+\int\limits_{t_0}^te^{(t-\tau)\boldsymbol{H}}\boldsymbol{f}(\tau,\boldsymbol{u}(\tau))d\tau.
\end{equation}
Depending on how the integral term  in \eqref{eq_const_var} is approximated, different types of time exponential integrators  are defined (see \citet{hochbruck2010exponential} for a review). For our particular application to seismic imaging, $\boldsymbol{f}$ represents the source term, which we will utilize later in its Taylor expanded form. In such case, \citet{higham2005scaling} shows that \eqref{eq_const_var} can be transformed into the calculation of the exponential of a slightly larger matrix,
\begin{equation*}
    \tilde{\boldsymbol{H}}=\begin{pmatrix}
    \boldsymbol{H}& \boldsymbol{W}\\
    \boldsymbol{0} & \boldsymbol{J}_{p-1}
    \end{pmatrix},
\end{equation*}
where the columns of the matrix $\boldsymbol{W}$ are formed by the values of the function $\boldsymbol{f}$, and the approximations of the first $p-1$ derivatives of $\boldsymbol{f}$ at $t_0$, and $\boldsymbol{J}_{p-1}$ is a square matrix of dimensions $p\times p$ with value one in the upper diagonal and all the other elements zero. Consequently, the information about the source term is contained in the matrix $\boldsymbol{W}$, and the solution of the system can be written as
\begin{equation}\label{eq_matrix_ampli}
    \boldsymbol{u}(t)=\begin{bmatrix}\boldsymbol{I}_{n\times n}&\boldsymbol{0}\end{bmatrix}e^{(t-t_0)\tilde{\boldsymbol{H}}}\begin{bmatrix}\boldsymbol{u}_0\\\boldsymbol{e}_p\end{bmatrix},
\end{equation}
where $\boldsymbol{e}_p\in\mathbb{R}^p$ is a vector with zero in its firsts $p-1$ elements and one in its last element, and $\boldsymbol{I}_{n\times n}$ is the identity matrix of dimension $n$.

As with classical methods, we define a time-step size $\Delta t$ and calculate the solution at $t_k=t_{k-1}+\Delta t$ as
\begin{equation*}
    \boldsymbol{u}(t_k)=\begin{bmatrix}\boldsymbol{I}_{n\times n}&0\end{bmatrix}e^{\Delta t\tilde{\boldsymbol{H}}}\begin{bmatrix}\boldsymbol{u}(t_{k-1})\\\boldsymbol{e}_p\end{bmatrix},
\end{equation*}
and the sub-matrix $\boldsymbol{W}$ of $\tilde{\boldsymbol{H}}$ now relates to an evaluation of $\boldsymbol{f}$ and its derivatives at time $t_{k-1}$.

The calculation of the matrix exponential is one of the core steps of exponential integrators and several families of methods to approximate the application of the exponential of a matrix $\boldsymbol{H}$ onto a vector $\boldsymbol{u}_0$ have been proposed (see \citet{moler2003nineteen}). We continue by investigating such evaluations with the help of Faber's polynomials, where $\boldsymbol{H}$ would be $\tilde{\boldsymbol{H}}$ if a source term is include.

In this research, we will use a Taylor expansion of the term $\boldsymbol{f}$ of the same order as the matrix exponential approximation. Moreover, taking into account that $\boldsymbol{f}$ is usually a known function, the partial derivatives of the Taylor expansion will be calculated symbolically. This ensures that the effective order of the temporal scheme is in agreement between exponential and source parts, without much added computational cost, since the order of the expansion ($p$) is usually orders of magnitude smaller than the dimension of the discrete operator $\boldsymbol{H}$.

\subsection{Faber polynomials}\label{sec_faber_polynomials}

Given a degree $j$, and a square matrix $\boldsymbol{H}$, Faber's polynomials are defined as $\boldsymbol{F}_j(\boldsymbol{H})$, with
\begin{align}
\boldsymbol{F}_0(\boldsymbol{H})&=\boldsymbol{I}_{n\times n},\quad \boldsymbol{F}_1(\boldsymbol{H})=\boldsymbol{H}/\gamma-c_0\boldsymbol{I}_{n\times n},\\  \boldsymbol{F}_2(\boldsymbol{H})&=\boldsymbol{F}_1(\boldsymbol{H})\boldsymbol{F}_1(\boldsymbol{H})-2c_1\boldsymbol{I}_{n\times n},\label{eq_faber_pol_0}\\
\boldsymbol{F}_j(\boldsymbol{H})&=\boldsymbol{F}_1(\boldsymbol{H})\boldsymbol{F}_{j-1}(\boldsymbol{H})-c_1\boldsymbol{F}_{j-2}(\boldsymbol{H}),\quad j\geq 3,\label{eq_faber_pol}
\end{align} 
with
\begin{equation}
\gamma=\frac{a+b}{2},\quad c_0=\frac{d}{\gamma},\quad c_1=\frac{c_f^2}{4\gamma^2},
\end{equation}
where the parameters $a$, $b$, $c_f$, and $d$, are set according to the spectrum of the operator $\boldsymbol{H}$. Faber's polynomials are considered a generalization of Chebyshev's polynomials because they are stretched and translated Chebyshev's polynomials \citep{starke1993hybrid}.

In \citet{bergamaschi2000efficient}, it is shown that given an ellipse $\mathcal{E}(d,c_f,a)$, symmetric with respect to the real axis (with center $d$, focuses $d\pm c_f$, and semi-major axis $a$), if the spectrum $\sigma(\boldsymbol{H})$ of the matrix $\boldsymbol{H}$ is contained in the ellipse $\mathcal{E}(d,c_f,a)$, the Faber partial sums
\begin{equation}\label{eq_faber_sum}
S_m(\boldsymbol{H})=\sum\limits_{j=0}^m a_j\boldsymbol{F}_j(\boldsymbol{H})
\end{equation}
are maximally convergent to $e^{\boldsymbol{H}}$ in $\mathcal{E}(d,c_f,a)$, i.e.,
\begin{equation}\label{eq_maximal_convergence}
\lim\limits_{m\rightarrow\infty}\sup \|e^{\boldsymbol{H}}-\boldsymbol{S}_m(\boldsymbol{H})\|_{\mathcal{E}(d,c_f,a)}^{1/m}=\lim\limits_{m\rightarrow\infty}\sup \|e^{\boldsymbol{H}}-\boldsymbol{p}_m^*(\boldsymbol{H})\|_{\mathcal{E}(d,c_f,a)}^{1/m},
\end{equation}
where $\boldsymbol{p}_m^*$ is the polynomial of degree $m$ that optimally approximates $e^{\boldsymbol{H}}$ in the infinite norm $\|\cdot\|_{\mathcal{E}(d,c_f,a)}$,

\begin{equation}
\label{eq_el_norm}
    \|\boldsymbol{h}(.)\|_{\mathcal{E}}=\max\limits_{x\in\mathcal{E}}|\boldsymbol{h}(x)|.
\end{equation}

In such a case, the polynomial coefficients $a_j$ are obtained from 
\begin{equation}\label{eq_intro_faber_coeff}
a_j=\int\limits_0^1{exp\left(\left(\gamma+\frac{c_f^2}{4\gamma}\right)\cos(2\pi\theta)+d+i\left(\gamma-\frac{c_f^2}{4\gamma}\right)\sin{2\pi\theta}\right)e^{-ij2\pi\theta}d\theta}.
\end{equation}

As pointed out in \citet{bergamaschi2003efficient}, if the ellipse where the polynomials are defined is large, it can be difficult to compute the coefficients \eqref{eq_intro_faber_coeff} accurately. When the ellipse capacity $\gamma$ increases, the magnitude of the term inside the integral \eqref{eq_intro_faber_coeff} grows several orders of magnitude, introducing significant numerical errors. Therefore, the calculation of the coefficients can require higher arithmetic precision in their calculations. Here we used double precision in all calculation, which was enough to ensure an adequate representation. 

Optimization of exponential time integrator methods, to reduce computational cost, can be achieved with the use of the scaling and squaring technique for large matrices \citep{higham2005scaling,moler2003nineteen}.
It consists on the selection of the parameters $s\geq1$, and $z\in\mathbb{N}$, such that the error of the truncated Faber series
\begin{equation}\label{intro_ss}
e^{\boldsymbol{H}}=\left(e^{s^{-1}\boldsymbol{H}}\right)^s\approx\left(\sum\limits_{j=0}^ma_j\boldsymbol{F}_j(s^{-1}\boldsymbol{H})\right)^s
\end{equation}
remains under a fixed threshold, while the amount of matrix-vector operations (MVOs) is minimized. The amount of MVOs required by Eq.\,\eqref{intro_ss} is $s\times m$, where $m$ is the degree of the polynomial used in the approximation of $e^{s^{-1}\boldsymbol{H}}$. This approach was successfully implemented for Padé approximations of the exponential in the works of \citet{higham2005scaling}, and \citet{al2010new}, connecting it later to exponential integrators in \citet{al2011computing}. However, to estimate the optimal $s$ and $m$, reliable error bounds must be obtained, as is shown in works of \citet{higham2005scaling}, and \citet{al2010new}, where sharp bounds for the Padé approximation were required. For Faber polynomials, this is still an open problem, and the purpose of the next subsection is to provide further insights in this direction.

\subsubsection{Faber exponential error bounds}

A general expression for the error (on the usual Euclidean norm for $\mathbb{C}^n$, $\|.\|_2$) is stated in \citet{bergamaschi2000efficient},
\begin{equation}\label{eq_error_1}
\|e^{\boldsymbol{H}}-\boldsymbol{S}_{m}(\boldsymbol{H})\|_2\leq \text{cond}_2(\boldsymbol{P})\|\text{exp}(.)-\boldsymbol{S}_{m}(.)\|_{\mathcal{E}},
\end{equation}
where $\boldsymbol{S}_{m}$ is the partial sum of the first $m+1$ terms of the series in the recurrence \eqref{eq_faber_pol_0}-\eqref{eq_faber_pol}, $\boldsymbol{P}$ is the diagonalization  matrix of $\boldsymbol{H}$ (i.e., $\boldsymbol{P}^{-1}\boldsymbol{H}\boldsymbol{P}$ is a diagonal matrix), $\text{cond}_2(\boldsymbol{P})=\|\boldsymbol{P}\|_2\|\boldsymbol{P}^{-1}\|_2$, and $\|.\|_{\mathcal{E}}$ is the infinite norm over the ellipse $\mathcal{E}$ (see \eqref{eq_el_norm}), enclosing the spectrum of $\boldsymbol{H}$.

Expression \eqref{eq_error_1} has the complication that determining $\text{cond}_2(\boldsymbol{P})$ is a difficult challenge for large matrices, as the ones resulting after discretizing spatial derivatives in partial differential equations.
In addition, to the authors' best knowledge, a general bound for the term $\|\text{exp}(.)-\boldsymbol{S}_{m}(.)\|_{\mathcal{E}}$  has not been reported in the literature, yet.

In \citet{bergamaschi2000efficient}, they present a bound for the second term on the right of \eqref{eq_error_1},
\begin{equation}\label{eq_error_2}
\|\exp(.)-\boldsymbol{S}_m(.)\|_{\mathcal{E}}\leq\left\{\begin{array}{ll}\frac{8\gamma}{m}\exp\left(\frac{4\gamma^2}{4\gamma-m}+d-\frac{m^2}{4\gamma}+\frac{c^2(4\gamma-m)}{16\gamma^2}\right), & m\leq 2\gamma,\\
4 \exp\left(d+\frac{c^2}{4m}\right)\left(\frac{e\gamma}{m}\right)^m, & m>2\gamma,
\end{array}
\right.
\end{equation}
where $\gamma$ is the mean of the semi-axis of the ellipse $\mathcal{E}$, $d$ is its center, $c$ is the eccentricity, $e$ is the base of the natural logarithm, and $m$ is the polynomial degree used. The second part of \eqref{eq_error_2}, when $m>2\gamma$, represents an asymptotic super-linear convergence. But when $m\leq 2\gamma$, the expression in (\ref{eq_error_2}) is very pessimistic and is not appropriate for practical applications (see \citet{bergamaschi2000efficient}). Furthermore, although it was not stated in their research, the inequality \eqref{eq_error_2} is only guaranteed if the ellipse is strictly contained in the right half plane (see \citet{moret2001computation}). However, for many problems of interest, including the wave problems to be discussed in this work, there are eigenvalues with negative real part and, therefore, expression \eqref{eq_error_2} is not fulfilled.

Here, we propose an error bound that improves \eqref{eq_error_2} when the ellipses are in the positive half plane, and is also valid for any kind of ellipses. Let $m_{\epsilon/2}$ be the polynomial degree such that the Faber approximation error is at most $\epsilon/2$. Then, for all $m<m_{\epsilon/2}$ we have,

\begin{flalign}
\|\exp(.)-\boldsymbol{S}_{m}(.)\|_{\mathcal{E}}&\leq \|S_{m_{\epsilon/2}}(.)-\boldsymbol{S}_{m}(.)\|_{\mathcal{E}}+\|\exp(.)-\boldsymbol{S}_{m_{\epsilon/2}}(.)\|_{\mathcal{E}}\nonumber\\
&=\left\lVert\sum\limits_{j={m+1}}^{m_{\epsilon/2}}a_j\boldsymbol{F}_j(.)\right\rVert_\mathcal{E}+\|\exp(.)-\boldsymbol{S}_{m_{\epsilon/2}}(.)\|_{\mathcal{E}}.\nonumber\\
&\leq \sum\limits_{j=m+1}^{m_{\epsilon/2}}|a_j|\|\boldsymbol{F}_j(.)\|_\mathcal{E}+\|\exp(.)-\boldsymbol{S}_{m_{\epsilon/2}}(.)\|_{\mathcal{E}}\nonumber\\
&\leq \sum\limits_{j=m+1}^{m_{\epsilon/2}}|a_j|\|\boldsymbol{F}_j(.)\|_\mathcal{E}+\frac{\epsilon}{2}\label{eq_error_bound_FA}.
\end{flalign}

The polynomial coefficients $a_j$ can be calculated by explicit formulas using \eqref{eq_intro_faber_coeff}. Thus, the unknowns in \eqref{eq_error_bound_FA} are the norm over the ellipse $\mathcal{E}$ of the polynomials $\boldsymbol{F}_j$, and $m_{\epsilon/2}$ discussed next.

From Section \ref{sec_faber}, we know that $\boldsymbol{F}_j$ are stretched and translated Chebyshev polynomials, then, all its roots are located at the points
\begin{equation*}
    r_k=d+c_f\cos\left(\frac{1+2k}{2j}\pi\right),\quad k=0,...,j-1.
\end{equation*}
Therefore, as stated in \citet{munch2019chebyshev}, the extreme values of the polynomials are attained on points over the ellipse $\mathcal{E}$. Moreover, one of those points are where the ellipse cuts the line passing by its focuses. Hence, we can compute $\|\boldsymbol{F}_j(.)\|_\mathcal{E}$ in a straightforward way by evaluating the polynomial at these two points.

To estimate $m_{\epsilon/2}$ we may calculate the series in \eqref{eq_error_bound_FA} until the new terms are smaller than $\epsilon/2$. In practice, due to the extremely fast convergence to zero of Faber coefficients $a_j$, the rest of the series $\|\exp(.)-\boldsymbol{S}_{m_{\epsilon/2}}(.)\|_{\mathcal{E}}$ will be negligible when compared with the terms already computed.

We compare bounds \eqref{eq_error_2} and \eqref{eq_error_bound_FA}, by means of numerical experiments, using arbitrary diagonalizable matrices $\boldsymbol{H}=\boldsymbol{P}\boldsymbol{D}\boldsymbol{P}^{-1}$. To ensure that the condition number of the diagonalization matrix $\boldsymbol{P}$ is $1$, we take real normal matrices $\boldsymbol{H}$, since in this case $\boldsymbol{P}$ is unitary and therefore the condition number of $\boldsymbol{P}$ in the Euclidean norm is 1, i.e., $\text{cond}_2(\boldsymbol{P})=1$.

\begin{figure}[hbt]
	\subfloat[Ellipse of minimum $\gamma$ containing the \linebreak $n=60$ eigenvalues of $H$.]{\includegraphics[scale=0.35]{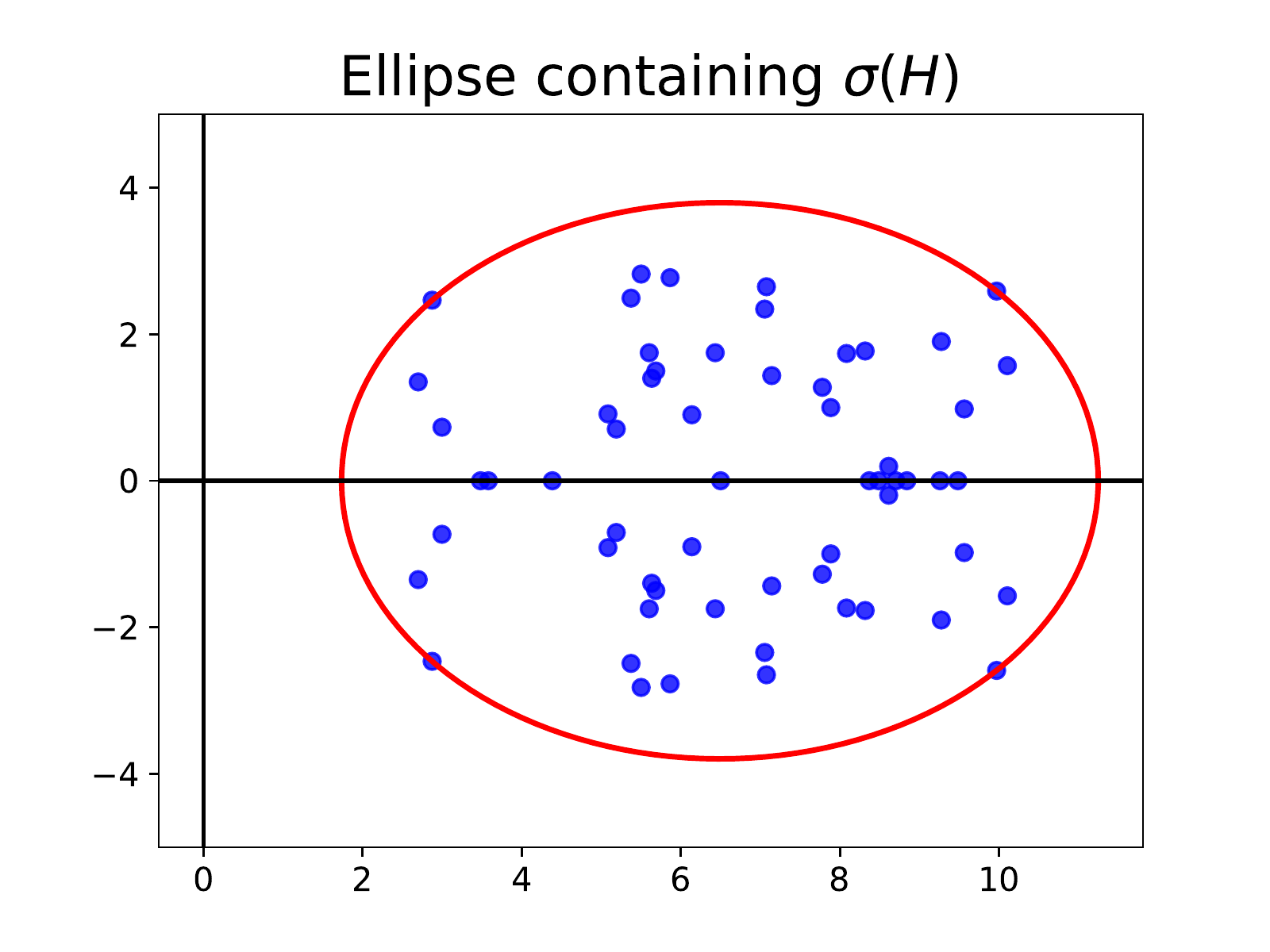}}
	\hfill
	\subfloat[Error bounds for the exponential of matrix $H$.]{\includegraphics[scale=0.35]{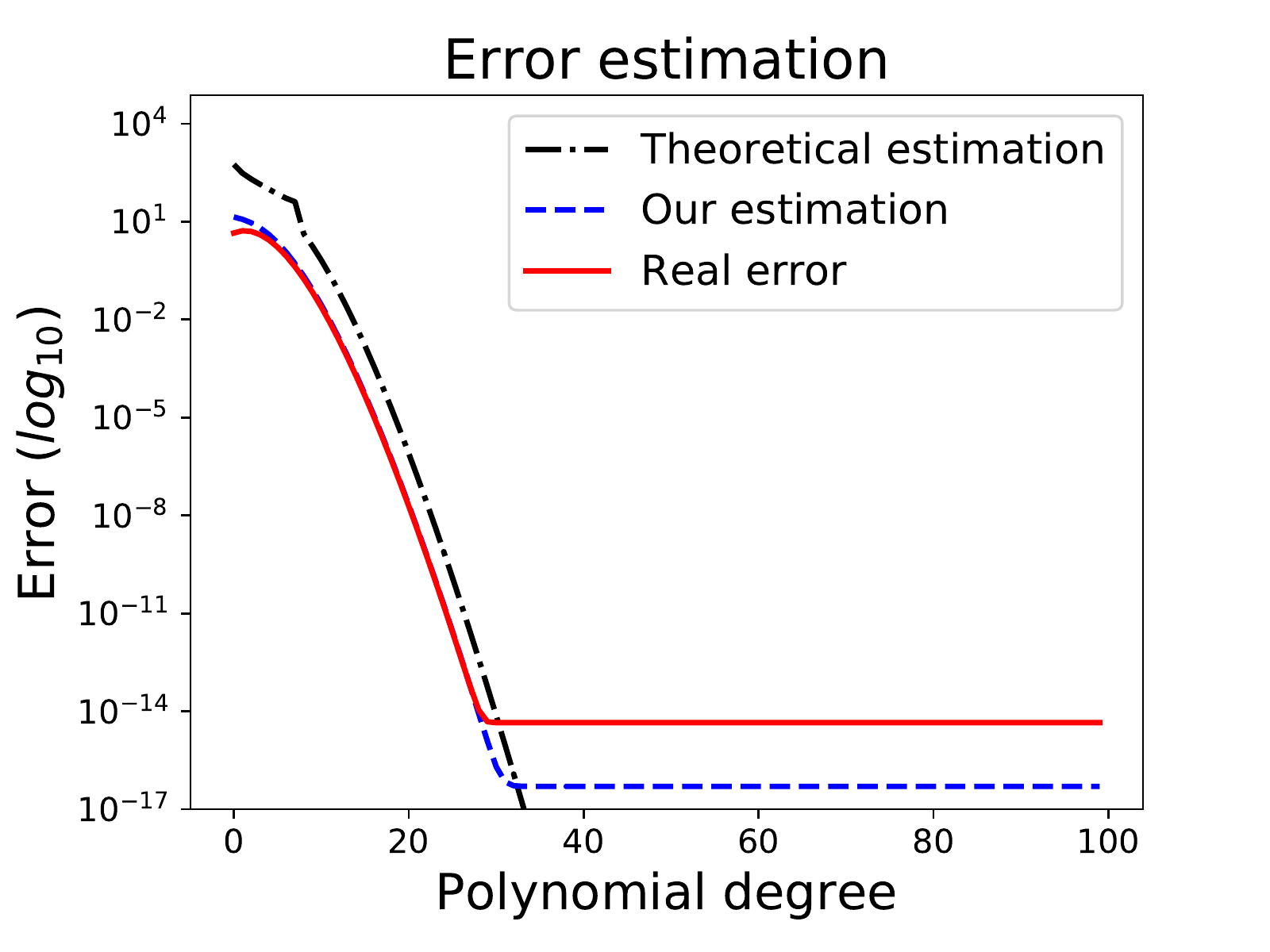}}\\
	\subfloat[Ellipse of minimum $\gamma$ containing the \linebreak $n=60$ eigenvalues of $H$.]{\includegraphics[scale=0.35]{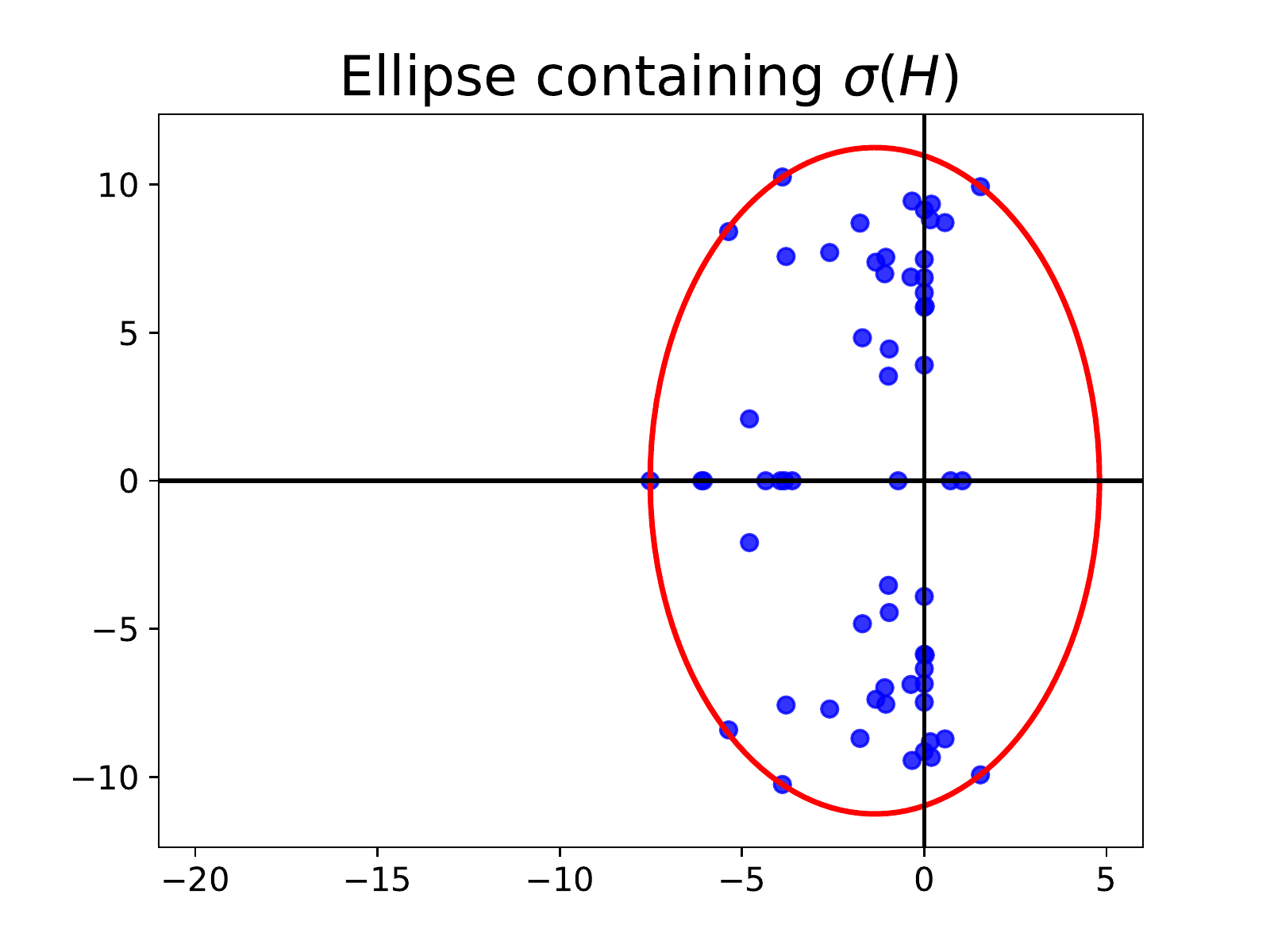}}
	\hfill
	\subfloat[Error bounds for the exponential of matrix $H$.]{\includegraphics[scale=0.35]{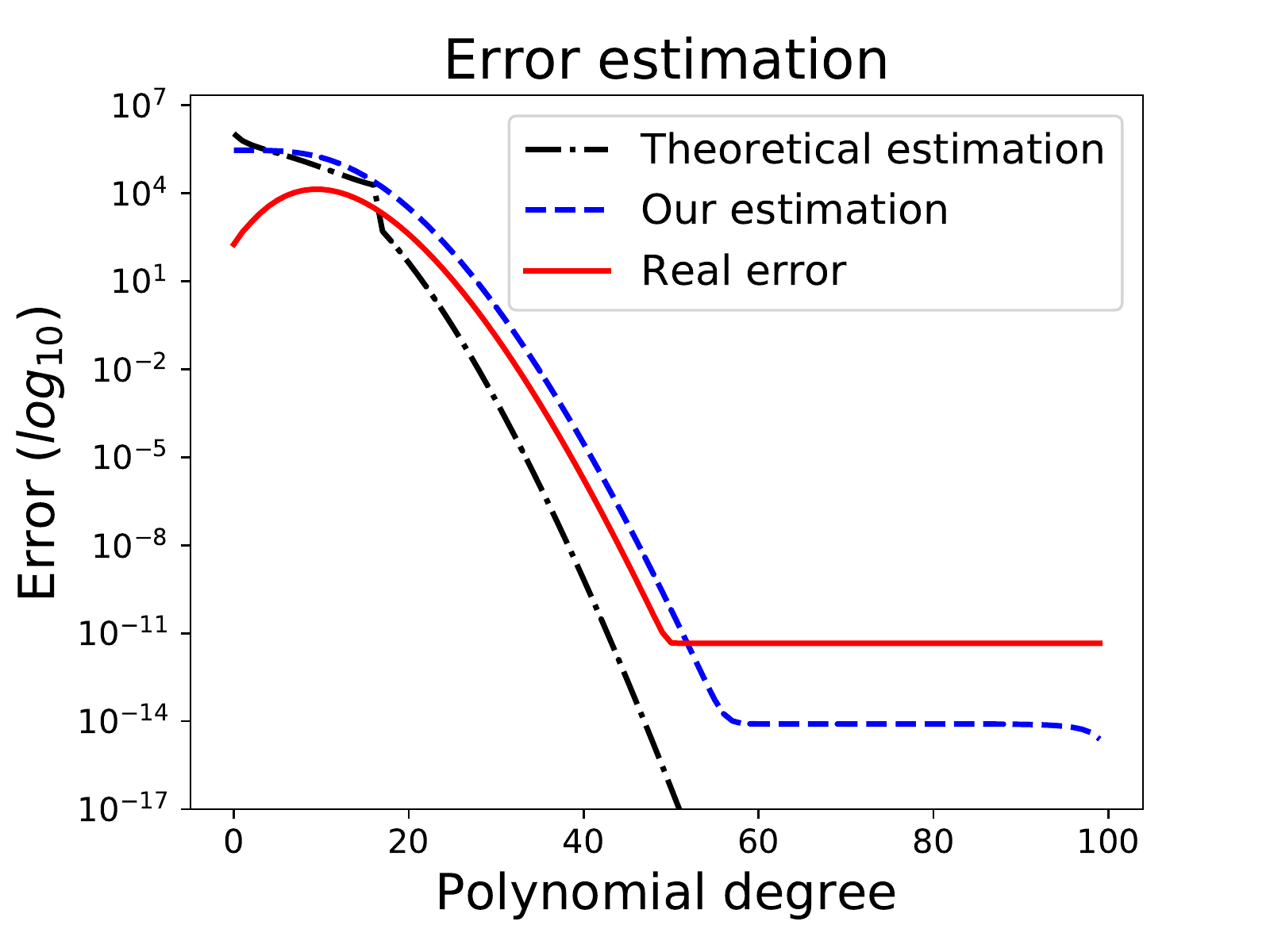}}
	\caption{(Upper images) Error estimation of $e^{\boldsymbol{H}}$ for a real normal matrix $\boldsymbol{H}$ with randomly generated eigenvalues on the right half plane. (Lower images) Same as before, but on both half planes. (Right images) The theoretical estimation is based on \citet{bergamaschi2000efficient}, our approximation is based on inequality \eqref{eq_error_bound_FA}, and the real error are the effectively calculated errors, subject to rounding errors effects, when approximating the exponential. The bound proposed in this research is sharper than the one reported in the literature.}\label{fig_error_bound1}
\end{figure}

In the following illustrative examples, we set $\boldsymbol{P}$ as a random orthonormal matrices (therefore unitary)\footnote{We used the python function \texttt{ortho\_group.rvs} from the \texttt{Python}'s package \texttt{scipy}} with dimensions $60 \times 60$. Then we define the eigenvalues of the diagonal matrices $\boldsymbol{D}$ as the composition of 10 randomly generated real numbers, and 50 randomly sampled complex eigenvalues, symmetric respect to the real axis. We present results of two experiments, for which the ranges of the randomly generated eigenvalues are given by:
\begin{enumerate}
    \item Experiment 1 (Figures \ref{fig_error_bound1}a and \ref{fig_error_bound1}b): real numbers between $[2.8,\;10.8]$ and complex values in the domain $[2.8,\;10.8]\times i[-3,\;3]$.
    \item Experiment 2 (Figures \ref{fig_error_bound1}c and \ref{fig_error_bound1}d): real numbers between $[-8,\;2]$ and complex values in the domain $[-8,\;2]\times i[-11,\;11]$.
\end{enumerate}

The spectrum and the comparison between the errors bound are shown in Figure \ref{fig_error_bound1}.

From Figure \ref{fig_error_bound1}, we can observe that the error bound in Eq.\,\eqref{eq_error_bound_FA} is sharper than in Eq.\,\eqref{eq_error_2}.
In fact, it can be even used as an estimator of the error of high polynomial degrees.
This behavior is also observed in several other simulations for normal matrices generated using a random number of eigenvalues with a uniform distribution. Moreover, for ellipses not contained in the right half plane, the error bound \eqref{eq_error_2} is not particularly reliable (see Fig.\,\ref{fig_error_bound1}d).

Although these results provide valuable theoretical insight into Faber polynomials for normal matrices, in wave propagation cases with absorbing boundary conditions, which will be investigated later on, the matrix operator is not necessarily normal and, thus, the bound \eqref{eq_error_1} cannot be readily applied.

\subsubsection{Faber's polynomials in conics}\label{sec_conics}

The performance of Faber polynomials depends on the conics they are defined on. This dependence is very strong, in the sense that different ellipses can lead to an enhancement or deterioration of the error in several orders of magnitude. Therefore, we dedicate this Section to discuss the approximation performance for different types of conics.

\begin{figure}[hbt]
	\subfloat[Five ellipses, one optimal (green), four others varying in sizes, and the eigenvalues of a randomly generated normal matrix.]{\includegraphics[scale=0.36]{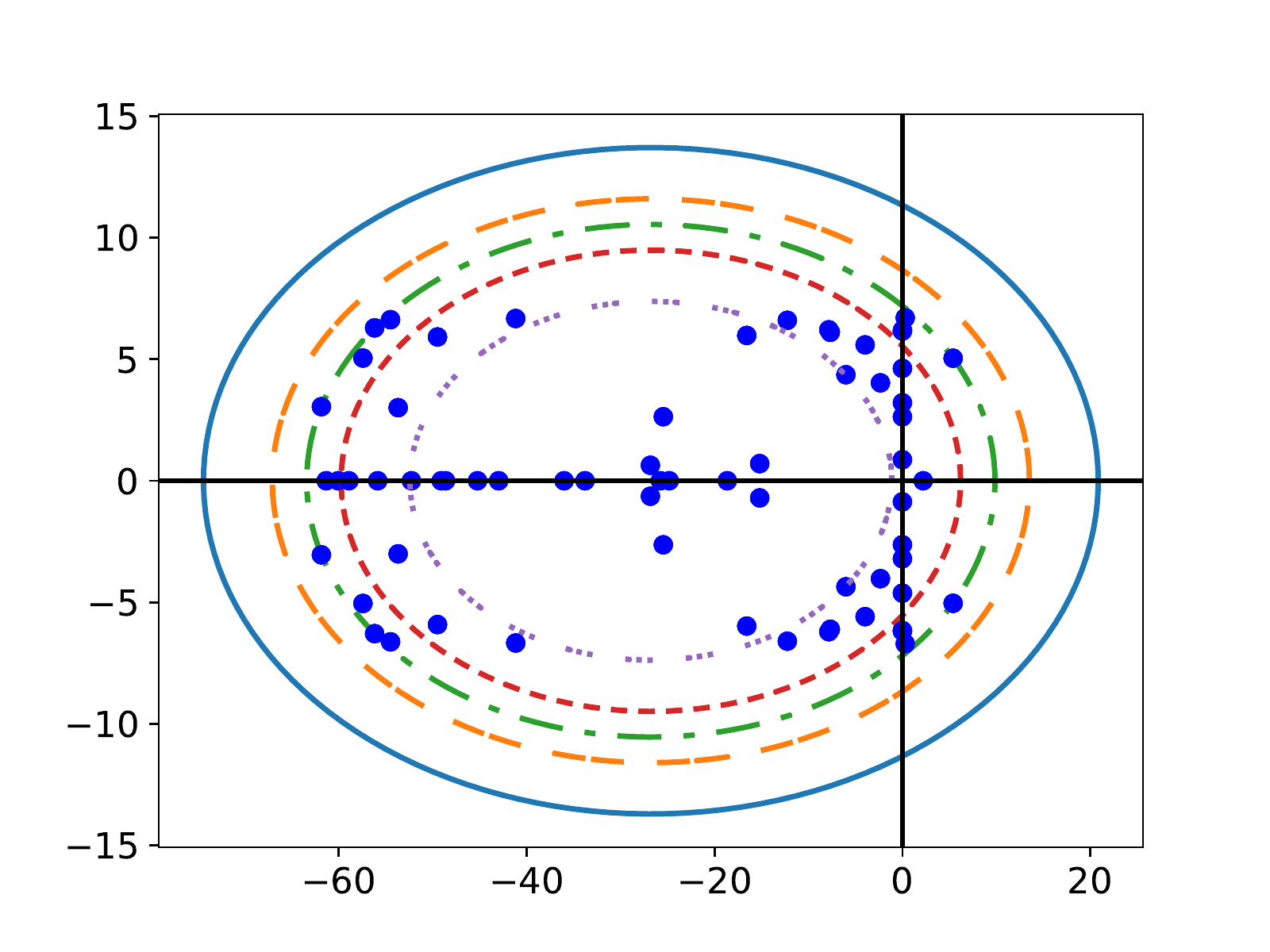}}
	\hfill
	\subfloat[Error of Faber approximation using each of the ellipses on the left to approximate the exponential of the matrix.]{\includegraphics[scale=0.36]{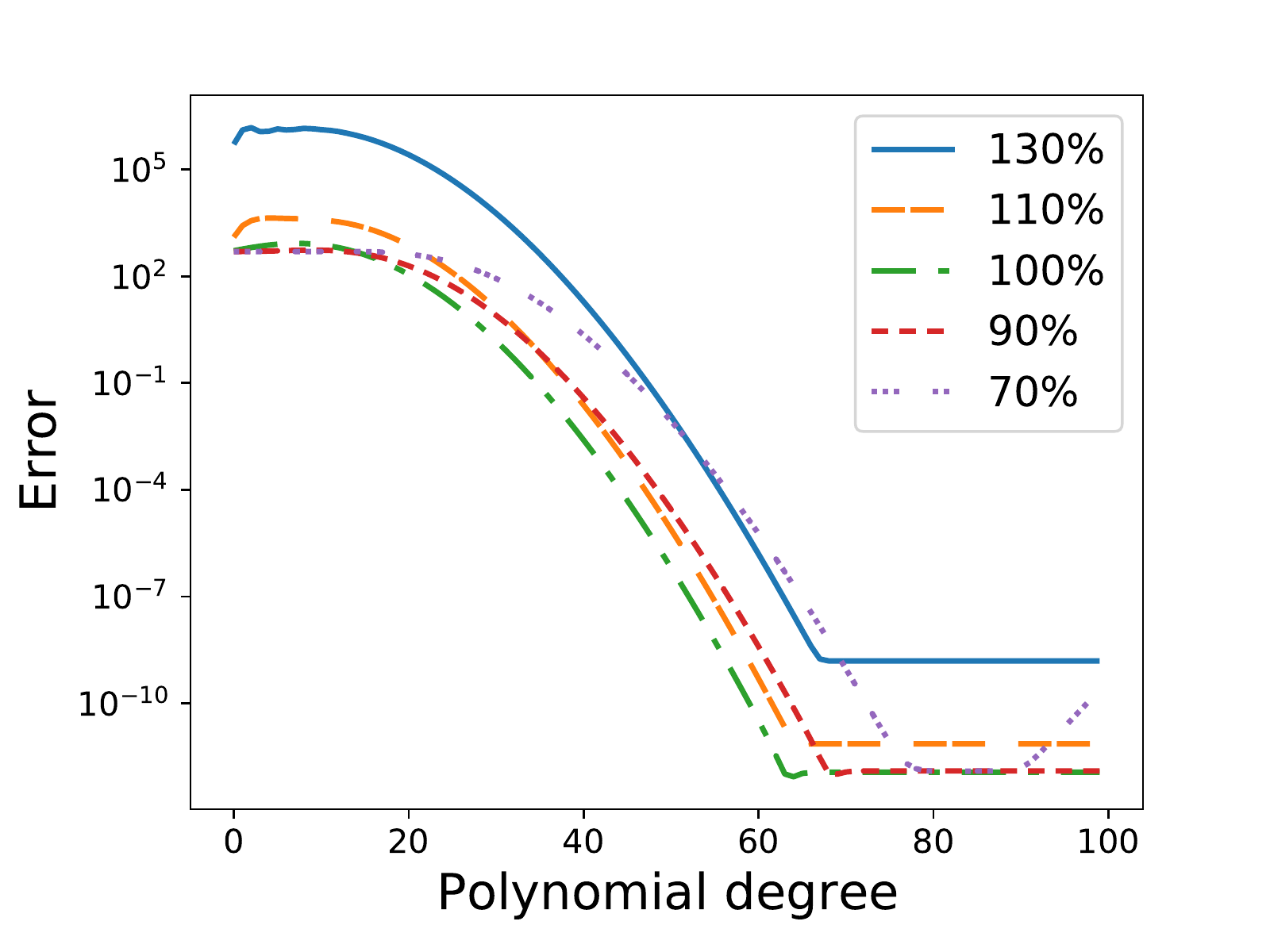}}
	\caption{Error of Faber polynomials on the approximation of the exponential of a randomly generated normal matrix by using the optimal ellipse (green), and four other ellipses with varying sizes in a percentage of the optimal ellipse. If the ellipse does not enclose all eigenvalues or if the ellipse is too large, errors can be considerably larger.}\label{fig_five_ellipses}
\end{figure}

The requirement in the Faber series expansions of matrix functions of having $\sigma(\boldsymbol{H})$ enclosed by the ellipse $\mathcal{E}(d,c_f,a)$ is not only a necessary condition from the theory.
Violations of this requirement in numerical experiments show considerably worsened solutions.
Furthermore, if we choose an ellipse larger than necessary, the convergence series slows down. Fig.\,\ref{fig_five_ellipses} depicts the error in an exponential approximation with respect to the polynomial degree considering ellipses of different sizes. 

As we notice, the best approximation is given by the smallest ellipse that still covers all eigenvalues of the operator. For the other cases, the error increases for the lower polynomial degrees at least of an order of magnitude. Hence, a good estimation of a small-as-possible convex cover of $\sigma(\boldsymbol{H})$ to construct the ellipse seems to be of utmost importance when dealing with computational efficiency.

From bound \eqref{eq_error_bound_FA}, we notice that the amplitude of Faber coefficients $a_j$ influence the convergence speed. When they are smaller, the error bound \eqref{eq_error_bound_FA} is lower and, therefore, a faster convergence may be achieved. The principal constants influencing the magnitude of $a_j$ coefficients are the ellipse capacity $\gamma$, and the ellipse eccentricity $c_f$. In Fig.\,\ref{fig_five_ellipses}, we observe that when the capacity of the ellipse increases, the velocity of the convergence decreases.

Next, we also investigate decreasing the ellipse eccentricity $c_f$ at the cost of increasing the ellipse capacity. Starting from the ellipse with minimum capacity, we construct the optimal ellipses with a predefined $c_f$, followed by reducing step-by-step $c_f$ with each ellipse until we get a circle ($c_f=0$).

\begin{figure}[H]
    \subfloat[Five ellipses with decreasing $c_f$ and the eigenvalues of a randomly generated normal matrix.]{\includegraphics[scale=0.36]{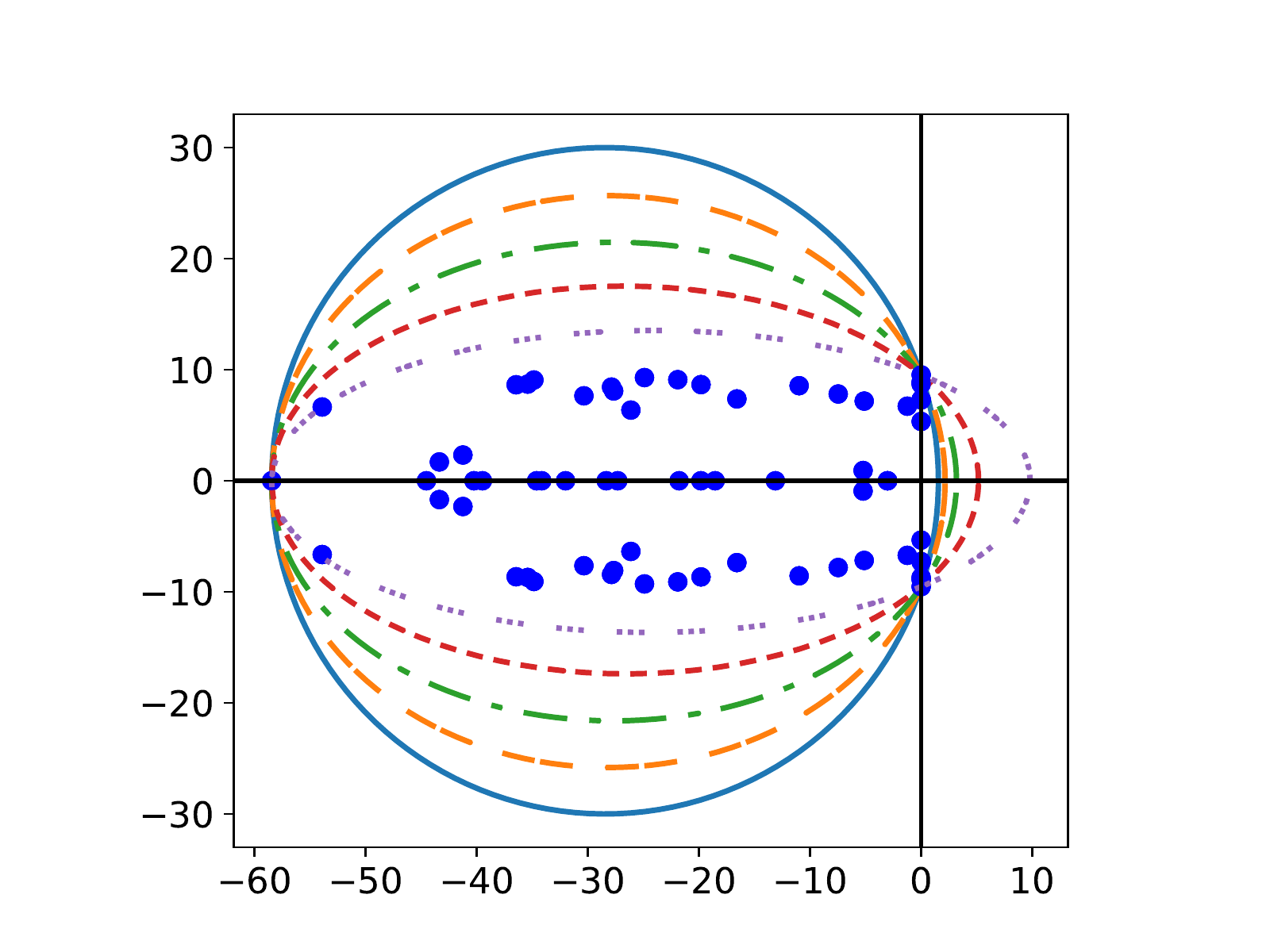}}
    \hfill
    \subfloat[Error of Faber approximation for each of the ellipses on the left.]{\includegraphics[scale=0.36]{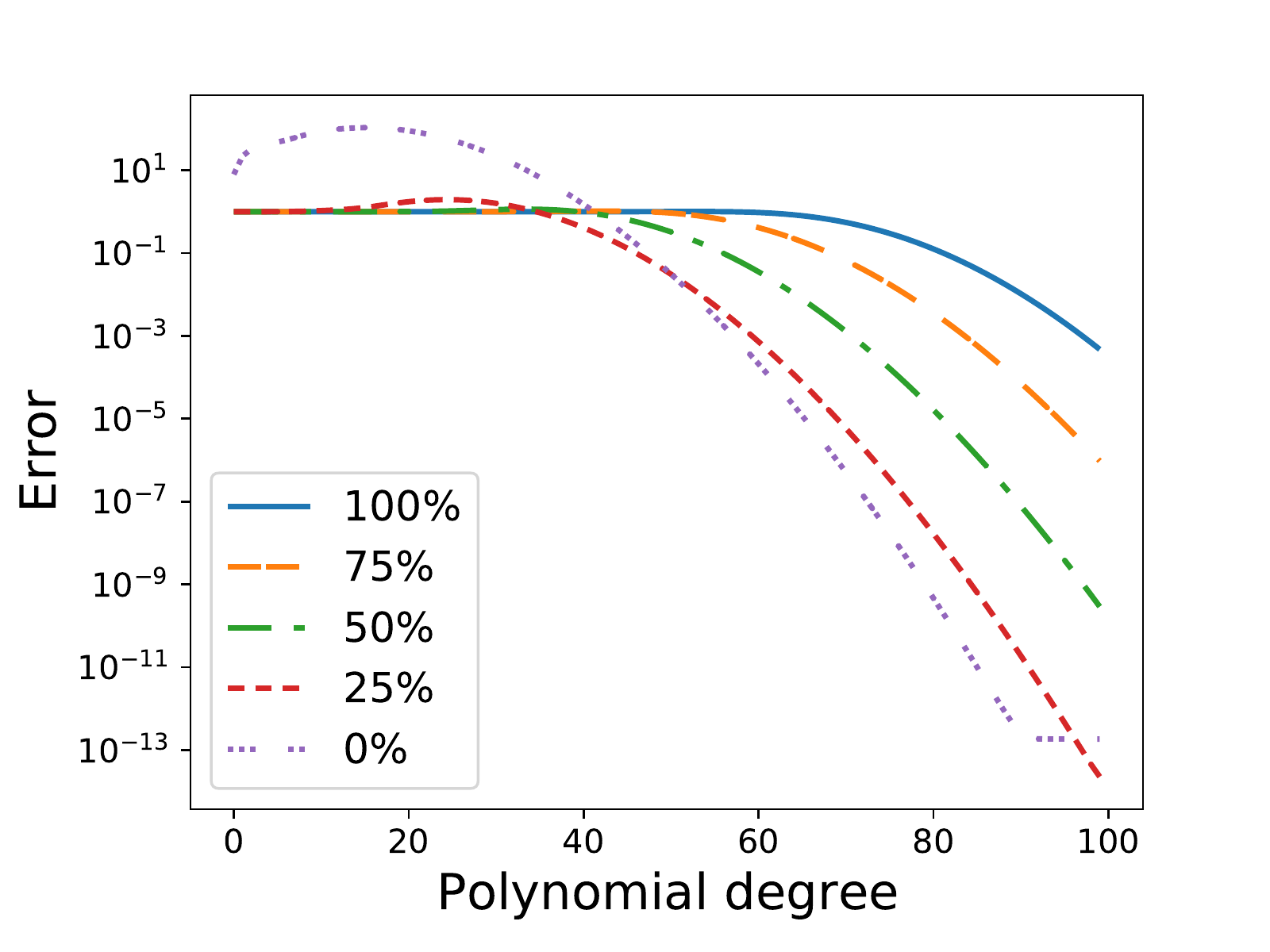}}
\caption{Error of Faber polynomials on the approximation of the exponential of a normal matrix by using five different ellipses with decreasing eccentricity. When the eccentricity diminishes, the error is lower for some degrees and higher for others.}\label{fig_five_ellipses_1}
\end{figure}

Fig.\,\ref{fig_five_ellipses_1} shows the errors of calculating the matrix exponential for each of the conics using the matrix infinity norm.
For low polynomial degrees, the approximations using conics with less eccentricity is better, but the magnitude of the errors are too high to be considered as a good approximation. In addition, the ellipse with minimum capacity has the lowest error for high degrees.

\begin{figure}[H]
    \subfloat[Ellipse and circle used in the Faber exponential approximation.]{\includegraphics[scale=0.36]{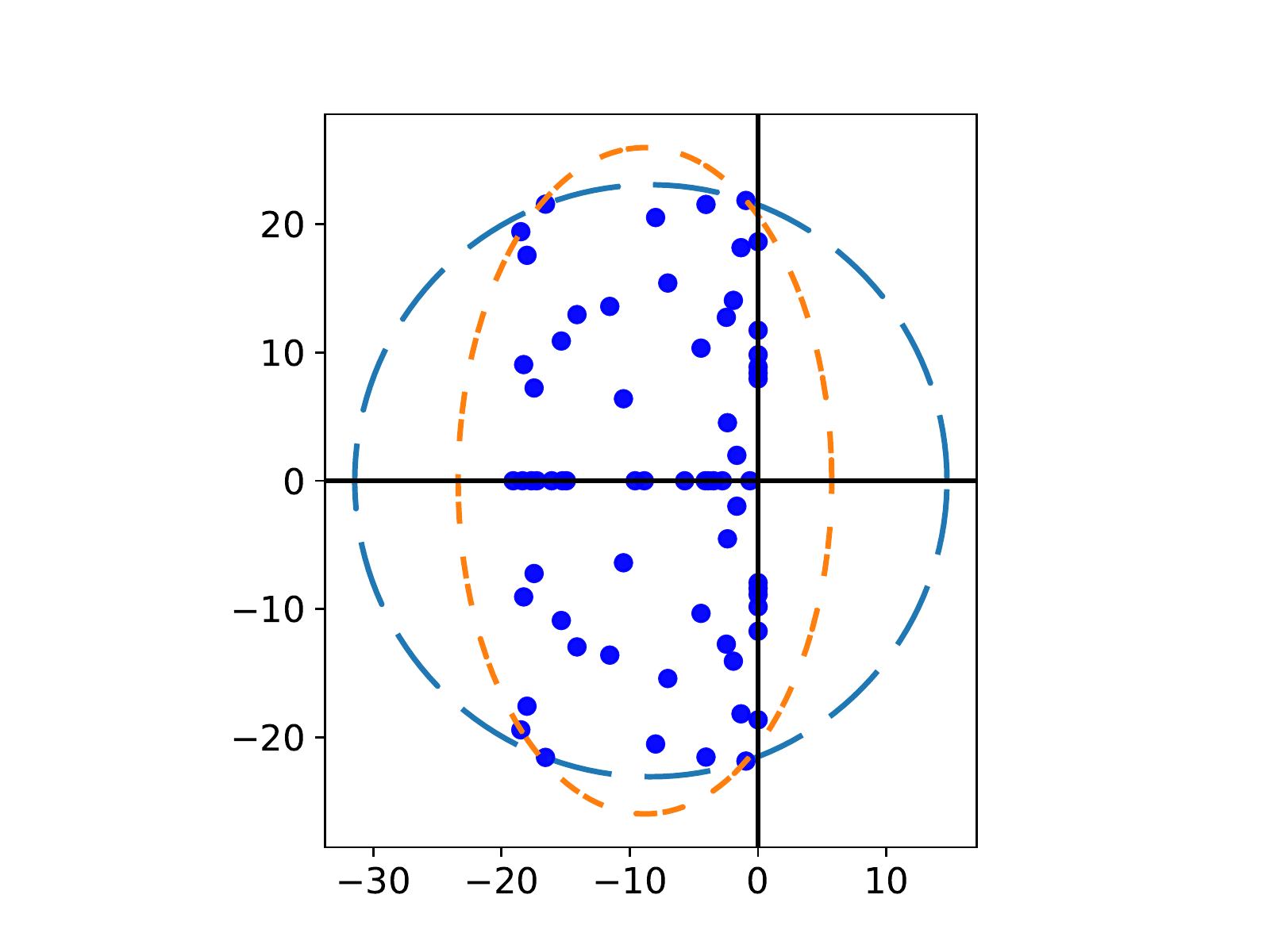}}
    \hfill
    \subfloat[Error of Faber approximation for ellipse and circle shown on the left.]{\includegraphics[scale=0.36]{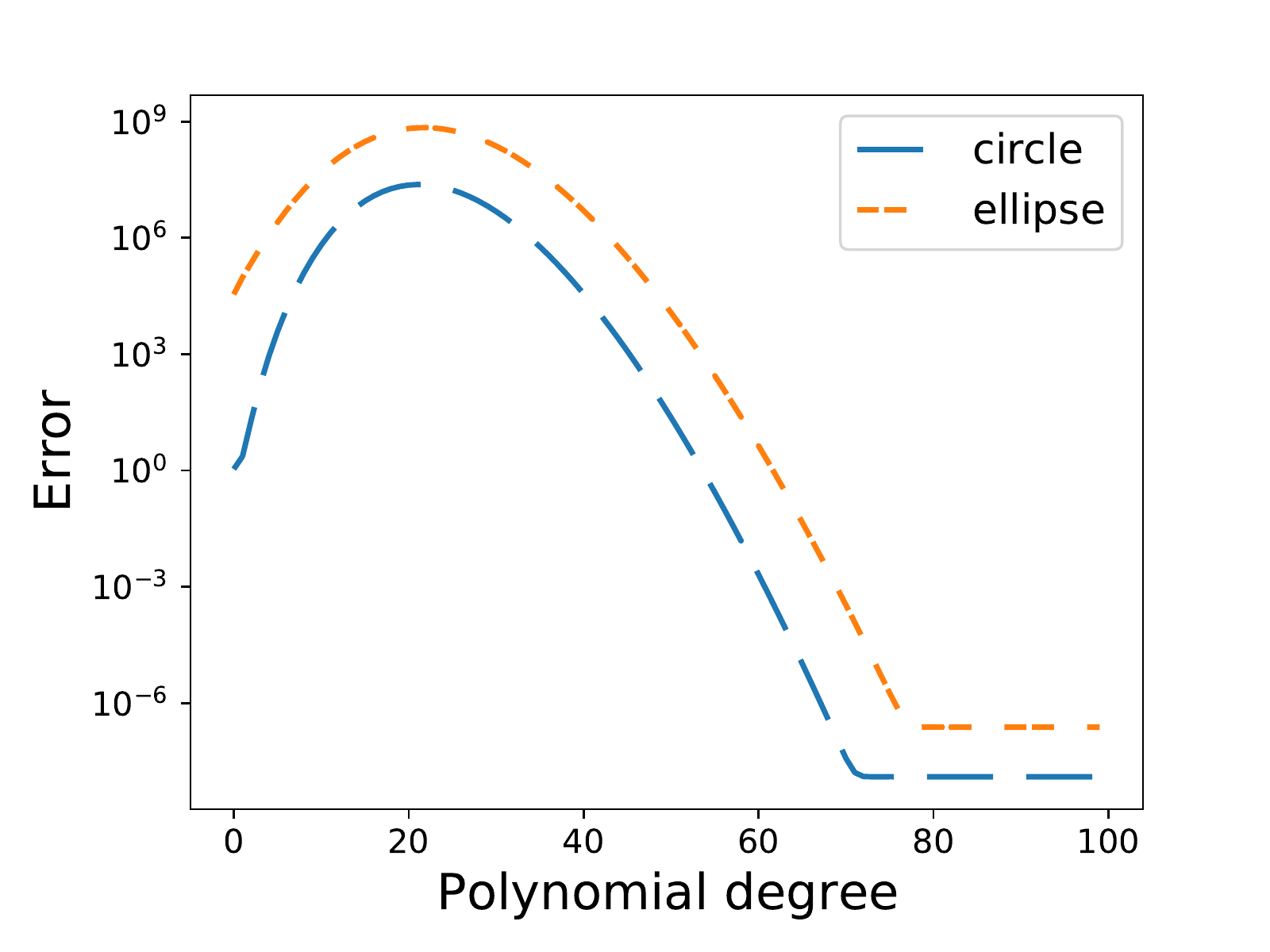}}
\caption{Error of Faber polynomials on the approximation of the exponential of a normal matrix by using the optimal ellipse and an optimal circle. In this particular case, a circle, with higher capacity, outperforms the use of the optimal ellipse.}\label{fig_five_ellipses_2}
\end{figure}

In general, the conclusion that an ellipse with minimum capacity has the lowest error for high degrees is not always valid.
If we have a different distribution of eigenvalues, similar to the ones appearing in the wave equations discussed in the next sections, we may get different relations as illustrated in Fig.\,\ref{fig_five_ellipses_2}. As such, a circle can provide a better alternative than an ellipse with minimum capacity under special circumstances. Moreover, maybe even other conics could be more suited to allow faster convergence, depending on the eigenvalue distribution. Further investigations on the Faber convergence with respect to eccentricity and conic form may be of interest, but goes beyond the scope of this paper and will be discussed elsewhere.
Hence, through this paper, we will study only the method over an ellipse with minimum capacity.

\section{Application to seismic waves}\label{sec_application}

This section is dedicated to the mathematical formulations of the wave equations with an absorbing boundary condition, its numerical discretization, and the description of the experiments used in this work. Taking all these different aspects into account is of utmost importance for the Faber approximation, since they define the discrete operator and the exponential method relies on the spectrum of the operator.

\subsection{Formulations of the wave equations with PML}

The eigenvalue distribution of a discrete operator is strongly influenced by the continuous equations formulations. In addition, the wave equations with PML can be expressed with different formulations as a first order system in time.
Therefore, to study the spectrum $\sigma(\boldsymbol{H})$ of the matrix operator $\boldsymbol{H}$, and the characteristics of the approximation using Faber polynomials, we take different formulations into account to obtain better numerical conclusions.
All formulations include an absorbing boundary layer using the Perfectly Matching Layer (PML) as discussed next.

The PML is one of the most popular absorbing boundary conditions used for the wave equations and related areas.
While termed as a boundary condition, it is in fact an extension of the problem to a larger domain, containing an absorbing layer, together with a set of additional variables and equations acting on this layer. It is very effective for most seismic imaging applications, but it is also of relative complex implementation, due to an increase in the number of equations. It is build from a transformation of the real domain to the complex plane, where the waves outside the physical region of interest (PML layers) are attenuated, while the others (inside the physical domain) remains unchanged \citep{berenger1994perfectly,assi2017compact}.

For the acoustic wave equations, we have two formulations in one and two dimensions, where the main difference is about the spatial derivatives to be of second (2SD) or first order (1SD), given as follows.\\

\noindent
\textbf{One dimensional form:} ($x \in \Omega=[a_1,a_2]$, $t>t_0$)
\begin{align}
	    \frac{\partial }{\partial t}
	    \begin{pmatrix}
	    u\\v\\w
	    \end{pmatrix}
	        &=\begin{pmatrix}
	        0&1&0\\
	        c^2\frac{\partial ^2}{\partial x^2}&-\beta_x&c^2\frac{\partial }{\partial x}\\
	        -\beta_x\frac{\partial}{\partial x}&0&-\beta_x
	        \end{pmatrix}
	        \begin{pmatrix}
	    u\\v\\w
	    \end{pmatrix}\begin{pmatrix}
	    0\\f\\0
	    \end{pmatrix},\qquad \text{(2SD)}&&\label{eq_2SD_1D}\\
	    &\nonumber\\
	    \frac{\partial }{\partial t}
	    \begin{pmatrix}
	    u\\v\\w
	    \end{pmatrix}
	        &=\begin{pmatrix}
	        0&c^2\frac{\partial}{\partial x}&-c^2\\
	        \frac{\partial}{\partial x}&-\beta_x&0\\
	        0&\beta_x\frac{\partial}{\partial x}&-\beta_x
	        \end{pmatrix}
	        \begin{pmatrix}
	    u\\v\\w
	    \end{pmatrix}+\begin{pmatrix}
	    \int f\\0\\0
	    \end{pmatrix}.\qquad \text{(1SD)}\label{eq_1SD_1D}
\end{align}

\noindent
\textbf{Two-dimensional form:}  ($(x,y) \in \Omega=[a_1,a_2]\times[b_1,b_2]$, $t>t_0$)

   \begin{equation}
           \frac{\partial}{\partial t}\begin{pmatrix}
    u\\v\\w_x\\w_y
    \end{pmatrix}
        =\begin{pmatrix}
        0&1&0&0\\
        -\beta_x\beta_y+c^2\left(\frac{\partial^2}{\partial x^2}+\frac{\partial^2}{\partial y^2}\right) & -(\beta_x+\beta_y)&c^2\frac{\partial}{\partial x}&c^2\frac{\partial}{\partial y}\\
        (\beta_y-\beta_x)\frac{\partial}{\partial x}& 0& -\beta_x&0\\
        (\beta_x-\beta_y)\frac{\partial}{\partial y}& 0& 0&-\beta_y
        \end{pmatrix}
        \begin{pmatrix}
    u\\v\\w_x\\w_y
    \end{pmatrix}
 	+\begin{pmatrix}
    0\\f\\0\\0
    \end{pmatrix}, \text{(2SD)}
   \end{equation}

\begin{flalign}
	\frac{\partial }{\partial t}
	\begin{pmatrix}
		u\\v_x\\v_y\\w_x\\w_y
	\end{pmatrix}
	&=\begin{pmatrix}
		0&c^2\frac{\partial}{\partial x}&c^2\frac{\partial}{\partial y}&-c^2&-c^2\\
		\frac{\partial}{\partial x}&-\beta_x&0&0&0\\
		\frac{\partial }{\partial y}&0&-\beta_y&0&0\\
		0&\beta_x\frac{\partial}{\partial x}&0&-\beta_x&0\\
		0&0&\beta_y\frac{\partial}{\partial y}&0&-\beta_y
	\end{pmatrix}
	\begin{pmatrix}
		u\\v_x\\v_y\\w_x\\w_y
	\end{pmatrix}+\begin{pmatrix}
		\int f\\0\\0\\0\\0
	\end{pmatrix}.\qquad \text{(1SD)}
\end{flalign}

Here, $u=u(t,x)$ (or $u=u(t,x,y)$ in 2D) is the displacement, $c=c(x)$ (or $c=c(x,y)$ in 2D) is the given velocity distribution in the medium, $v=v(t,x)$ (or $(v_x,v_y)=(v_x(t,x,y),v_y(t,x,y))$ in 2D) is the wave velocity, $f=f(x,t)$ (or $f=f(x,y,t)$ in 2D) is the source term, and its time integral $\int f$ is calculated over the interval $[t_0,t]$, where $t_0$ is the initial time. The $w$-functions \linebreak ($w=w(t,x)$ in 1D and $(w_x,w_y)=(w_x(t,x,y),w_y(t,x,y))$ in 2D) are the auxiliary variables of the PML approach and the $\beta$-functions are known and control the damping factor in the absorbing layer. 

The spatial domain is decomposed into two parts: a main physical domain of interest and an outer domain layer surrounding the physical one, used to place the wave absorbing conditions (PML).
On the physical domain (outside the PML layer) the $\beta$-functions are zero and the system of equations coincides with the classic wave propagation without absorbing boundary conditions. Therefore, the auxiliary $w$-functions are different from zero only in the PML domain. On the other hand, since the displacement is attenuated in the PML layer, arbitrary conditions can be set at the boundary of the PML region, where zero-Dirichlet conditions are adopted. 
The differential equations are then well-defined once the initial conditions are given in conjunction with the Dirichlet (null displacement) outer boundary conditions, see \citet{assi2017compact}.

For the sake of readability, equations and details for the two-dimensional elastic wave problem are only provided in Appendix \ref{sec_appendix_continuous_framework}, together with further description of the continuous equations for the acoustic problem.

\subsection{Numerical discretization by finite differences}\label{sec_numerical_framework}

\begin{figure}[H]
\centering
\begin{tikzpicture}
\draw[blue, very thick] (0,0) -- (0,3) -- (6,3) -- (6,0) -- (0,0);
\node at (6.8,2.4) {\Huge\textcolor{blue}{$\Omega$}};
\filldraw[color=red, fill=red!10, very thick](0,0) circle (0.07);
\filldraw[color=red, fill=red!10, very thick](1,0) circle (0.07);
\filldraw[color=red, fill=red!10, very thick](2,0) circle (0.07);
\filldraw[color=red, fill=red!10, very thick](3,0) circle (0.07);
\filldraw[color=red, fill=red!10, very thick](4,0) circle (0.07);
\filldraw[color=red, fill=red!10, very thick](5,0) circle (0.07);
\filldraw[color=red, fill=red!10, very thick](6,0) circle (0.07);
\filldraw[color=black, fill=black!80, very thick](0.5,0) circle (0.07);
\filldraw[color=black, fill=black!80, very thick](1.5,0) circle (0.07);
\filldraw[color=black, fill=black!80, very thick](2.5,0) circle (0.07);
\filldraw[color=black, fill=black!80, very thick](3.5,0) circle (0.07);
\filldraw[color=black, fill=black!80, very thick](4.5,0) circle (0.07);
\filldraw[color=black, fill=black!80, very thick](5.5,0) circle (0.07);
\filldraw[color=black, fill=black!80, very thick](5.5,0) circle (0.07);

\filldraw[color=red, fill=red!10, very thick](0,1) circle (0.07);
\filldraw[color=red, fill=red!10, very thick](1,1) circle (0.07);
\filldraw[color=red, fill=red!10, very thick](2,1) circle (0.07);
\filldraw[color=red, fill=red!10, very thick](3,1) circle (0.07);
\filldraw[color=red, fill=red!10, very thick](4,1) circle (0.07);
\filldraw[color=red, fill=red!10, very thick](5,1) circle (0.07);
\filldraw[color=red, fill=red!10, very thick](6,1) circle (0.07);
\filldraw[color=black, fill=black!80, very thick](0.5,1) circle (0.07);
\filldraw[color=black, fill=black!80, very thick](1.5,1) circle (0.07);
\filldraw[color=black, fill=black!80, very thick](2.5,1) circle (0.07);
\filldraw[color=black, fill=black!80, very thick](3.5,1) circle (0.07);
\filldraw[color=black, fill=black!80, very thick](4.5,1) circle (0.07);
\filldraw[color=black, fill=black!80, very thick](5.5,1) circle (0.07);
\filldraw[color=black, fill=black!80, very thick](5.5,1) circle (0.07);

\filldraw[color=red, fill=red!10, very thick](0,2) circle (0.07);
\filldraw[color=red, fill=red!10, very thick](1,2) circle (0.07);
\filldraw[color=red, fill=red!10, very thick](2,2) circle (0.07);
\filldraw[color=red, fill=red!10, very thick](3,2) circle (0.07);
\filldraw[color=red, fill=red!10, very thick](4,2) circle (0.07);
\filldraw[color=red, fill=red!10, very thick](5,2) circle (0.07);
\filldraw[color=red, fill=red!10, very thick](6,2) circle (0.07);
\filldraw[color=black, fill=black!80, very thick](0.5,2) circle (0.07);
\filldraw[color=black, fill=black!80, very thick](1.5,2) circle (0.07);
\filldraw[color=black, fill=black!80, very thick](2.5,2) circle (0.07);
\filldraw[color=black, fill=black!80, very thick](3.5,2) circle (0.07);
\filldraw[color=black, fill=black!80, very thick](4.5,2) circle (0.07);
\filldraw[color=black, fill=black!80, very thick](5.5,2) circle (0.07);
\filldraw[color=black, fill=black!80, very thick](5.5,2) circle (0.07);

\filldraw[color=red, fill=red!10, very thick](0,3) circle (0.07);
\filldraw[color=red, fill=red!10, very thick](1,3) circle (0.07);
\filldraw[color=red, fill=red!10, very thick](2,3) circle (0.07);
\filldraw[color=red, fill=red!10, very thick](3,3) circle (0.07);
\filldraw[color=red, fill=red!10, very thick](4,3) circle (0.07);
\filldraw[color=red, fill=red!10, very thick](5,3) circle (0.07);
\filldraw[color=red, fill=red!10, very thick](6,3) circle (0.07);
\filldraw[color=black, fill=black!80, very thick](0.5,3) circle (0.07);
\filldraw[color=black, fill=black!80, very thick](1.5,3) circle (0.07);
\filldraw[color=black, fill=black!80, very thick](2.5,3) circle (0.07);
\filldraw[color=black, fill=black!80, very thick](3.5,3) circle (0.07);
\filldraw[color=black, fill=black!80, very thick](4.5,3) circle (0.07);
\filldraw[color=black, fill=black!80, very thick](5.5,3) circle (0.07);
\filldraw[color=black, fill=black!80, very thick](5.5,3) circle (0.07);

\filldraw[color=Plum, fill=Plum!10, very thick](0,0.5) circle (0.07);
\filldraw[color=Plum, fill=Plum!10, very thick](1,0.5) circle (0.07);
\filldraw[color=Plum, fill=Plum!10, very thick](2,0.5) circle (0.07);
\filldraw[color=Plum, fill=Plum!10, very thick](3,0.5) circle (0.07);
\filldraw[color=Plum, fill=Plum!10, very thick](4,0.5) circle (0.07);
\filldraw[color=Plum, fill=Plum!10, very thick](5,0.5) circle (0.07);
\filldraw[color=Plum, fill=Plum!10, very thick](6,0.5) circle (0.07);
\filldraw[color=Plum, fill=Plum!10, very thick](0,1.5) circle (0.07);
\filldraw[color=Plum, fill=Plum!10, very thick](1,1.5) circle (0.07);
\filldraw[color=Plum, fill=Plum!10, very thick](2,1.5) circle (0.07);
\filldraw[color=Plum, fill=Plum!10, very thick](3,1.5) circle (0.07);
\filldraw[color=Plum, fill=Plum!10, very thick](4,1.5) circle (0.07);
\filldraw[color=Plum, fill=Plum!10, very thick](5,1.5) circle (0.07);
\filldraw[color=Plum, fill=Plum!10, very thick](6,1.5) circle (0.07);
\filldraw[color=Plum, fill=Plum!10, very thick](0,2.5) circle (0.07);
\filldraw[color=Plum, fill=Plum!10, very thick](1,2.5) circle (0.07);
\filldraw[color=Plum, fill=Plum!10, very thick](2,2.5) circle (0.07);
\filldraw[color=Plum, fill=Plum!10, very thick](3,2.5) circle (0.07);
\filldraw[color=Plum, fill=Plum!10, very thick](4,2.5) circle (0.07);
\filldraw[color=Plum, fill=Plum!10, very thick](5,2.5) circle (0.07);
\filldraw[color=Plum, fill=Plum!10, very thick](6,2.5) circle (0.07);
\node at (0,0.5) {+};
\node at (1,0.5) {+};
\node at (2,0.5) {+};
\node at (3,0.5) {+};
\node at (4,0.5) {+};
\node at (5,0.5) {+};
\node at (6,0.5) {+};
\node at (0,1.5) {+};
\node at (1,1.5) {+};
\node at (2,1.5) {+};
\node at (3,1.5) {+};
\node at (4,1.5) {+};
\node at (5,1.5) {+};
\node at (6,1.5) {+};
\node at (0,2.5) {+};
\node at (1,2.5) {+};
\node at (2,2.5) {+};
\node at (3,2.5) {+};
\node at (4,2.5) {+};
\node at (5,2.5) {+};
\node at (6,2.5) {+};

\filldraw[color=red, fill=red!10, very thick](-4,2) circle (0.07);
\filldraw[color=black, fill=black!80, very thick](-4,1.5) circle (0.07);
\filldraw[color=Plum, fill=Plum!10, very thick](-4,1) circle (0.07);
\node at (-4,1) {+};

\node at (-2.63,2) {\large $u,\;w_x,\;w_y,\;c$};
\node at (-3.5,1.5) {\large $v_x$};
\node at (-3.5,1) {\large $v_y$};
\draw[black] (0,-0.1) -- (0,-0.3) -- (0.5,-0.3) -- (0.5,-0.1);
\node at (0.25,-0.7) {\Large$\frac{\Delta x}{2}$};
\draw[black] (-0.1,0.5) -- (-0.3,0.5) -- (-0.3,1) -- (-0.1,1);
\node at (-0.7,0.75) {\Large$\frac{\Delta x}{2}$};
\end{tikzpicture}
\caption{Staggered grid in 2D with the relative positions of the (1SD) wave equations variables and parameters. $u,\;w_x,\;w_y$ and $c$ are collocated (centered) and $v_x$, $v_y$ are staggered in the grid.}\label{fig_staggered_2D}
\end{figure}
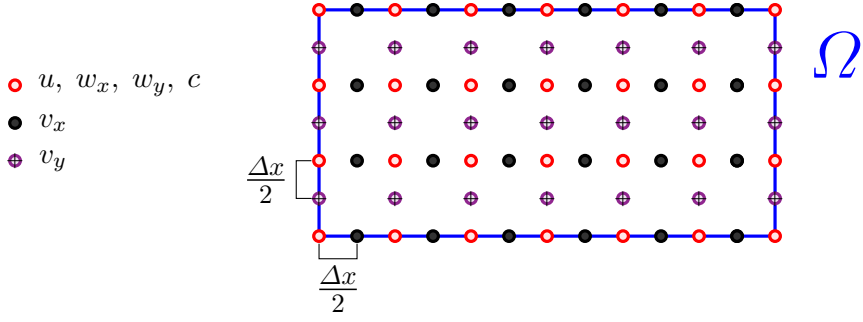

In this Section we present the basic information about the spatial discretization methods focused on classic finite difference schemes. To ensure an adequate representation of high frequency waves, we use a staggered grid, representing waves up to frequencies of $2/\Delta x$, improving spatial stability and dispersion properties \citep{wang2014seismic}.
Fig.\,\ref{fig_staggered_2D} depicts the variable positioning of 1SD equations for two-dimensions. Figures for other formulations are provided in App.\,\ref{sec_appendix_numerical_framework}.

We use discretizations of fourth and eighth order for the spatial derivatives.
Having a high order spatial discretization relates to a reduction of dispersion effects due to spatial discretizations \citep{liu2013globally}, supporting us to investigate errors in the time-domain of exponential integration schemes.

For collocated and staggered variables, the fourth-order scheme is given by

\begin{flalign}\label{eq_spatial_4th_1}
\frac{\partial u_{i+\frac{1}{2}}}{\partial x}&\approx \frac{u_{i-1}-27u_{i}+27u_{i+1}-u_{i+2}}{24\Delta x}&
\end{flalign}
and the eighth order scheme is given by

\begin{equation}
    \frac{\partial u_{i+\frac{1}{2}}}{\partial x}\approx \frac{1225}{1024\Delta x}\left(u_{i+1}-u_{i}-\frac{u_{i+2}-u_{i-1}}{15}+\frac{u_{i+3}-u_{i-2}}{125}-\frac{u_{i+4}-u_{i-3}}{1715}\right) \label{eq_spatial_8th_1}
\end{equation}
with analogous expressions for the $y$-coordinate in the 2D discretization.

For the points near the outer boundary, we use the same discretization formulas as for the interior points, but with zero-valued functions (Dirichlet Boundary condition) for points outside the domain.
After all, if a wave reaches the outer boundary it will be continuously weakened in its way to the boundary and back to the physical domain, attaining minimal energy.

For the PML thickness $\delta$, and the parameter $\beta_0$, \citet{assi2017compact} proposed a relation between them and the spatial grid space $\Delta x$. However, our objective in this paper is not the study of the PML absorbing boundary, but the solution of the wave equations with PML constraints.
Therefore, we choose suitable values so that $\delta$ is small, and a value of $\beta_0$ such that the waves reflections remains minimal, but without numerical errors because of large values of $\beta_0$. The values for these parameters and other numerical details can be found in Appendix \ref{sec_appendix_experiments_formulation}.

We continue with a description of various test cases used for further investigation of numerical experiments in one of the following sections.

\subsection{Test cases}\label{sec_test_cases}

To construct the operator, we also require the definition of the velocity field and other model parameters. Their description is organized in several numerical experiments, for the wave propagation equations with PML, that will be used throughout this paper. The numerical tests comprehend  scenarios with variable difficulty, changing the dimension of the wave equations, the characteristics of the medium, the initial conditions, and the use of a source term. The general features of the numerical experiments used in the remainder of this work are summarized in Table \ref{tab_num_experiments}, with further details provided in Appendix \ref{sec_appendix_experiments_formulation}.
The source term for the test cases consists in a Ricker's wavelet, since this is the one of the most frequently used source terms in seismic imaging \citep{ikelle2018introduction}.

{\renewcommand{\arraystretch}{1.2}
\begin{table}[tbh]
	\begin{center}
		\begin{tabular}{|c|c|c|c|c|c|}
		\hline
		\textbf{Test Case ID} &
		\textbf{Type} & \textbf{Dim} & \textbf{Medium} & \textbf{Initial cond.} & \textbf{Source term}\\
		\hline
		TC\#1 & acoustic & 1D & homog. & non-zero & -\\
		\hline
		TC\#2 & acoustic & 1D & heterog. & non-zero & -\\
		\hline
		TC\#3 & acoustic & 1D & heterog. & zero & yes\\
		\hline
		TC\#4 & acoustic & 2D & homog. & non-zero & -\\
		\hline
		TC\#5 & acoustic & 2D & heterog. & non-zero & -\\
		\hline
		TC\#6 & acoustic & 2D & heterog. & zero & yes\\
		\hline
		TC\#7 & elastic & 2D & heterog. & zero & yes\\
		\hline
		\end{tabular}
	\end{center}
	\caption{General features of the numerical experiments: type of equation, dimensions of the problem, medium heterogeneity (homogeneous or heterogeneous), initial conditions, and the use of the Ricker's source term.}\label{tab_num_experiments}
\end{table}}

The cases TC\#1 and TC\#4, where the wave is propagated in a \textit{homogeneous} medium \textit{without a source term} in one and two dimensions, are intended to analyze the Faber approximation of the wave equation with PML in a mathematical scenario of lowest complexity.
A non-zero initial condition is used in case of no source term.

\section{Spectrum of discrete operator}\label{sec_spectrum}

The Faber approximation requires an estimation of the convex hull of the spectrum of the matrix operator to construct the ellipse where the polynomials will be defined. This is a difficult task for the matrix operators derived from the spatial discretization of PDE systems, where the analytical expression of the eigenvalues is not generally known.
Since the matrix dimensions can be very large, the computational time to compute the eigenvalues can be inadequate for practical applications. In addition, when solving PDEs numerically, variations in the parameters of the equations often produce significant changes in the spectrum of the discrete operator. This is the cases in seismic imaging for the wave equations with PML, where the velocity field is constantly modified when solving the inverse problem.

Based on empirical results produced by the numerical considerations and the test experiments of the previous section, we propose an estimate of the spectrum of the discrete operator $\boldsymbol{H}$.
We start by studying the eigenvalue distribution of a lower-resolution discrete wave equations operator for which a computation of eigenvalues is possible.
By doing so, we aim for sharp bounds of the $\boldsymbol{H}$ operator spectrum $\sigma (\boldsymbol{H})$ for the construction of the optimal ellipse, which can be generalized to high-resolution discretizations.

\subsection{General properties}\label{sec_spectrum_prop}

\begin{figure}[H]
    \subfloat[1SD $\Delta x=0.105$]{\includegraphics[scale=0.185]{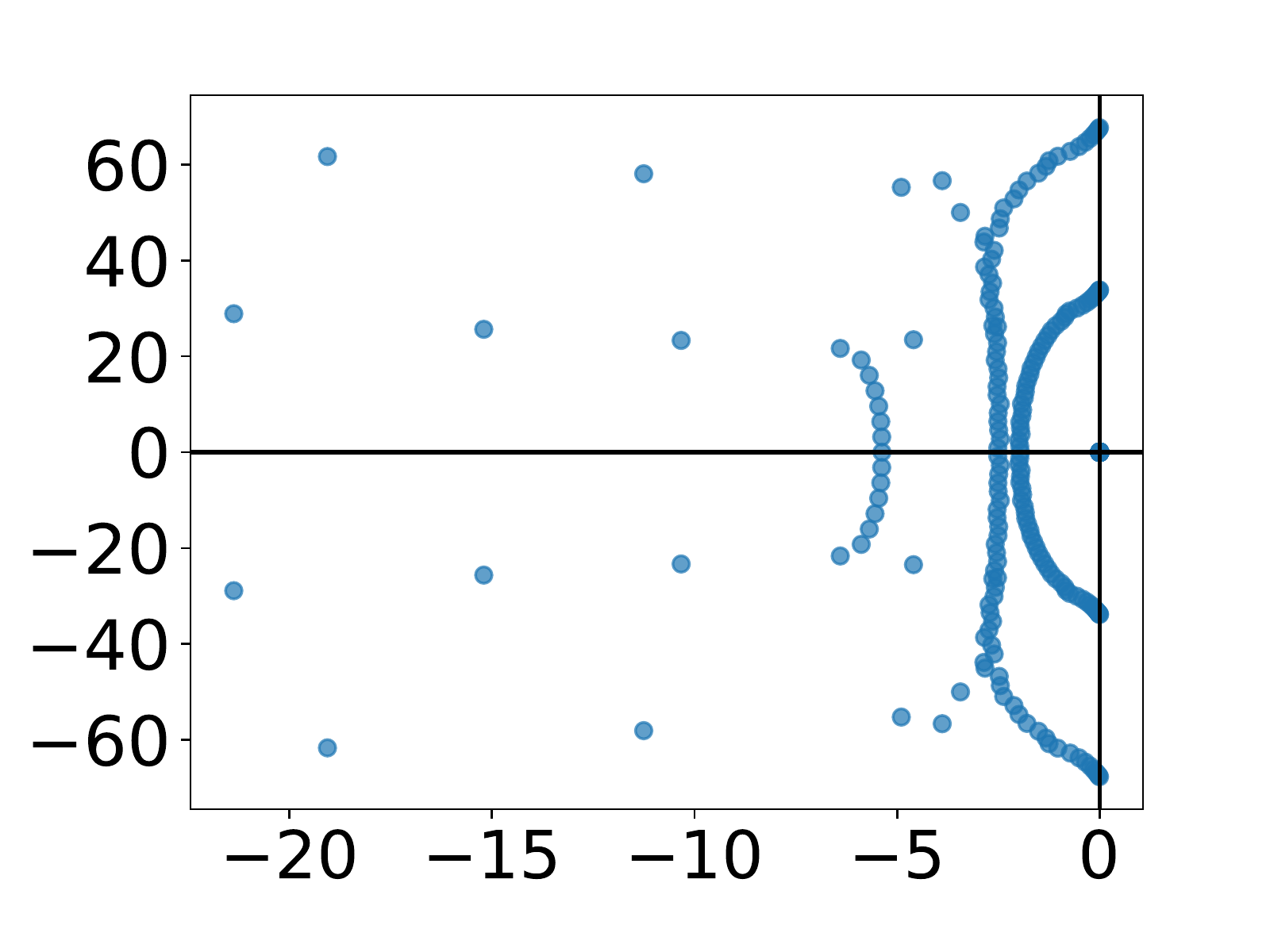}}
    \subfloat[1SD $\Delta x=0.021$]{\includegraphics[scale=0.185]{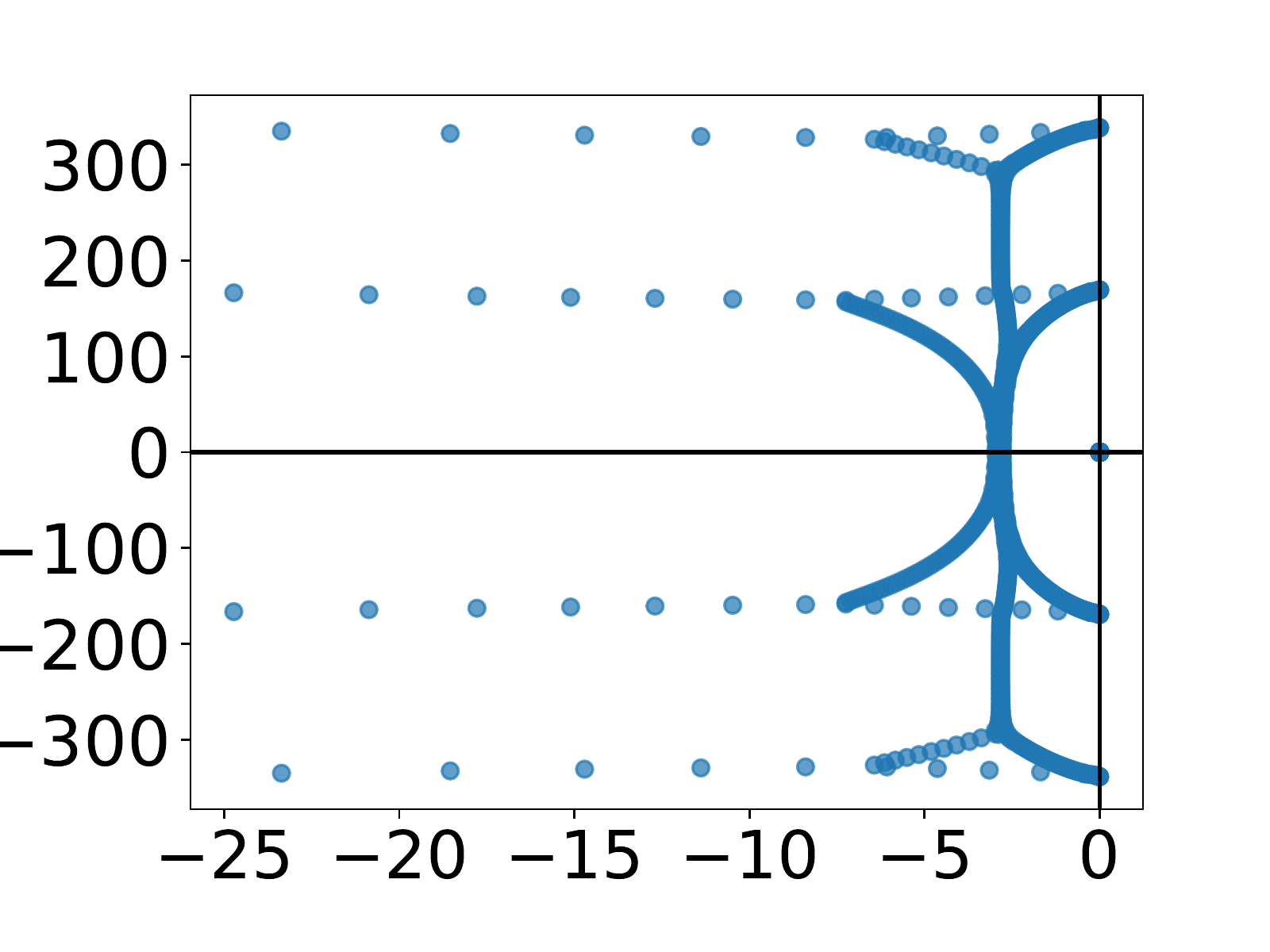}}
    \subfloat[1SD $\Delta x=0.0105$]{\includegraphics[scale=0.185]{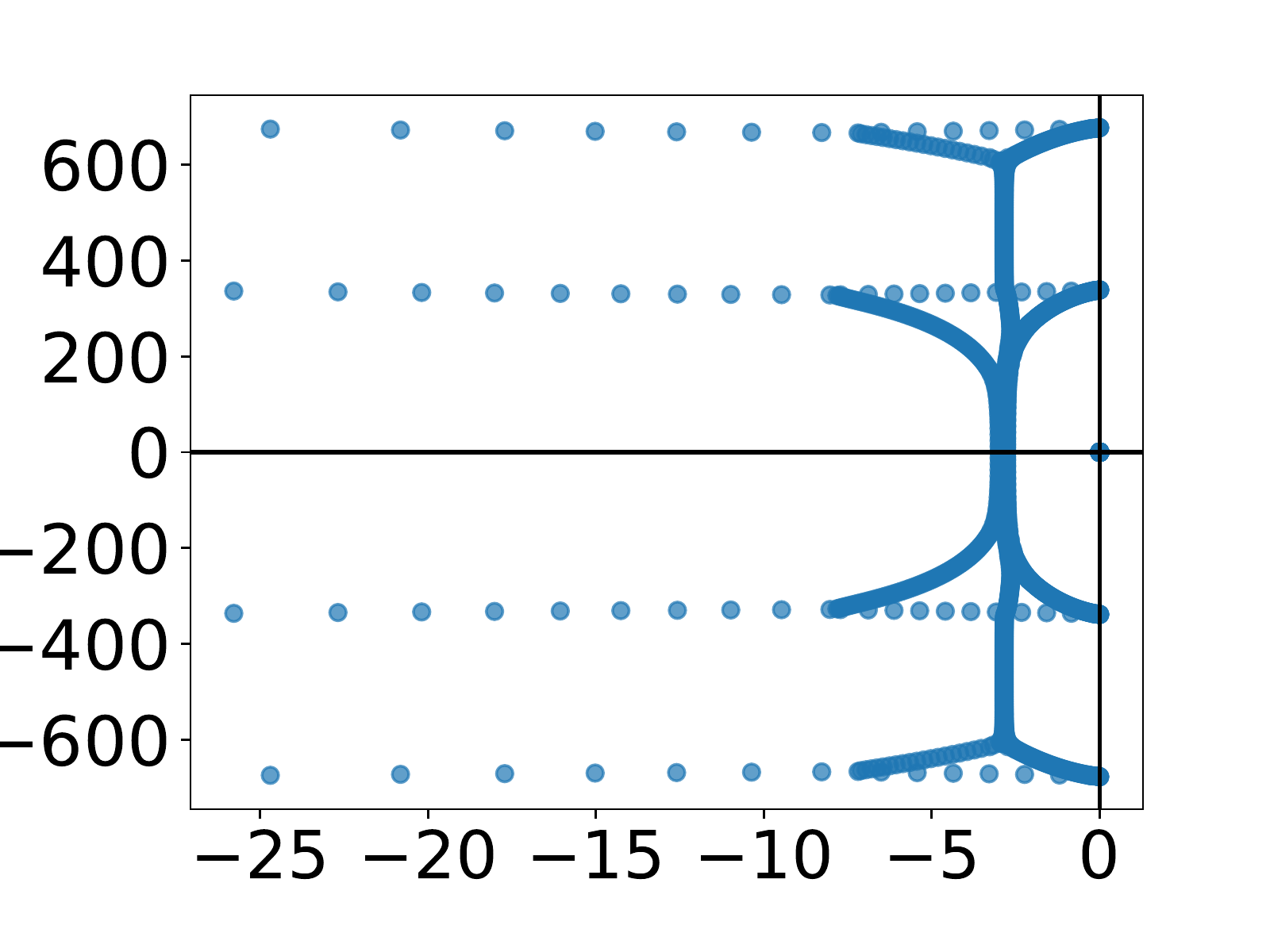}}
    \subfloat[1SD $\Delta x=0.0021$]{\includegraphics[scale=0.185]{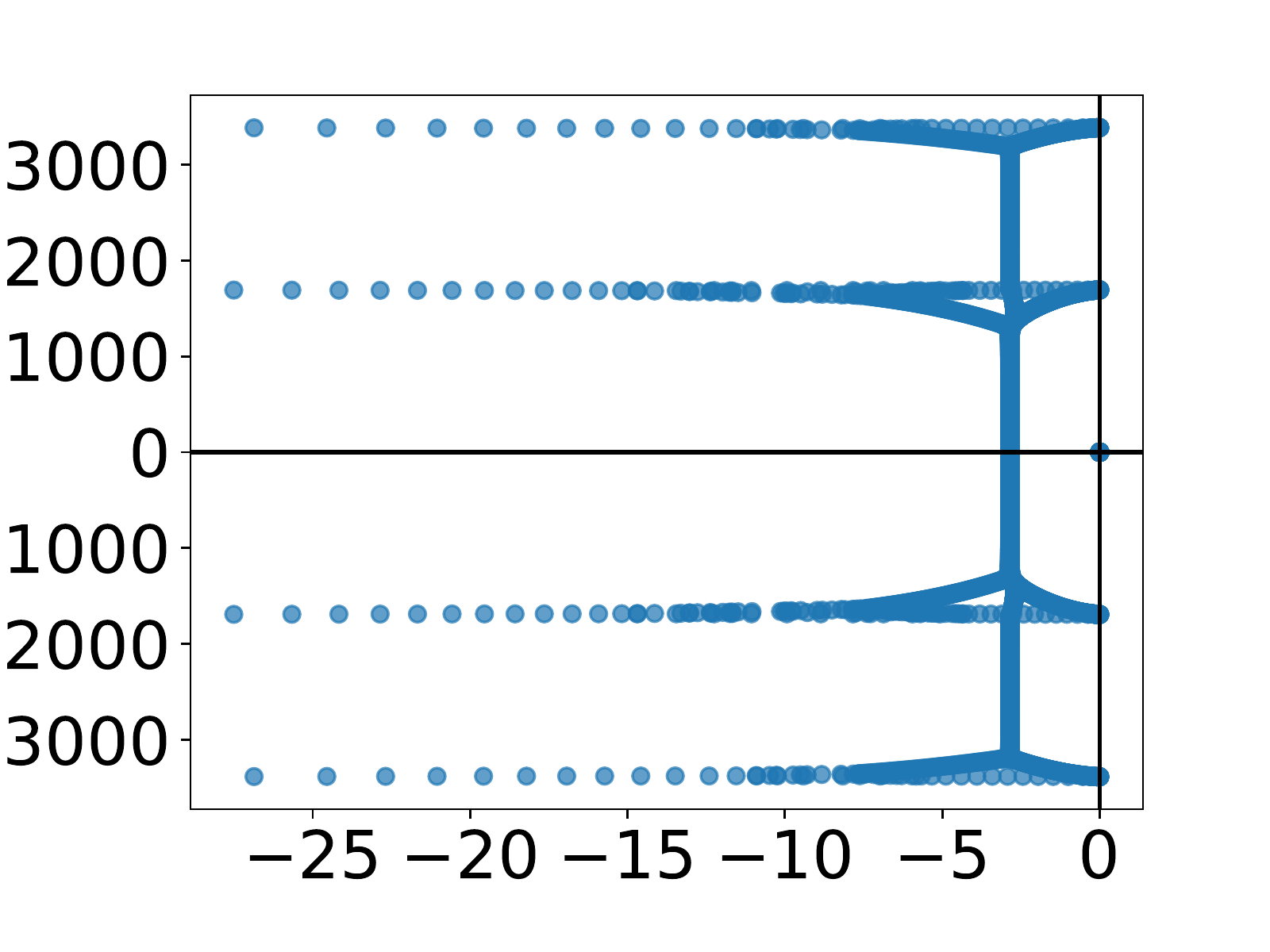}}\\
    \subfloat[2SD $\Delta x=0.105$]{\includegraphics[scale=0.185]{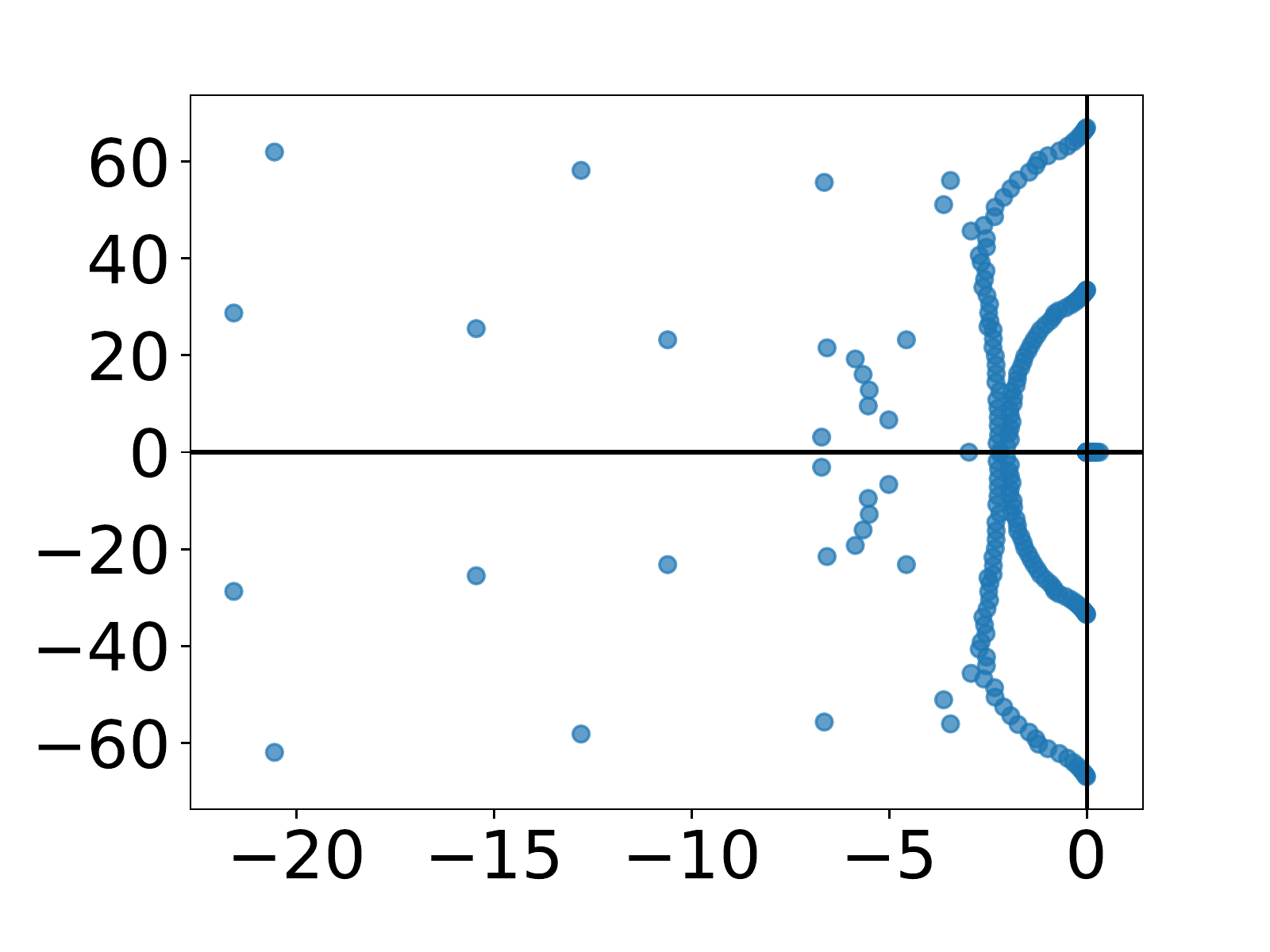}}
    \subfloat[2SD $\Delta x=0.021$]{\includegraphics[scale=0.185]{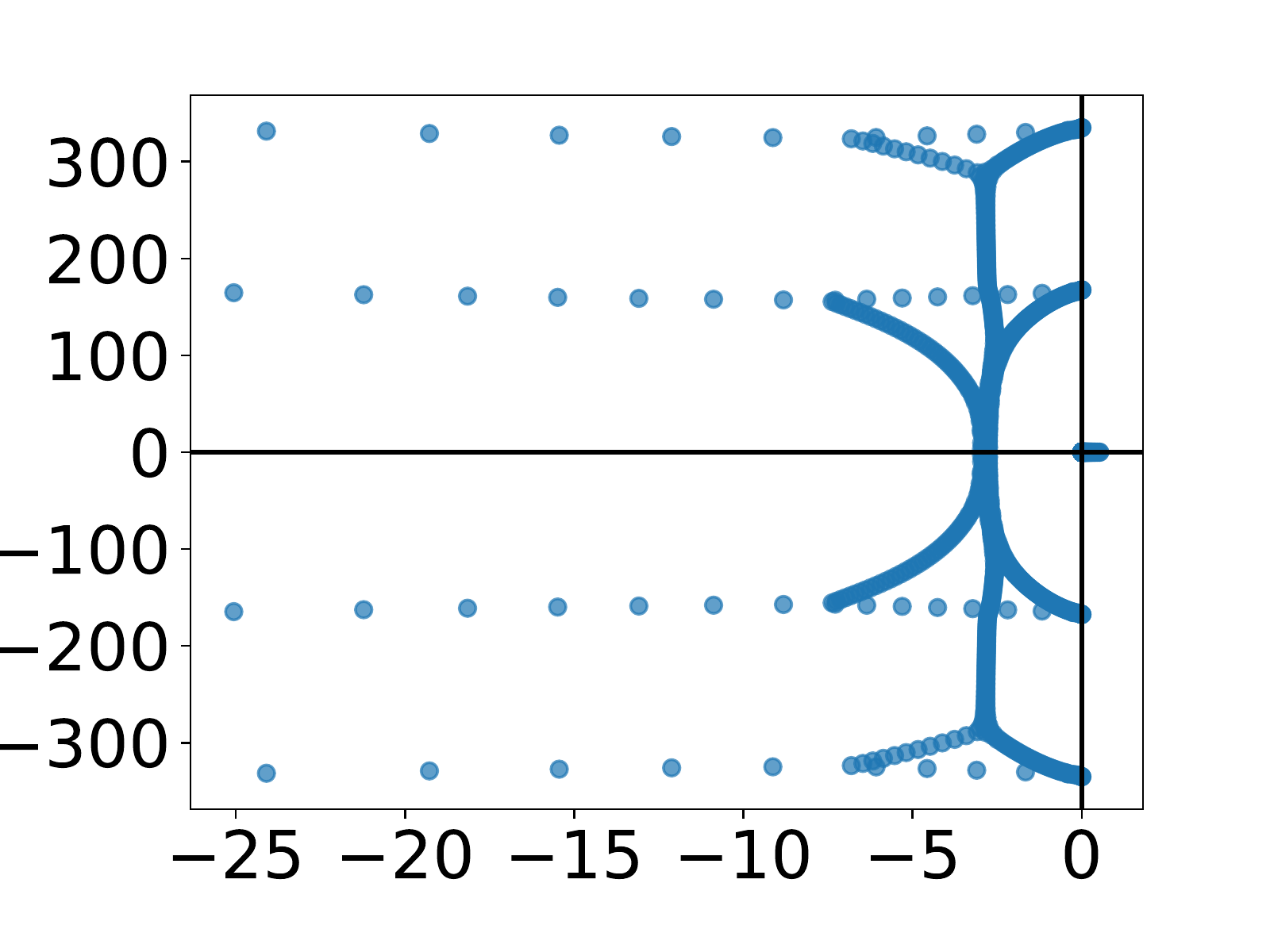}}
    \subfloat[2SD $\Delta x=0.0105$]{\includegraphics[scale=0.185]{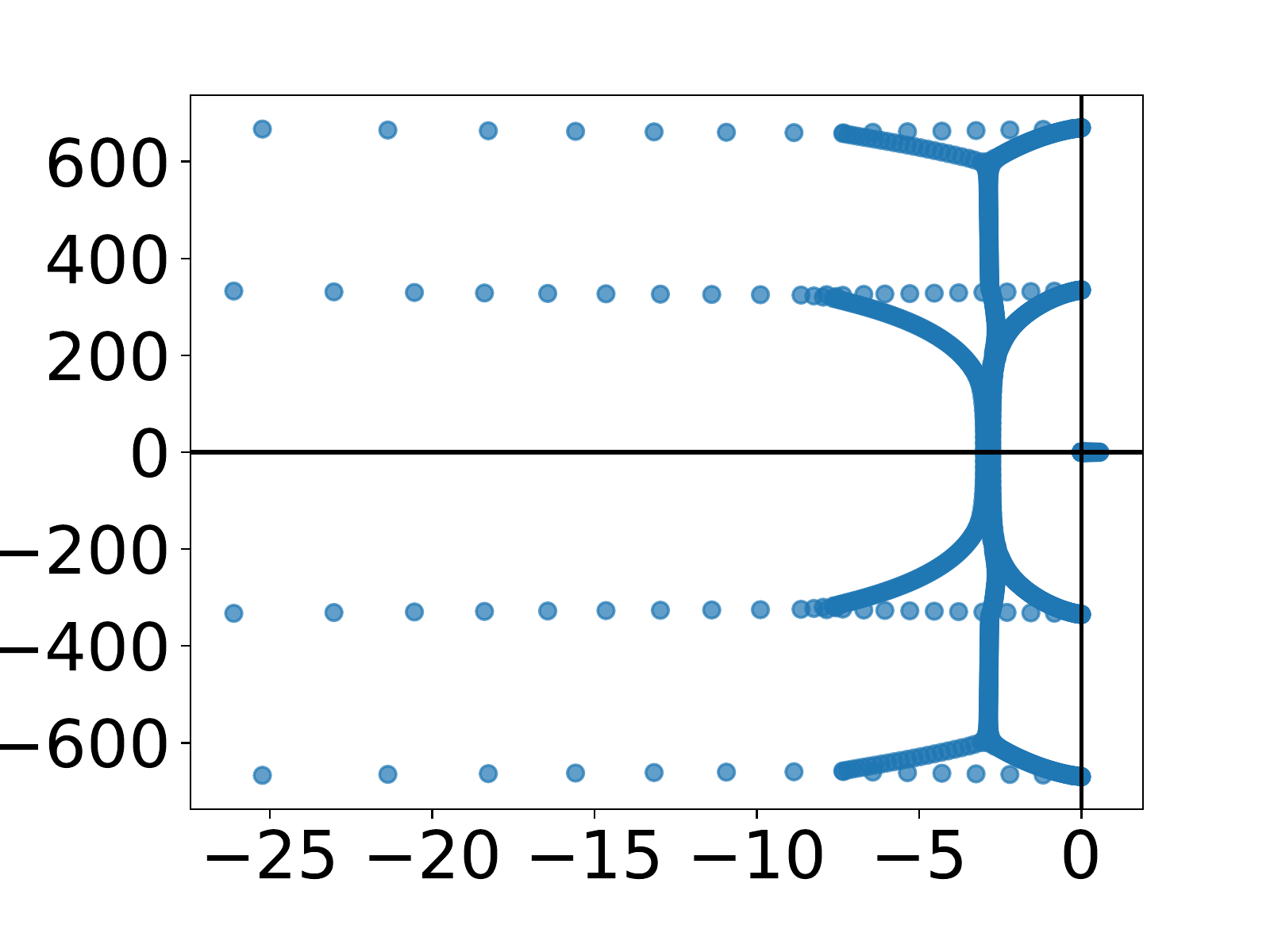}}
    \subfloat[2SD $\Delta x=0.0021$]{\includegraphics[scale=0.185]{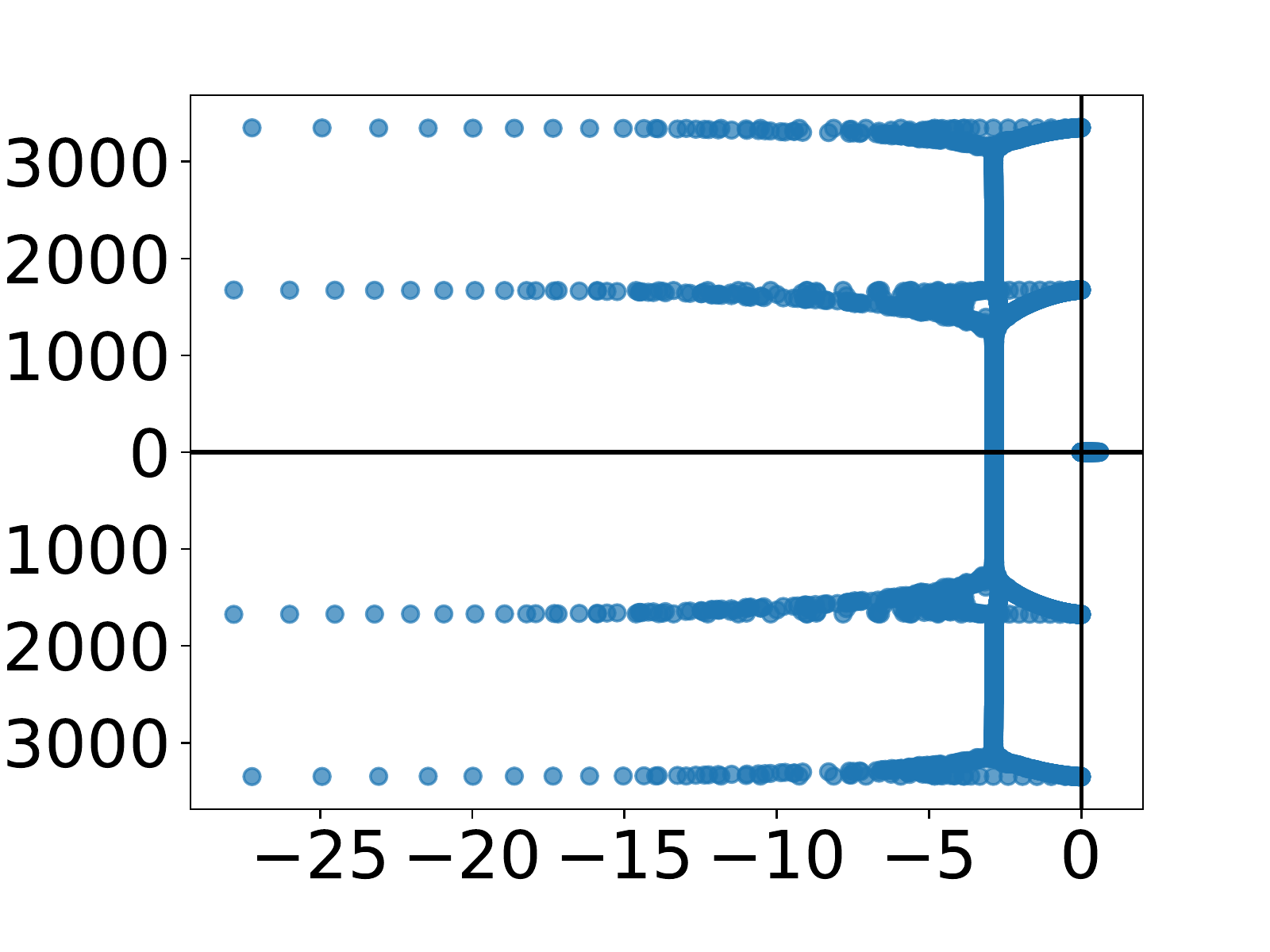}}
\caption{Eigenvalues on the complex plane (using coordinates $(Re(\lambda), Im(\lambda))$) using a fourth order spatial discretization of the acoustic wave equation in one dimension, with formulations 1SD (top row) and 2SD (bottom row) considering TC\#3, for $\Delta x=\{0.105,0.021,0.0105,0.0021\}$. Both formulations result in a similar spectrum, with a rectangular shaped convex hull for each $\Delta x$.}\label{fig_eigenvalues_distribution}
\end{figure}

We study the spectrum of $\boldsymbol{H}$ by calculating all of its eigenvalues for a finite decreasing sequence of $\Delta x$. Different distributions of eigenvalues on the complex plane for a 4th order spatial discretization with 1SD and 2SD formulations are given in Figure \ref{fig_eigenvalues_distribution}. When $\Delta x \rightarrow 0$ the convex hull of $\sigma(\boldsymbol{H})$ tends to a rectangle with sides parallel to the real axis. Thus, finding a relation between the rectangle sides and $\Delta x$ provides an estimator of the convex hull of $\sigma(\boldsymbol{H})$ for small $\Delta x$.

We also notice from Fig.\,\ref{fig_eigenvalues_distribution} that $\sigma(\boldsymbol{H})$ is symmetric with respect to the real axis. Since H is a real matrix, it is straightforward to verify that if\linebreak $\zeta\in\sigma(\boldsymbol{H})$, then $\overline{\zeta}\in\sigma(\boldsymbol{H})$. As a consequence, it is sufficient to only investigate the eigenvalues with non-negative imaginary part.
Moreover, Fig.\,\ref{fig_eigenvalues_distribution} indicates that the limits of the rectangle on the imaginary limit seems to relate linearly to $1/\Delta x$, but for the real part the relation is different, with an apparently constant negative limit on the left side, $-\beta_0$, for the PML parameter $\beta_0>0$, and a small positive number (for 2SD) or zero (for 1SD) on the right.

These bounds for the real axis are in agreement with the theoretical bounds for the continuous spectrum of the wave equations with PML conditions. For instance, performing a Fourier analysis of the eigenvalues of the continuum operator $\textsf{H}$ of the (1SD) formulation (Eq.\,\eqref{eq_1SD_1D}), in a unitary spatial domain ($[0,1]$), using
 \begin{equation}
     (u, \; v, \; w)^T=\sum\limits_{k=-\infty}^{\infty}(A_k,\; B_k, \; C_k)^Te^{ikx},
 \end{equation}
 and considering the linear map of the continuum operator $\textsf{H}$ on each term, we obtain the symbol of the operator $\textsf{H}$ \citep{strikwerda2004finite} as
\begin{equation}
    \textsf{H}
    \begin{pmatrix}
    	A_k\\ 
    	B_k\\
    	C_k
    \end{pmatrix}e^{ikx}=\begin{pmatrix}
    0\;\; & c^2ik & \;\;-c^2\\
    ik\;\; & -\beta_x(x) & \;\;0\\
    0\;\; & \beta_x(x)ik & \;\;-\beta_x(x)
    \end{pmatrix}\begin{pmatrix}
    A_k\\ 
    B_k\\
    C_k
    \end{pmatrix}e^{ikx}.
\end{equation}

Now, the eigenvalues of the operator symbol are given by 
 \begin{equation}
     \lambda_0=0,\;\lambda_{1,2}=-\beta_x(x)\pm ick.
 \end{equation}
 Therefore, if the wave number $k$ is constrained within $[-\frac{1}{\Delta x},\;\frac{1}{\Delta x}]$, as it is the case for discrete representations of the domain, the bounds of  $\sigma(\textsf{H})$ on the complex plane are given by $[-\beta_{\text{max}}, 0] \times [-\frac{c_{\text{max}}}{\Delta x},\frac{c_{\text{max}}}{\Delta x}]$, where $c_{\text{max}}$ is the maximum velocity in the medium and $\beta_{\text{max}}$ is the maximum of $\beta_x(x)$. So, we have a linear relation of the imaginary limits of the axis with respect to $\frac{1}{\Delta x}$, and the real part is within the interval $[-\beta_{\text{max}},0]$, where
 \begin{equation}\label{eq_bound_eigen_real}
     \beta_\text{max}=\beta_0\left(\frac{\delta-\Delta x/2}{\delta}\right)^2,
 \end{equation}
where $\beta_0$ and $\delta$ are PML parameters. 

We would like to remark here that although we only presented the full spectrum of the specific case of the 4th order spatial discretization and TC\#3. The previous analysis is also valid for all the other test cases, regardless if it used a 4th or an 8th order of spatial discretization, one or two dimensions, or the type of equation formulation, where we skipped these plots for sake of brevity.

In the next two subsections, we will investigate the dependency of the imaginary and real limits on the problem variables.

\subsection{Estimation of the imaginary limit}

\begin{figure}[tbh]
	\begin{tabular}{c|c|c}
		(a) TC\#3, in 1D, acoustic.
		&
		(b) TC\#5, in 2D, acoustic.
		&
		(c) TC\#7, in 2D, elastic.
		\\
	    \includegraphics[scale=0.24]{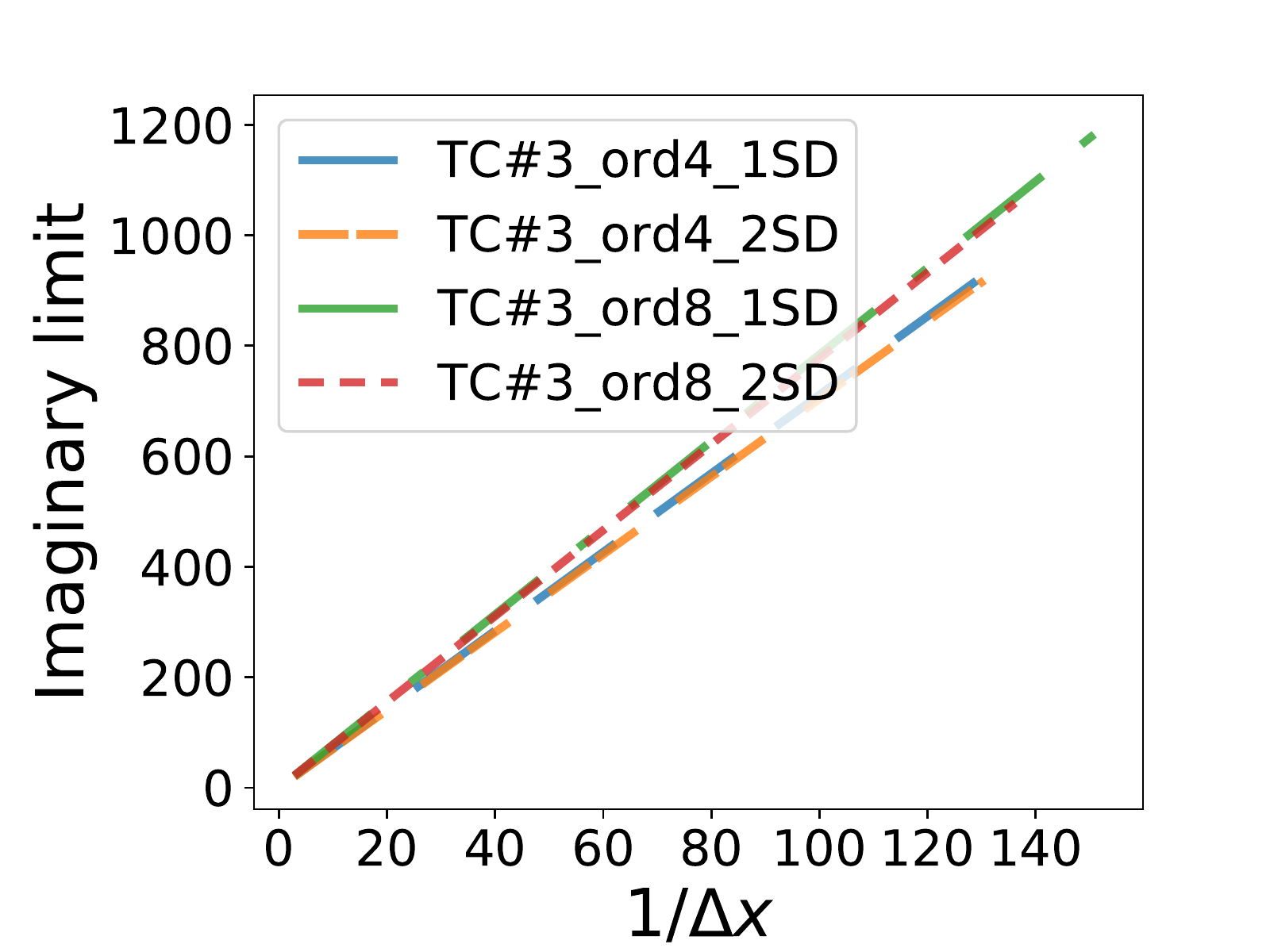}\hspace{-0.3cm}
	 	&
	    \includegraphics[scale=0.24]{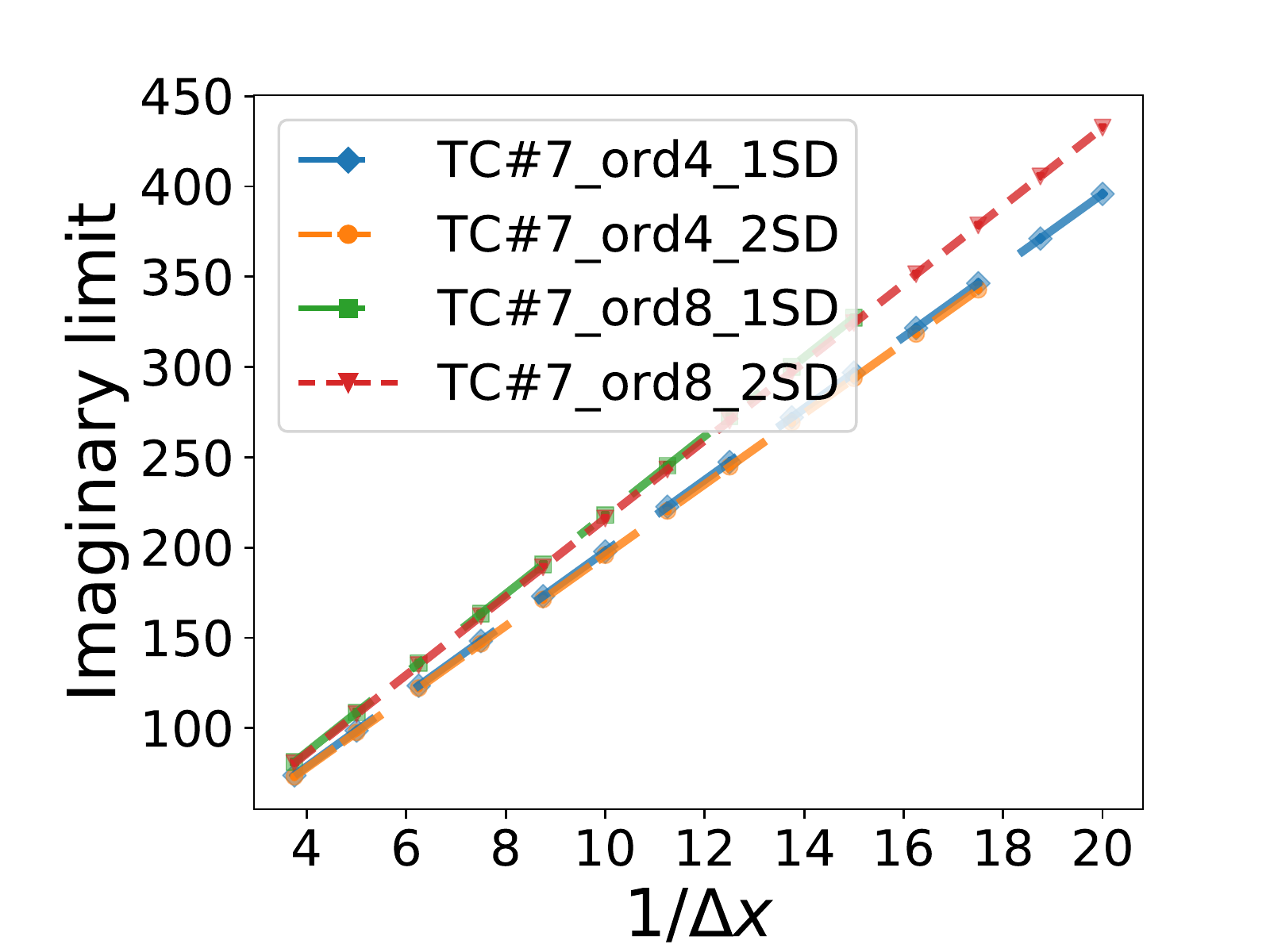}\hspace{-0.3cm}
	    &
	    \includegraphics[scale=0.24]{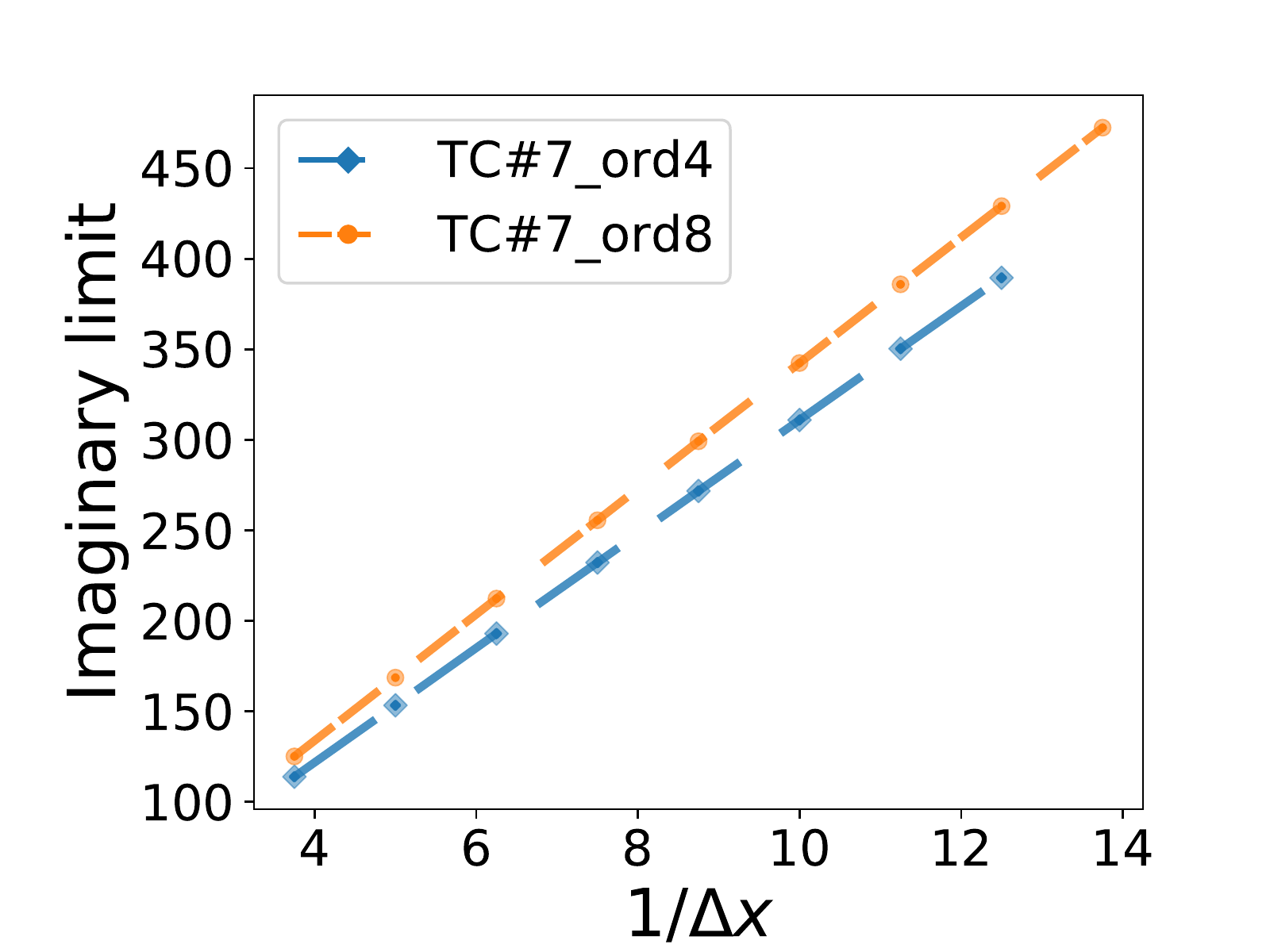}\hspace{-0.3cm}
    \end{tabular}
	\caption{Maximum imaginary parts of the eigenvalues of $\sigma(\boldsymbol{H})$ for varying $1/\Delta x$. The plot indicates a linear relation between the maximum imaginary part of the eigenvalues and $1/\Delta x$.
		\label{fig_imag_limits}}
\end{figure}

\begin{figure}[hbt]
	\subfloat[PML parameter variations of TC\#3.]{\includegraphics[scale=0.36]{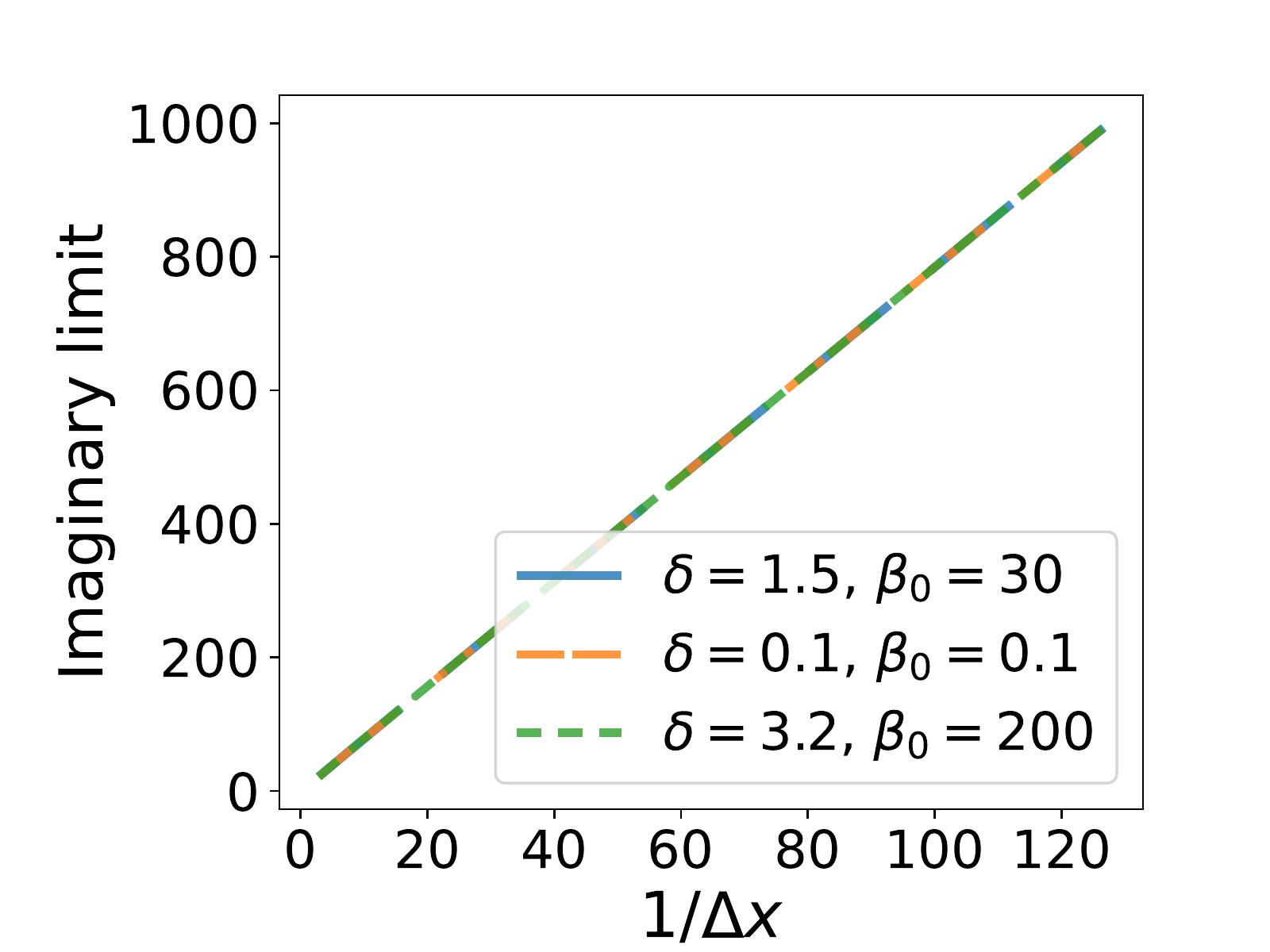}}
	\hfill
	\subfloat[Velocity field variations.]{\includegraphics[scale=0.36]{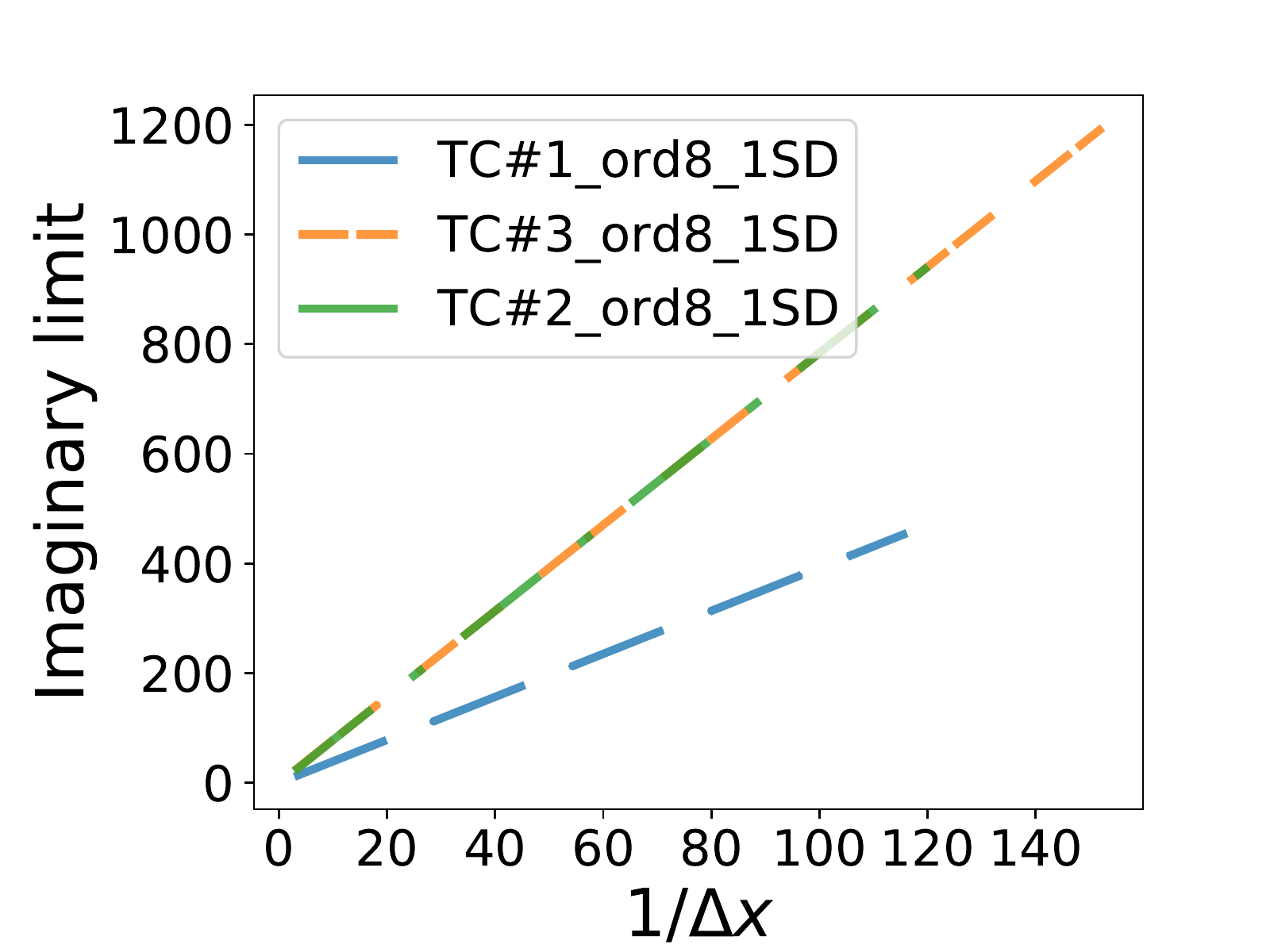}}
	\caption{Maximum imaginary parts of the eigenvalues of $\sigma(\boldsymbol{H})$ for varying $\Delta x$.
	(Left figure) Varying PML parameters $\delta$ and $\beta_0$. Variations of the PML parameters do not affect the maximum imaginary part of the discrete operator.
	(Right figure) Different velocity fields (test cases). If the maximum velocity $c_\text{max}$ does not change, the maximum imaginary parts remain unaltered.
}\label{fig_pml_variations}
\end{figure}

To empirically validate the linear relation between the maximum imaginary part of $\sigma(\boldsymbol{H})$ and $1/\Delta x$, we perform experiments varying the resolution and other parameters, such as the medium velocity field, the equation formulation, and the discretization order, with results shown in Figure \ref{fig_imag_limits}. The linear relations are clear in all experiments performed, with variations on the slopes, due to the different maximum velocities (from different test cases) and discretization orders, in agreement with the expected theory discussed in Sec \ref{sec_spectrum_prop}. However, we do not notice dependence on the model formulation (1SD or 2SD). The dependency on the model velocity is more clearly shown in the right panel of Fig.\,\ref{fig_pml_variations}.
In this case, the imaginary limit is connected linearly with the maximum medium velocity $c_\text{max}$.
The velocities $c_\text{max}$ are the same for TC\#2 and TC\#3, but differ by a factor of two with respect to $c_\text{max}$ of TC\#1.
This relation is also reflected by the imaginary limits of their respective operator spectrum, where the slope of the curves are 7.84 (for TC\#2 and TC\#3, which have maximum velocity given by $c_\text{max}=3.048$) and 3.92 (for TC\#1, with maximum velocity $c_\text{max}=1.524$), approximately.

Next, we investigate the influence of the PML parameters, with results given in the left image of Fig.\,\ref{fig_pml_variations}.
We can observe a superposition of the three lines, indicating a lack of dependence of the PML parameters. This was also observed for other test cases and parameter choices (not shown).

To briefly summarize, the particular slope of each curve depends, in decreasing importance, on the maximum medium velocity $c_\text{max}$, followed by the dimension of the problem, then the spatial discretization scheme and finally the formulation of the equations (for the acoustic case). 
Finally, the linear behavior between the maximum imaginary eigenvalue and $1/\Delta x$ allows determining the maximum imaginary eigenvalue for high resolution discretizations based an eigenvalue computation on a low resolution discretization.

\subsection{Estimation of the real limit}

\begin{figure}[H]
    \subfloat[TC\#3, in 1D, acoustic.]{\includegraphics[scale=0.25]{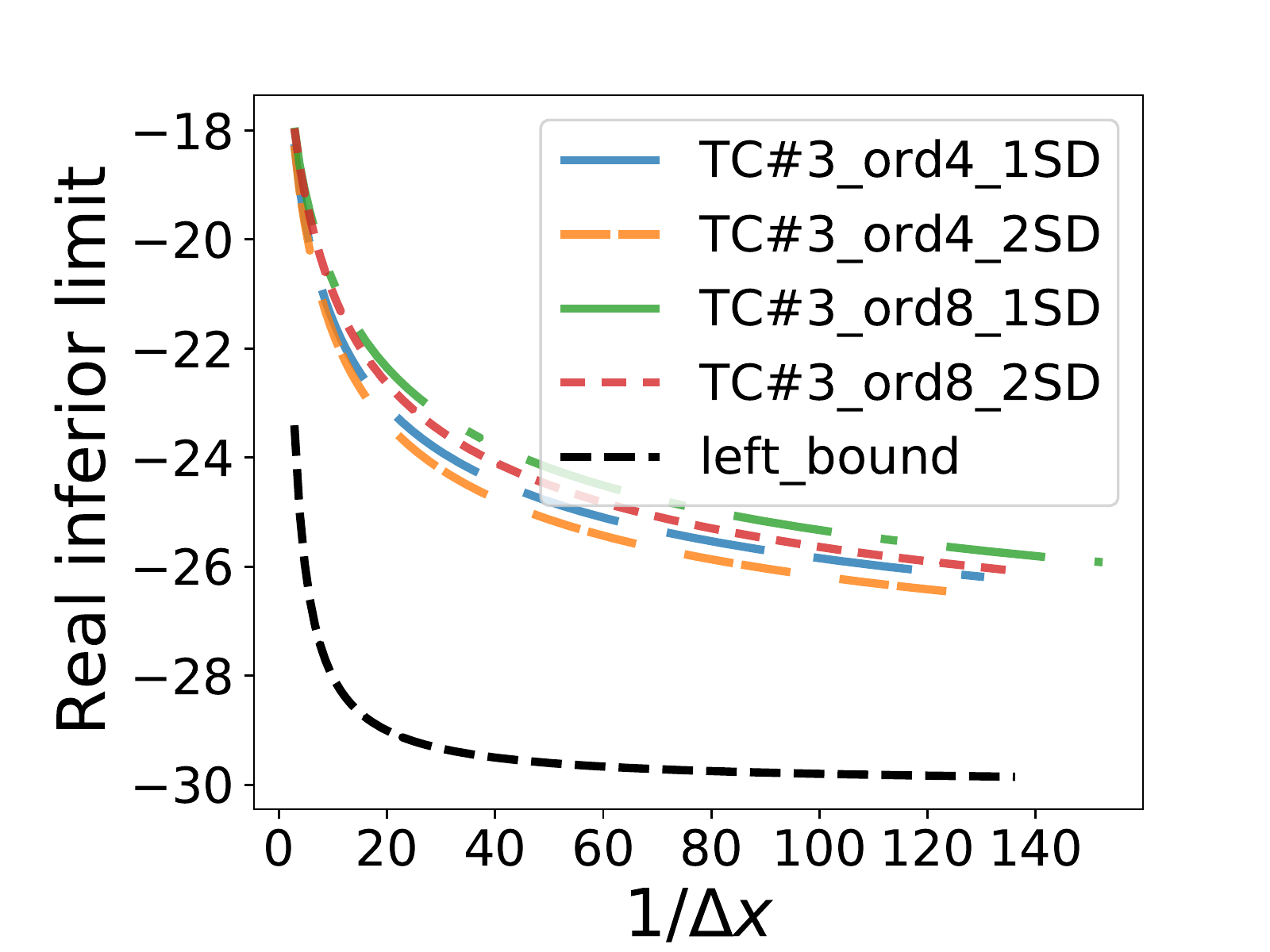}}
    \subfloat[TC\#5, in 2D, acoustic.]{\includegraphics[scale=0.25]{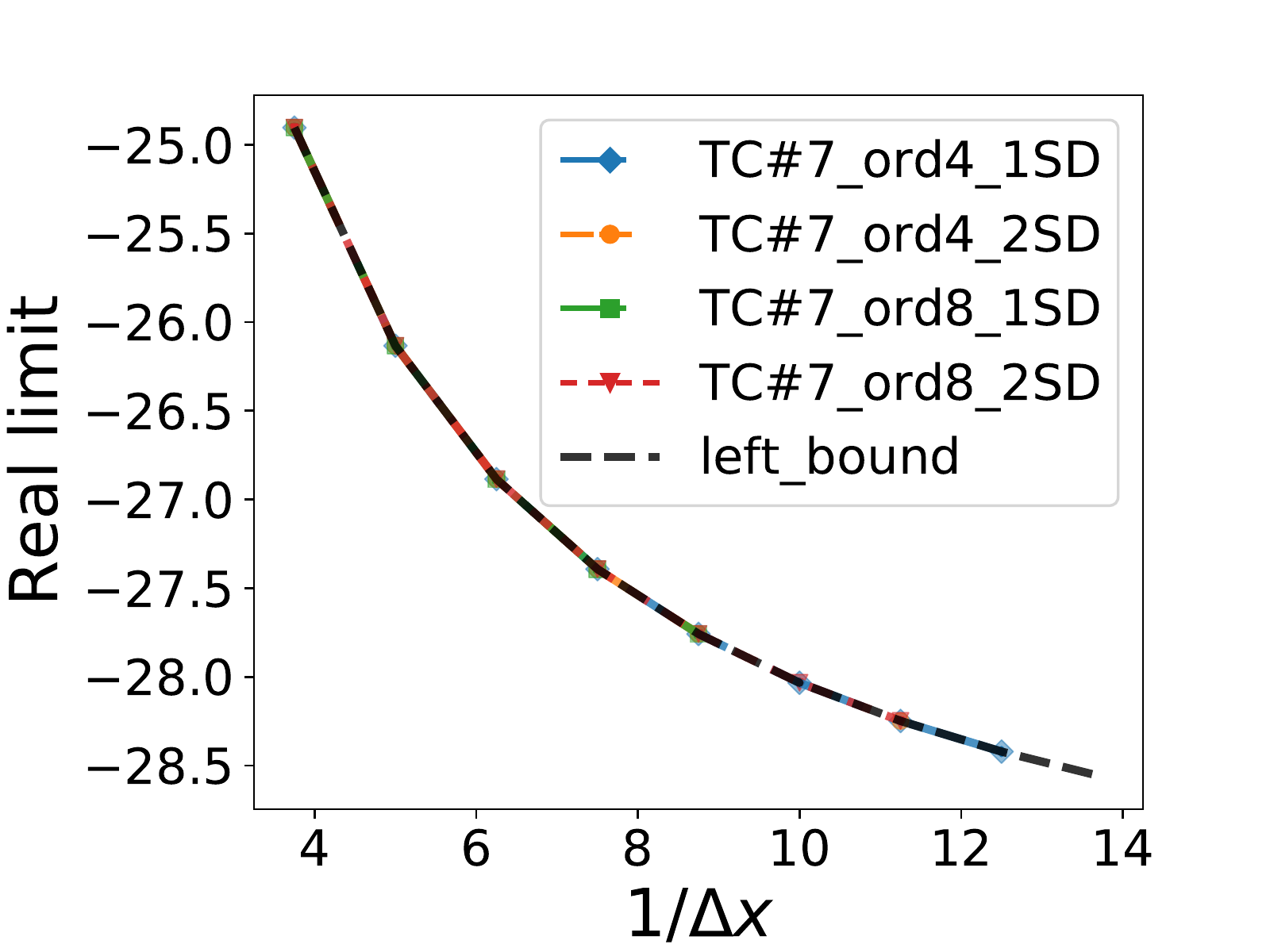}}
    \subfloat[TC\#7, in 2D, elastic.]{\includegraphics[scale=0.25]{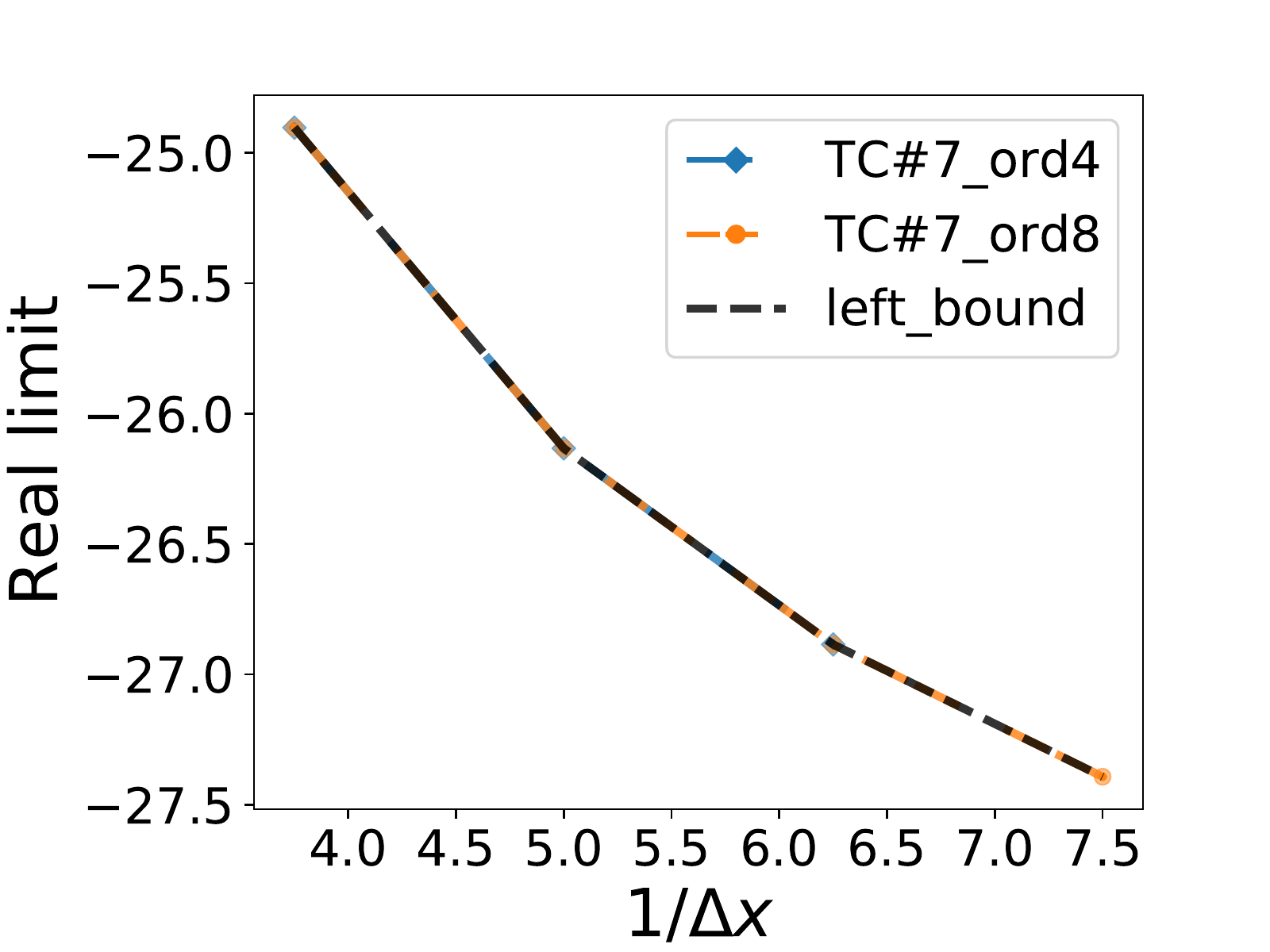}}
\caption{Lower bound of the real part of the eigenvalues for varying $\Delta x$, for 1SD, 2SD and elastic formulations, different dimensions and experimental sets.
The dashed line is our estimated lower bound given by Eq.\,\eqref{eq_bound_eigen_real}.}\label{fig_real_limits}
\end{figure}

Next, we compare the theoretical lower bound of the real part of the eigenvalues given by Eq.\,\eqref{eq_bound_eigen_real} with different test cases with results given in Fig.\,\ref{fig_real_limits}, where we vary $\Delta x$, and use different dimensions, formulations, spatial discretizations, and equations parameters.

For the 1D problems (left image of Fig.\,\ref{fig_real_limits}), we observe that our estimate provides an adequate lower bound, but it is overly pessimistic.

For the 2D problems (middle and right image), we observe very good matches, hence, providing a very sharp lower bound.

With respect to the  upper bound, we observed that the limit slightly depends on the particular model formulation (not shown), but always results in a value close to zero, or zero.
For the acoustic 1SD formulation and the elastic equations, the right limit is zero in all the experiments.
For the acoustic 2SD equations, empirical studies, as the one shown in Figure\,\ref{fig_eigenvalues_distribution}, indicate that the upper bound is a positive small number, smaller than 1. Since in all experiments this upper bound is always close to zero, this will not affect the estimates of the optimal ellipse for Faber polynomial approximations, as the ellipse size will be dominated by the imaginary axis bounds and the lower bound in the real axis. Therefore, a precise upper real bound will not be further required.

Overall, we have found a bound for the real part that only depends on the PML parameters, and which are independent on the velocity field, the equation formulation, and the spatial discretization scheme.

\subsection{Construction of enclosing rectangle}

We close this Section with some final remarks on the overall construction of the rectangle enclosing the spectrum $\sigma(\boldsymbol{H})$.
For the imaginary limit, we can use the linearity between the maximum imaginary part of the spectrum and $\Delta x$, where the linear relation can be computed with an Eigenvalue computation on a sufficiently low resolution discretization.
For the real part, we found an explicit lower bound, based on Eq.\,\eqref{eq_bound_eigen_real}, and an upper bound close to, or exactly, zero.
Finally, we can determine the enclosing rectangle, which will be the basis for the construction of the optimal ellipse with the Faber polynomials.

\section{Fourier's stability and dispersion results}\label{sec_stability_dispersion}

Next, we investigate two fundamentally important properties. First, the numerical stability using a von Neumann approach in the next subsection and, after this, the numerical dispersion, which is considered to be of high relevance in seismic imaging.
We analyze both aspects for several degrees of Faber polynomials separately, aiming to define and compute optimal criteria.

\subsection{Von Neumann stability analysis}\label{sec_stab}

In this Section we investigate the stability with Faber polynomials and estimate the CFL number ($c_{\text{\tiny CFL}}$) for each polynomial degree by performing a classical von Neumann analysis.
We start by replacing the absorbing boundary condition with periodic ones and drop the source term.
This reduces the equations to solve a purely hyperbolic (purely oscillatory, e.g. \citep{cox2002exponential}). We also assume a constant velocity profile (homogenous medium).
When expressed in Fourier series, the solution is given by
\begin{equation}\label{eq_error_stab}
    U(t,x)=\sum\limits_{m=0}^M{\begin{pmatrix}
    A_m(t)e^{i k_m x}\\
    B_m(t)e^{i k_m x}
    \end{pmatrix} },
\end{equation}
where $M$ is the number of frequencies considered in the solution, $k_m$ are the wave numbers, and the terms in the first and second rows stands for the solutions in the $u$ and $v$ variables (see Eqs.\,\eqref{eq_2SD_1D}-\eqref{eq_1SD_1D}), respectively.
This particular form of the linear wave equations allows an analysis of each spectral mode ($A_m, B_m$) separately.
Then, depending on the spatial scheme and formulation used, there will be a different matrix operator $\boldsymbol{G}$ known as amplification matrix (or stability function for this particular mode), such that
\begin{align}\label{eq_stability_op}
\begin{pmatrix}
A_m(t_{n+1})e^{ik x_i}\\
B_m(t_{n+1})e^{ikx_{i+1/2}}
\end{pmatrix}&=\boldsymbol{G}\begin{pmatrix}
A_m(t_{n})e^{ikx_i}\\
B_m(t_{n})e^{ikx_{i+1/2}}
\end{pmatrix}.
\end{align}

For the 1SD system, the amplification matrix $\boldsymbol{G}$ of Faber polynomial methods may be written as
\begin{equation}
    \label{eq_operator_D}
    \boldsymbol{G}=\sum\limits_{j=0}^m a_j(\Delta t \boldsymbol{H}) \boldsymbol{F}_j(\Delta t \boldsymbol{H}),
\end{equation}
where $\boldsymbol{H}$ is the right-hand side operator of equations formulation 1SD (see Eq.\,\eqref{eq_1SD_1D}), but without the PML term, and the coefficients $a_j$ are determined by Eq.\,\eqref{eq_intro_faber_coeff}.
As an intermediate step, we therefore investigate the application of $\Delta t\boldsymbol{H}$, which is a fundamental basic building block of polynomial approximations used here.

In what follows, we provide an example of a spatial fourth order finite-difference approximation \eqref{eq_spatial_4th_1}, where we obtain
\begin{align*}
    \Delta t\boldsymbol{H}\begin{pmatrix}
A_m(t)e^{i k x}\\
B_m(t)e^{(i+1/2) k x}
\end{pmatrix}&=\frac{\Delta t}{24\Delta x}\begin{pmatrix}
    0 & c^2g_1\\
    g_2 & 0
    \end{pmatrix}\begin{pmatrix}
A_m(t)e^{i k x}\\
B_m(t)e^{(i+1/2) k x}
\end{pmatrix}\\
&=\frac{\alpha}{24}\begin{pmatrix}
    0 & cg_1\\
    \frac{1}{c}g_2 & 0
    \end{pmatrix}\begin{pmatrix}
A_m(t)e^{i k x}\\
B_m(t)e^{(i+1/2) k x},
\end{pmatrix}
\end{align*}
where
\begin{align}
    \alpha&=\frac{c\Delta t}{\Delta x},\quad \theta=k\Delta x\in[0,\pi],\label{eq_alpha}\\
    g_1&=e^{-2i\theta}-e^{i\theta}+27(1-e^{-i\theta})\\
    g_2&=e^{-i\theta}-e^{2i\theta}+27(e^{i\theta}-1),
\end{align}
and we are considering the wavenumber $k\in[0,\pi/\Delta x]$. This leads to the representation of $\Delta t \boldsymbol{H}$
\begin{equation}\label{eq_deltaH}
    \Delta t\boldsymbol{H}=\frac{\alpha}{24}\begin{pmatrix}
    0 & g_1\\
    g_2 & 0
    \end{pmatrix},
\end{equation}
which also holds for other spatial discretization orders, but in other cases we may have different $g_1$ and $g_2$ values.

Given this form, we build the operator $\boldsymbol{G}$ by substituting \eqref{eq_deltaH} into expression \eqref{eq_operator_D}.
For the operators derived from the other formulations, spatial schemes, and dimensions, the reader can refer to Section \ref{sec_appendix_stability_dispersion}.

\begin{figure}[tbh]\hspace{-0.5cm}
	\begin{tabular}{c|c|c}
		(a) One dimension, acoustic.
		&
		(b) Two dimensions, acoustic.
		&
		(c) Two dimensions, elastic.
		\\
	    \includegraphics[scale=0.24]{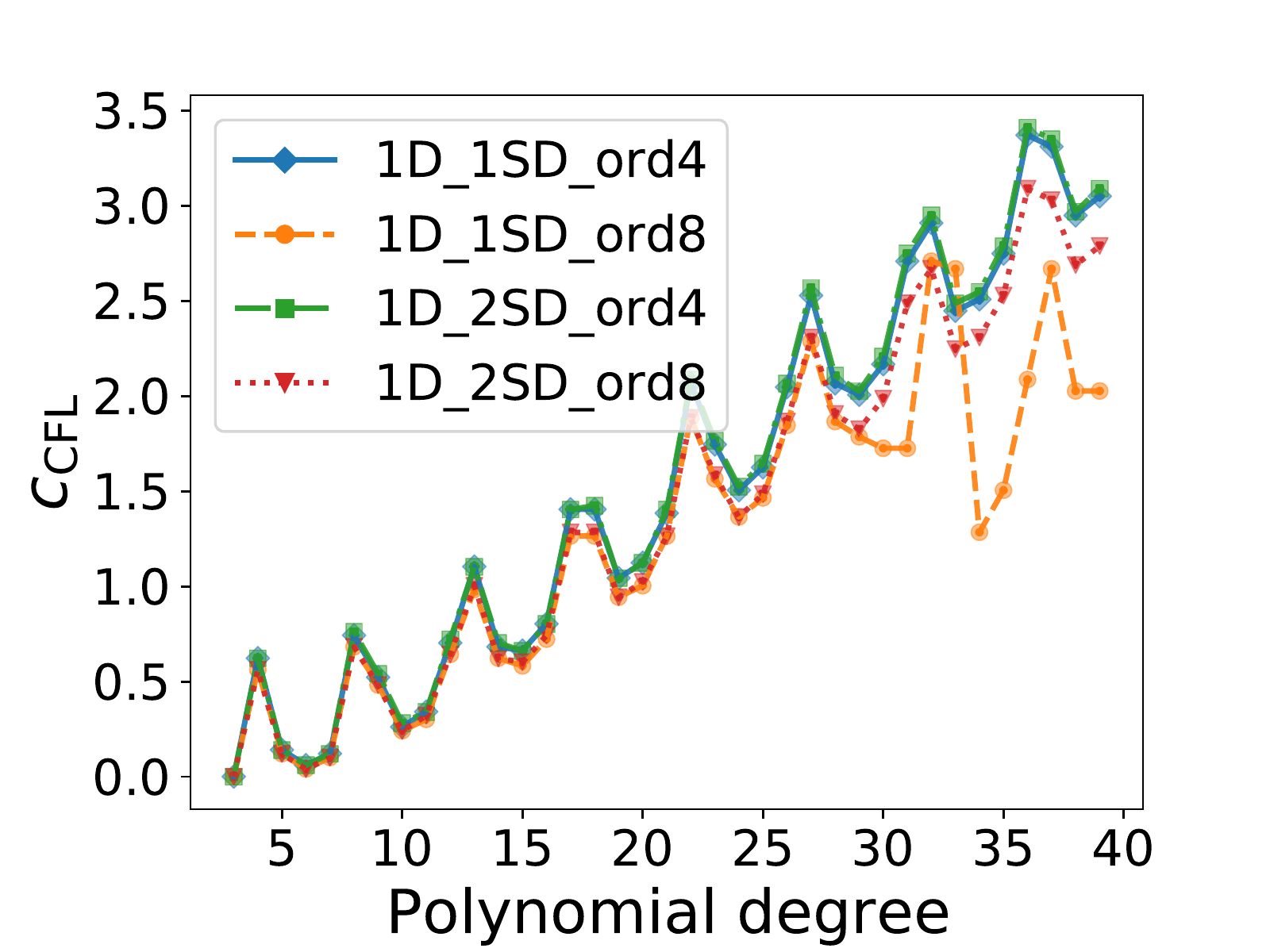}\hspace{-0.3cm}
	 	&
	    \includegraphics[scale=0.24]{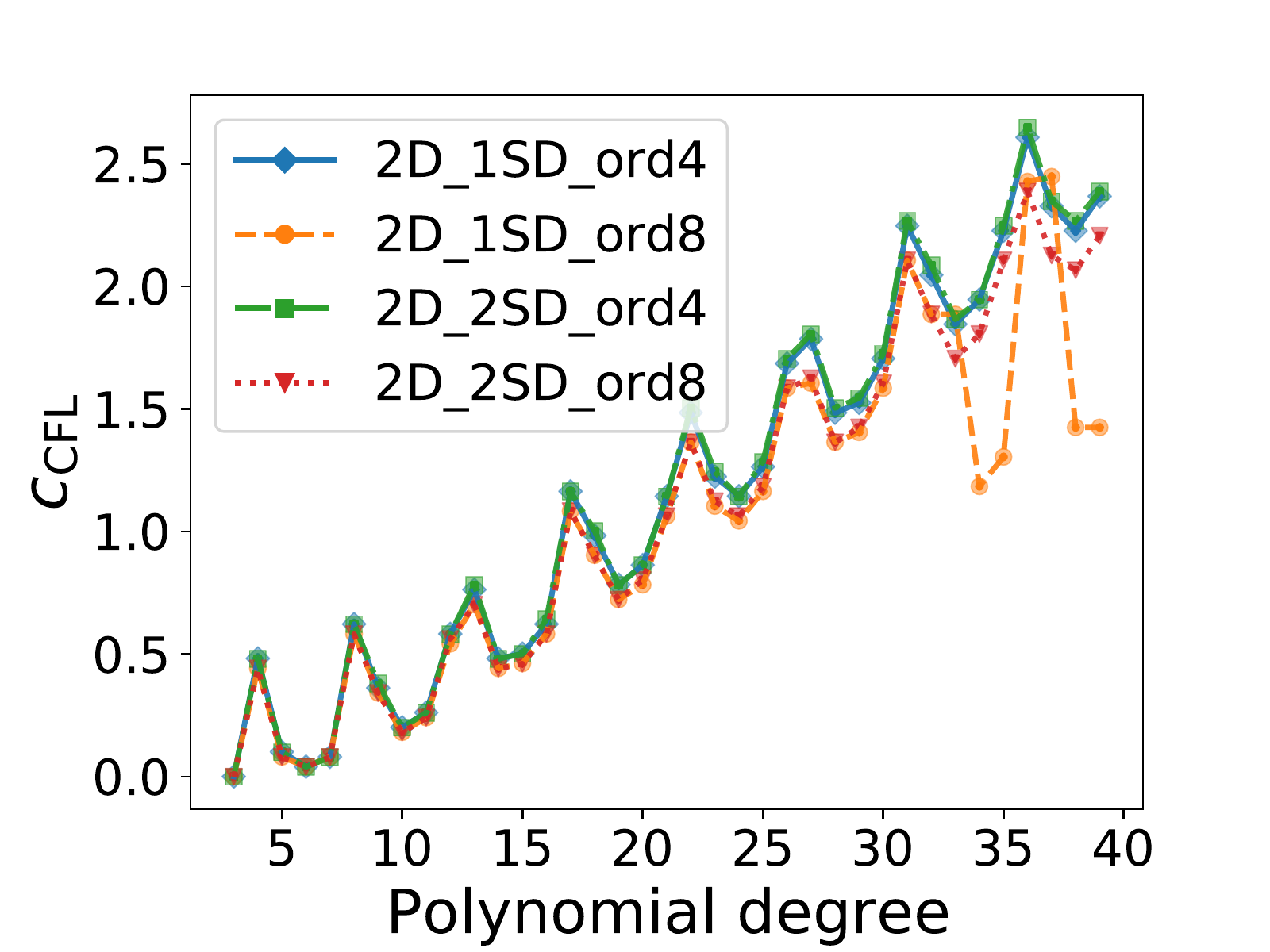}\hspace{-0.3cm}
	    &
	    \includegraphics[scale=0.24]{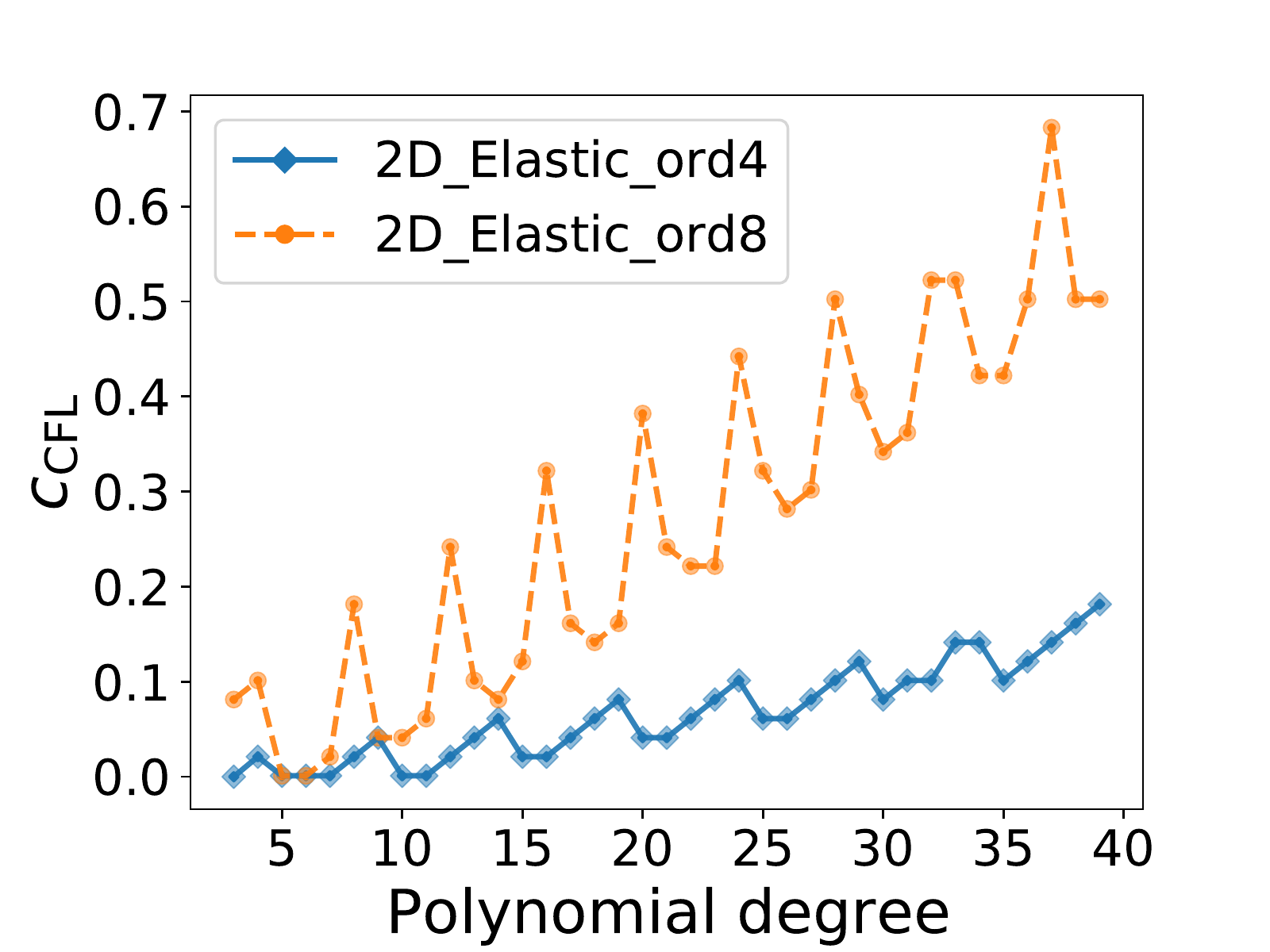}\hspace{-0.3cm}
	    \\
	    \includegraphics[scale=0.24]{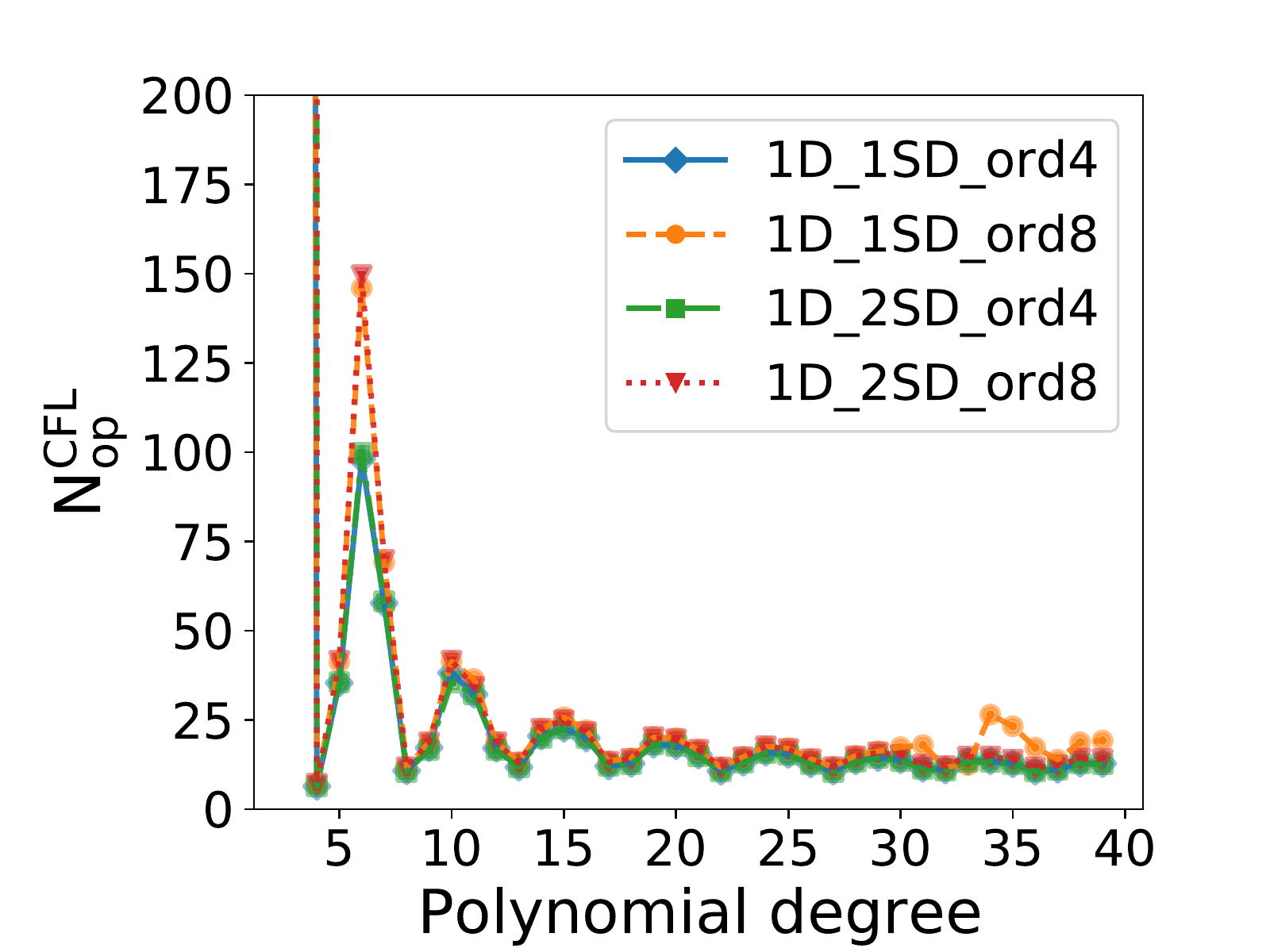}\hspace{-0.3cm}
	    &
	    \includegraphics[scale=0.24]{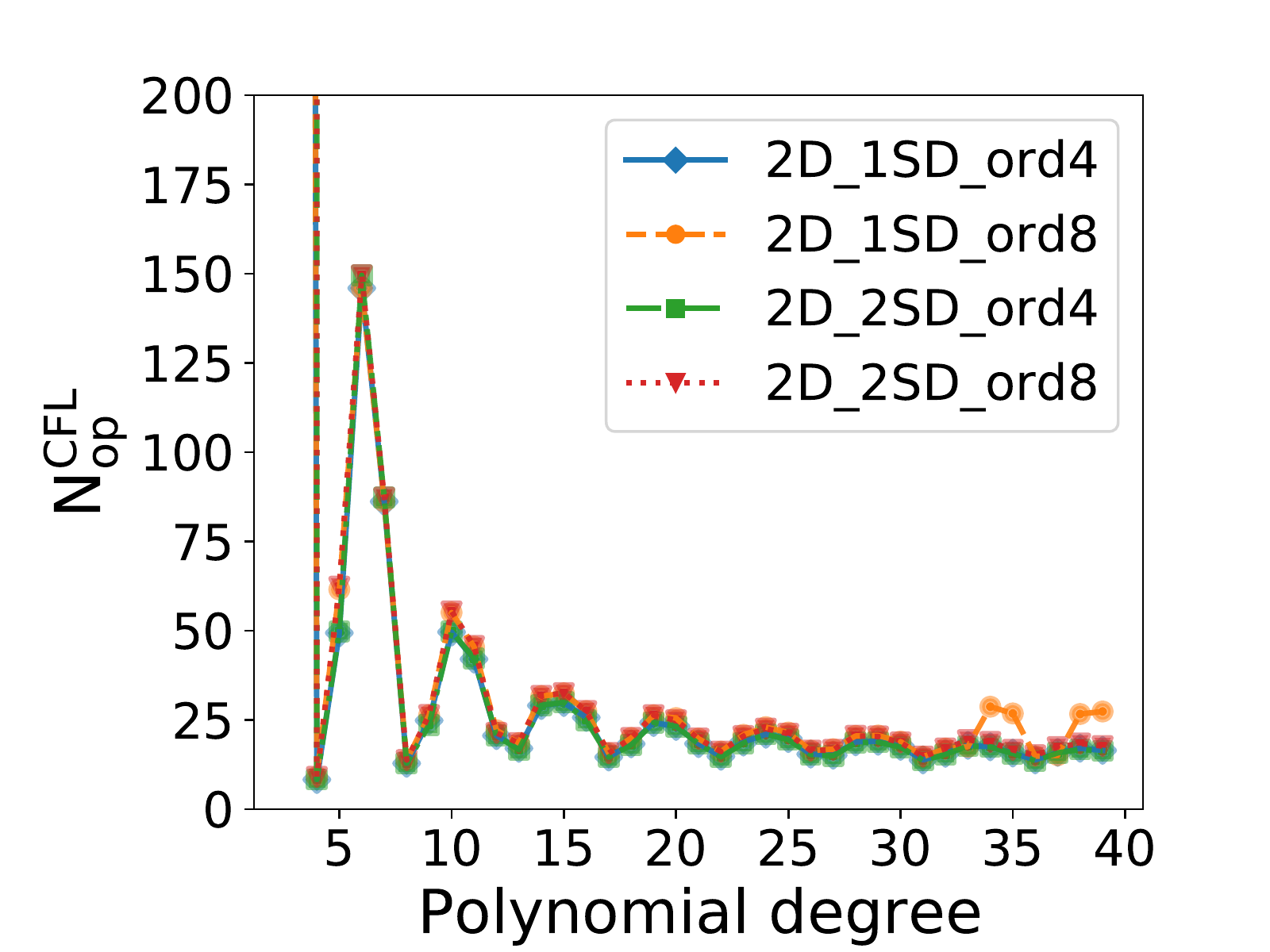}\hspace{-0.3cm}
	    &
	    \includegraphics[scale=0.24]{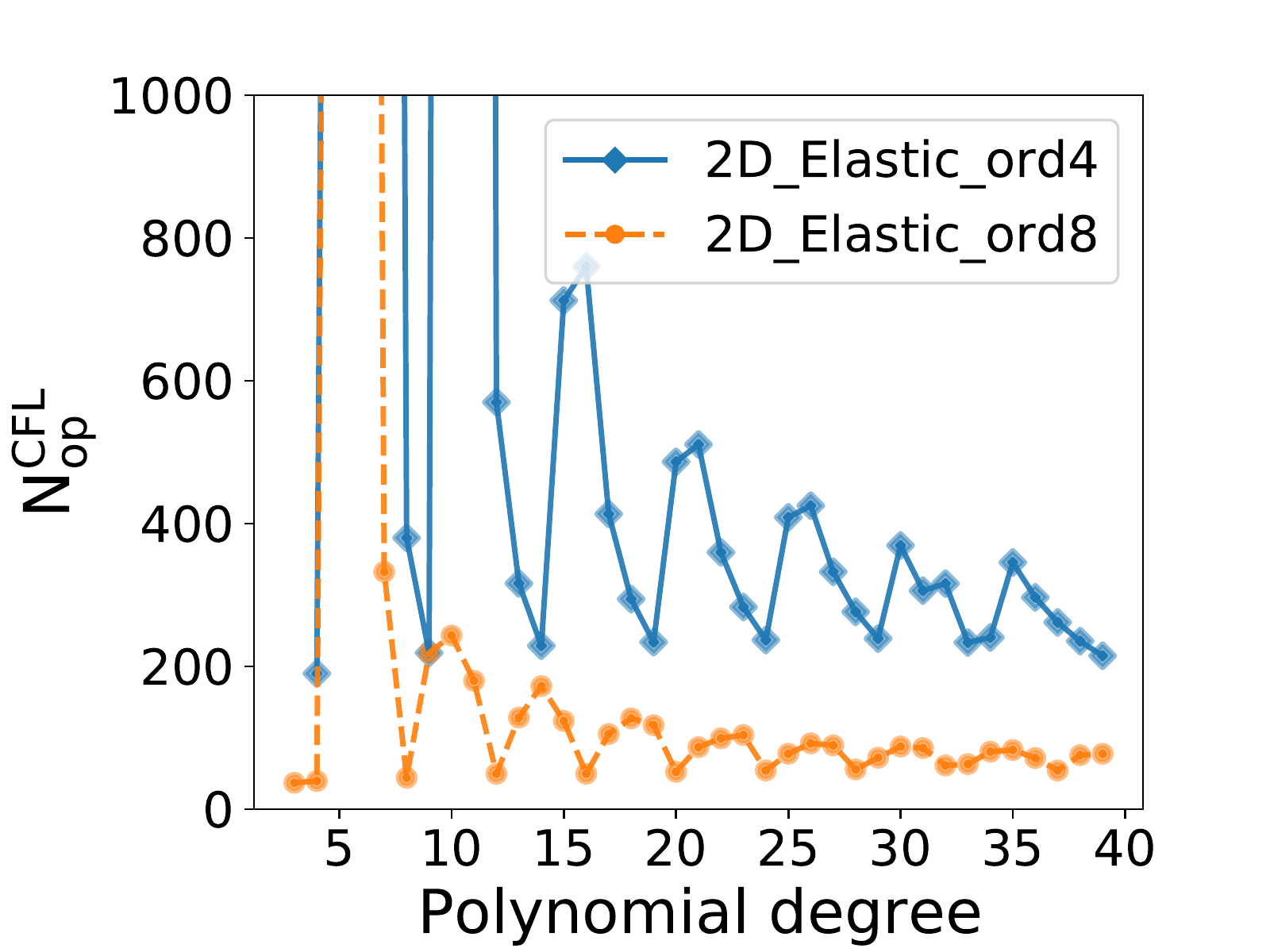}\hspace{-0.3cm}
    \end{tabular}
	\caption{Stability graphics of the Faber approximation method for different spatial discretizations, dimensions, equations formulations, and polynomial degrees $m=\{3,4,...,40\}$. The CFL number (top row), and the measure of the operations number of Eq. \eqref{eq_efficiency_s} (bottom row). Higher polynomial degrees implies in larger $c_{\text{\tiny CFL}}$, but this does not necessarily results in a decrease of the computations.		\label{fig_stab_eff_all}}
\end{figure}

Next, we compute the CFL number ($c_{\text{\tiny CFL}}$) as the largest $\alpha$ so that the spectral radius $\rho(\boldsymbol{G})$ is at most $1+\epsilon$, with $\epsilon=10^{-7}$ accounting for round-off errors. We compute this CFL number for several polynomial degrees and different numerical specifications with results shown in the upper row of Fig.\,\ref{fig_stab_eff_all}.
We observe that for the acoustic equations, the spatial discretization and equation formulations considered have little influence on the stability of the method. Yet, for the elastic equations this is not the case, and an important gain is observed when the spatial discretization order is increased. Overall, there is a small improvement on the performance using 2SD equations instead of 1SD, and the passing from 1D to 2D reduces the values of $c_{\text{\tiny CFL}}$. This reduction is in part because, in two dimensions, the CFL number is divided by a factor of square-root of two. Moreover, in all scenarios, the CFL number is enhanced with the increase of the polynomials degrees. 

Furthermore, we use the MVOs as reference for computational cost and define
\begin{equation}
    \text{N}^{\text{\tiny CFL}}_{\text{op}}=\frac{\text{\# MVOs}}{c_{\text{\tiny CFL}}},\label{eq_efficiency_s}
\end{equation}
representing the ratio between the computational requirements and a value relating to the CFL.
Hence, this scalar value represents the computational efficiency, where smaller values relate to better efficiency.
Since the number of MVOs coincides with the polynomial degree used in the approximation, for the graphics we will use polynomial degrees instead of MVOs with results given in Fig.\,\ref{fig_stab_eff_all}.
The oscillatory behavior indicates that higher polynomial degrees are not reflected by less computations, but rather some particular degrees are more fitted to improve the values of $\text{N}^{\text{\tiny CFL}}_{\text{op}}$.
The oscillations of the curves have a periodicity of a length of four and five degrees, which we account for by relating it to some sort of symmetry of the complex polynomials $\boldsymbol{F}_j$.

\subsection{Numerical dispersion}\label{sec_disp}

We continue studying the dispersion ($R$) given by the quotient between the velocity of the numerical solution $c_{\text{num}}$ and the real wave velocity $c$, hence
\begin{equation}\label{eq_dispersion}
    R=\frac{c_{\text{num}}}{c}=\frac{w_{\text{num}}}{k\,c},
\end{equation}
where $w_{\text{num}}$ is the numerical angular frequency, and $k$ is the wavenumber. Ideally, a method with no spurious numerical dispersion should have $R=1$, therefore, we define the dispersion error as $|R-1|$.

In Equation \eqref{eq_dispersion}, the velocity $c$ is known, $k$ is any wavenumber value contained in the interval $[0,\pi/\Delta x]$, and $w_{\text{num}}$ is calculated from the phase of the eigenvalues of $\boldsymbol{G}$.

\begin{figure}[tbh]\hspace{-0.5cm}
	\begin{tabular}{c|c|c}
		(a) One dimension, acoustic.
		&
		(b) Two dimensions, acoustic.
		&
		(c) Two dimensions, elastic.
		\\
	    \includegraphics[scale=0.24]{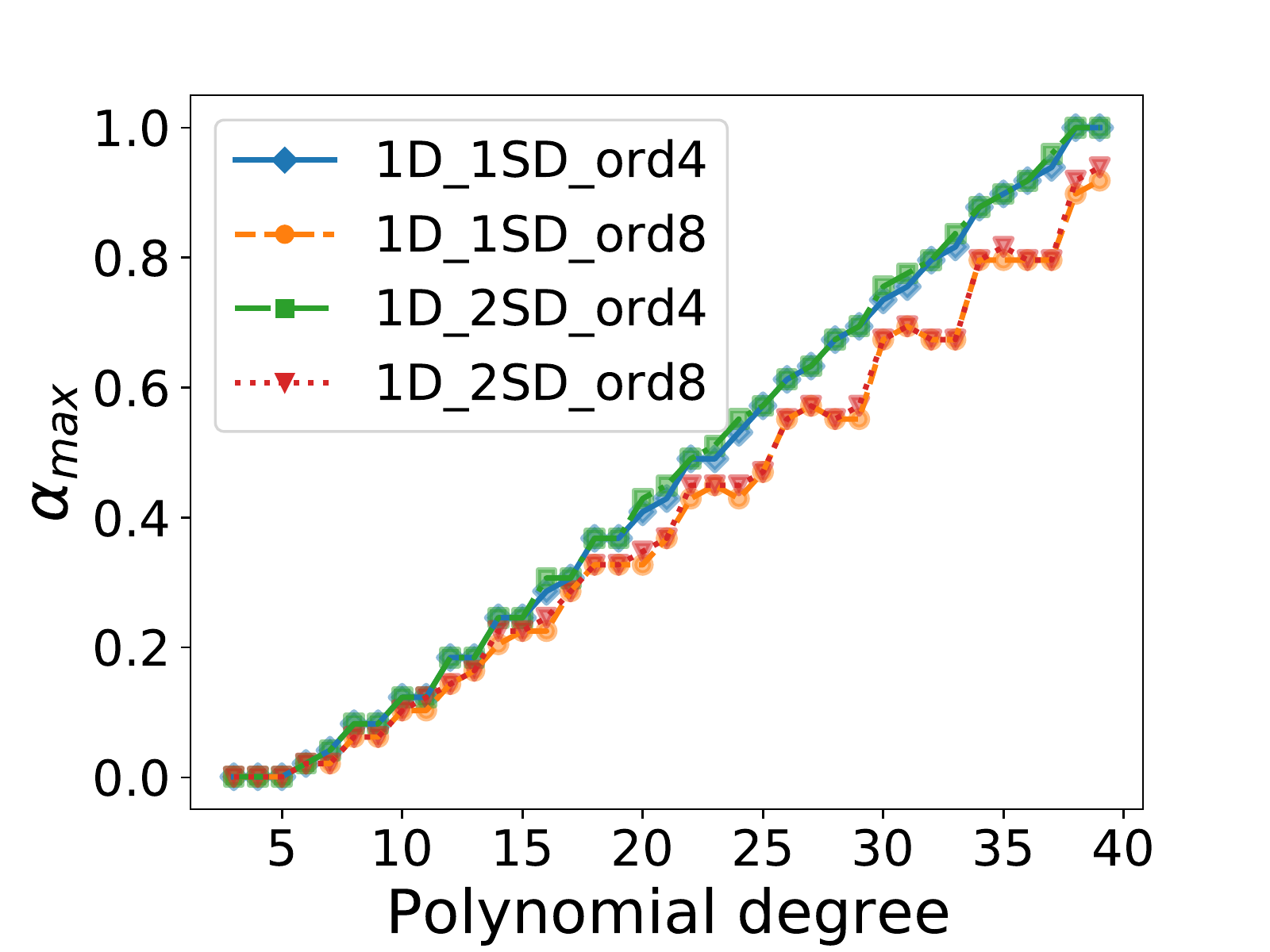}\hspace{-0.3cm}
	 	&
	    \includegraphics[scale=0.24]{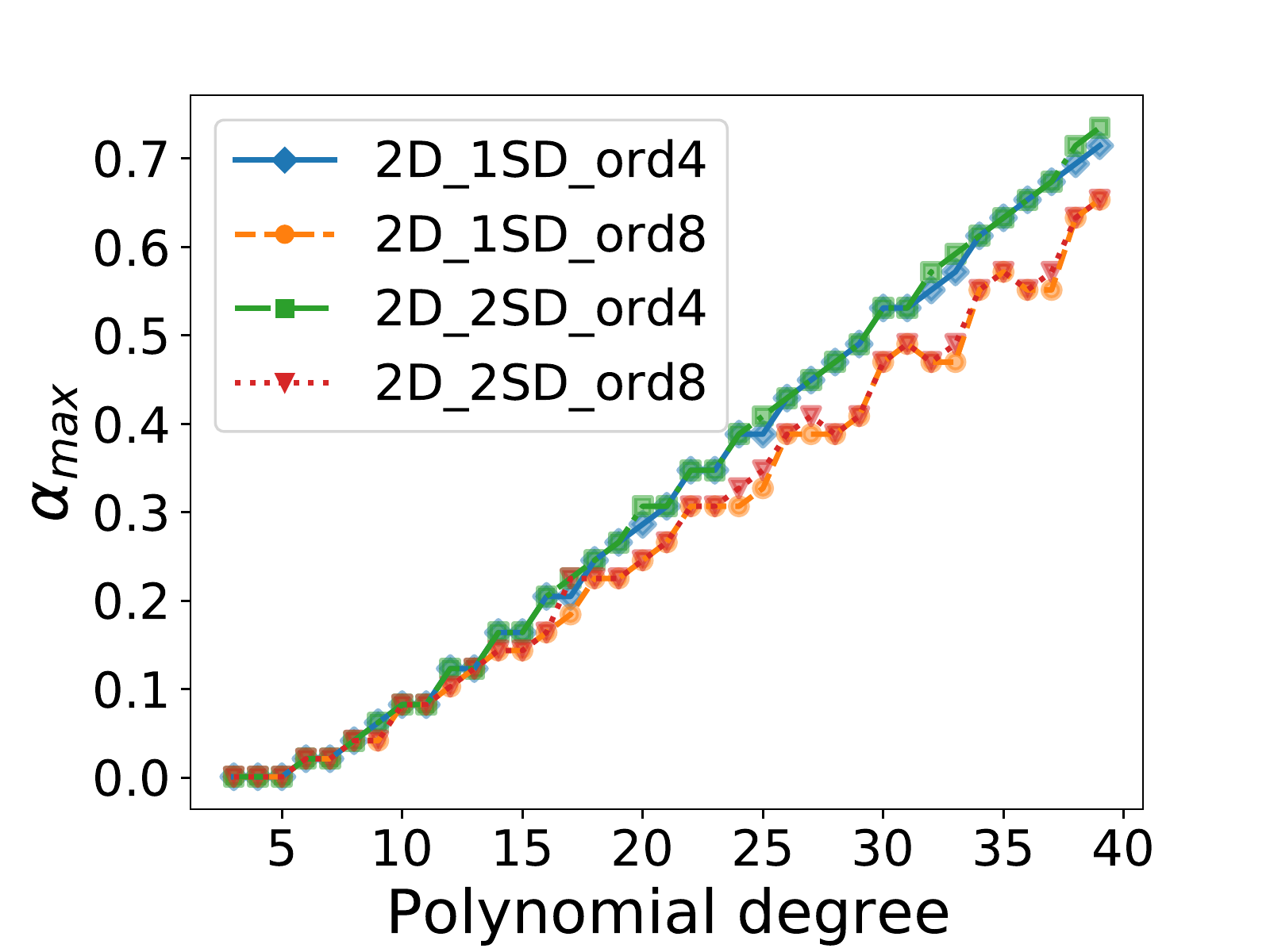}\hspace{-0.3cm}
	    &
	    \includegraphics[scale=0.24]{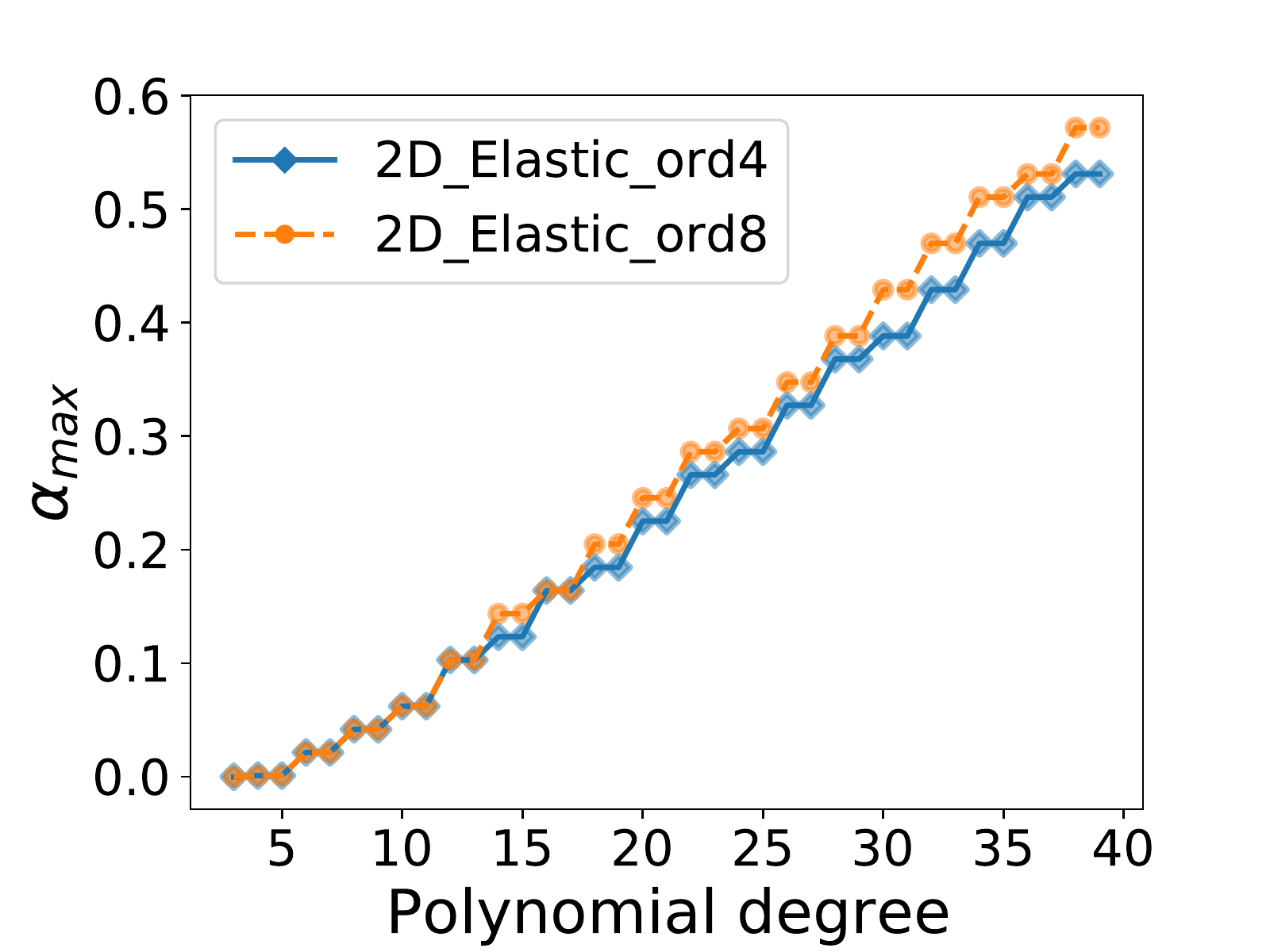}\hspace{-0.3cm}
	    \\
	    \includegraphics[scale=0.24]{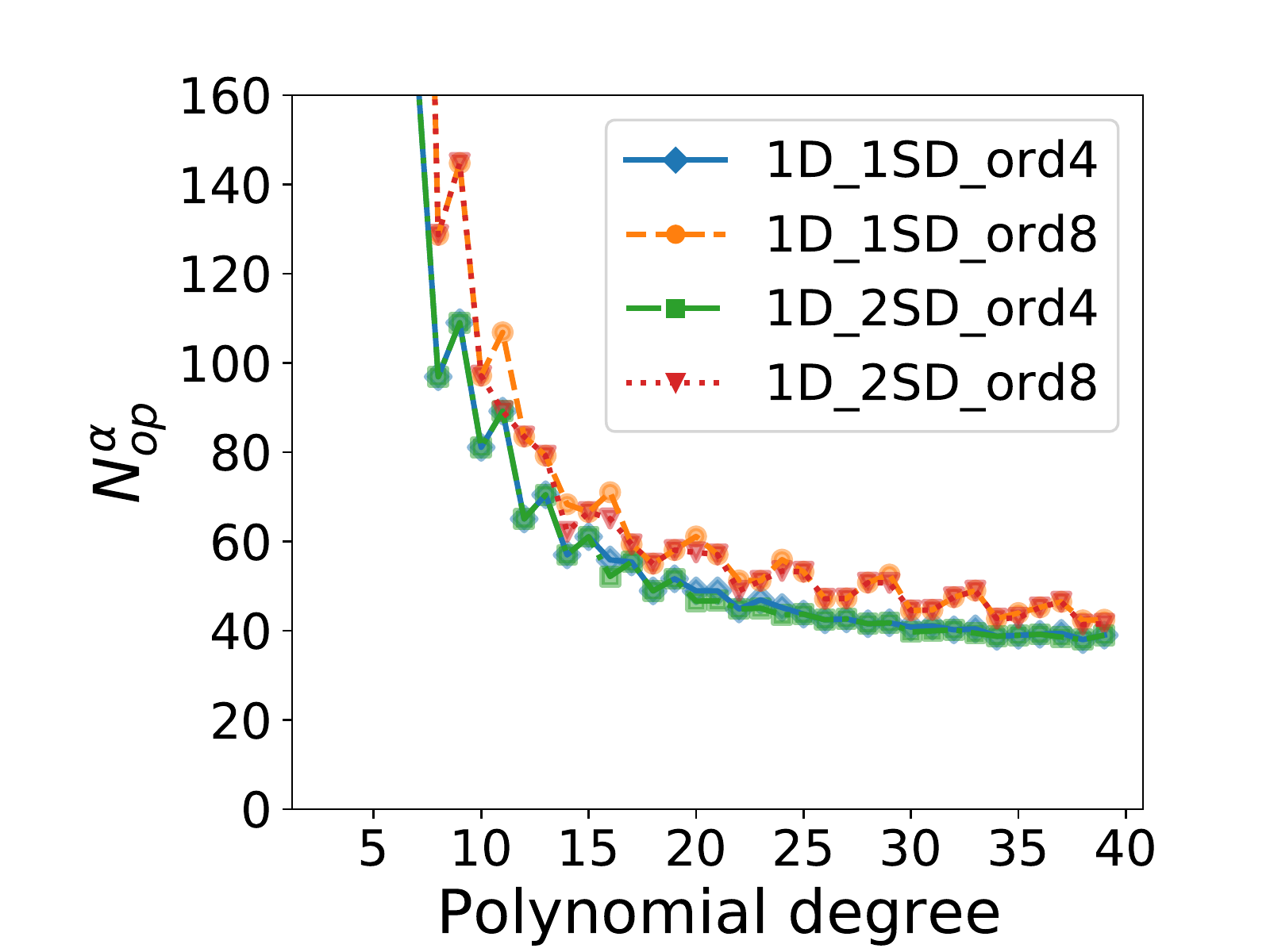}\hspace{-0.3cm}
	    &
	    \includegraphics[scale=0.24]{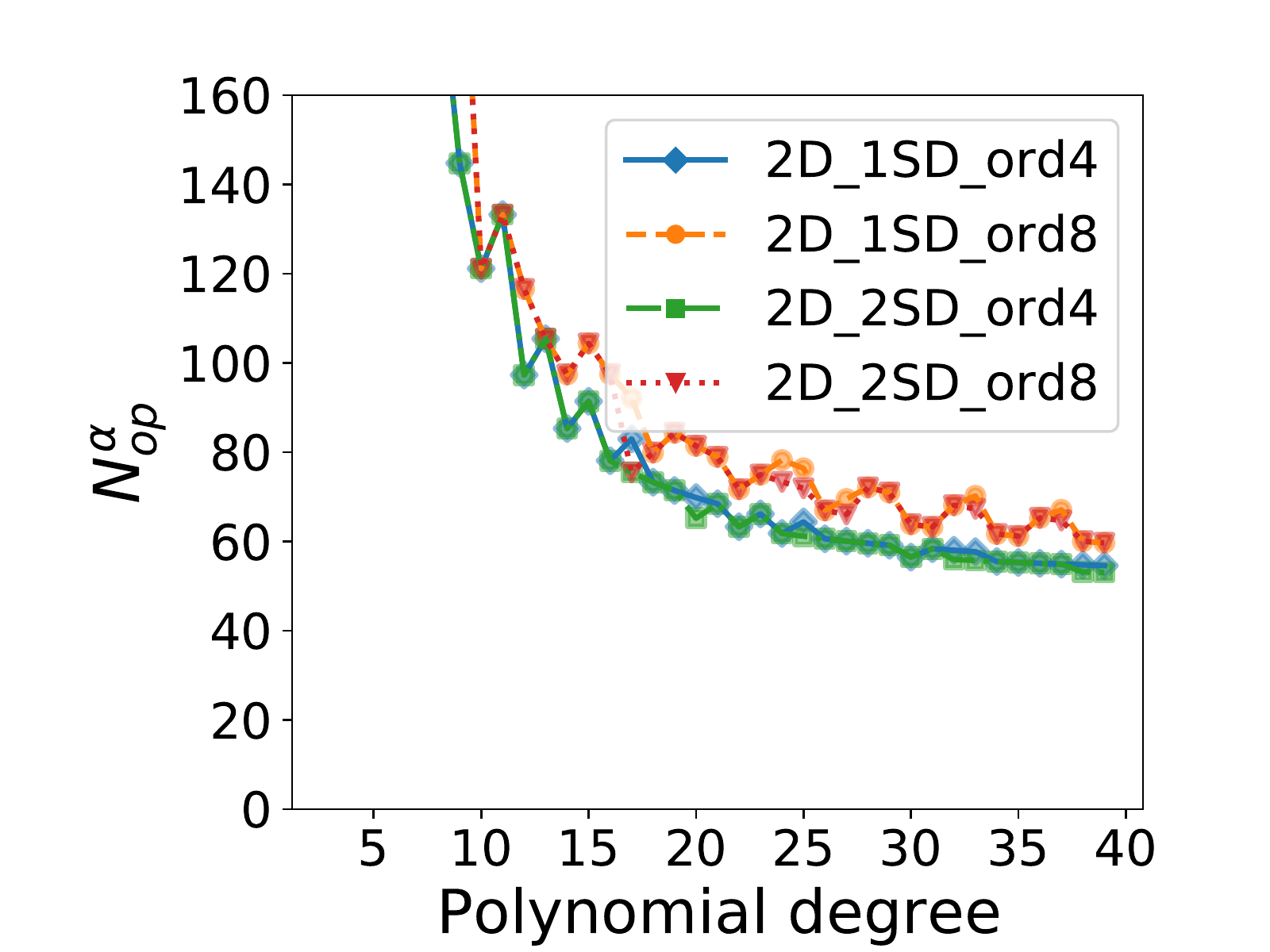}\hspace{-0.3cm}
	    &
	    \includegraphics[scale=0.24]{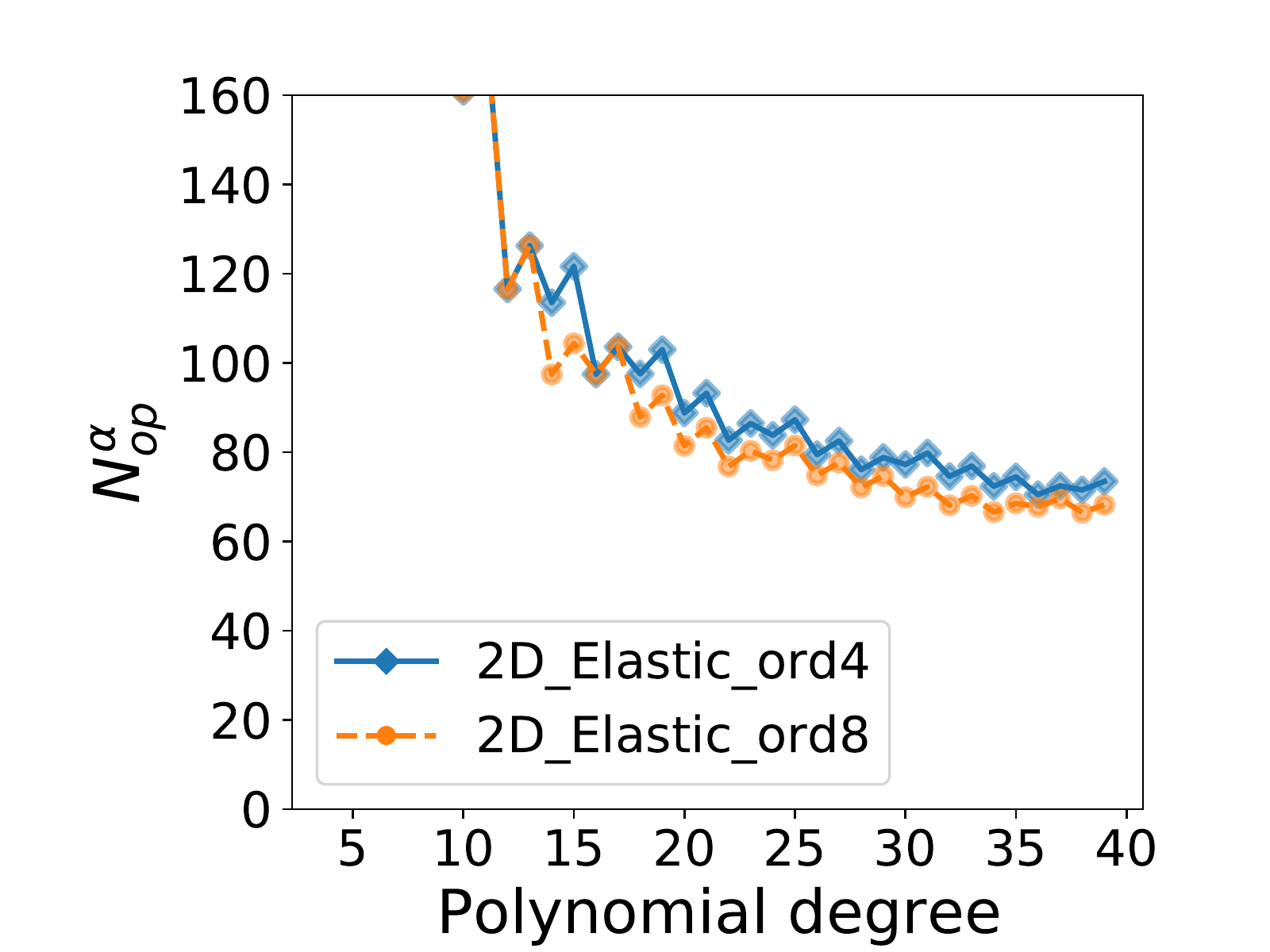}\hspace{-0.3cm}
    \end{tabular}
	\caption{Dispersion studies of Faber approximation method using different spatial discretizations, dimensions, equations formulations, and polynomial degrees $m=\{3,4,...,40\}$. Higher polynomial degrees implies larger $\alpha_{R}$ and a non-monotonous decrease of the operations number.
	\label{fig_disp_eff_all}}
\end{figure}

Now, we introduce $\alpha_{R}$, similar to the $c_{\text{CFL}}$ number, but relating to a maximum dispersion error instead.
We define $\alpha_{R}$ to be the maximum $\alpha$ (see Eq.\,\eqref{eq_alpha}) such that the dispersion error is less than $\epsilon_R=10^{-5}$, which is set to be the required dispersion accuracy.
We can then compute the computational efficiency with respect to dispersion by using
\begin{equation}
    \text{N}^{\alpha}_{\text{op}}=\frac{\text{\# MVOs}}{\alpha_R}.\label{eq_efficiency_d}
\end{equation}

Results are presented in Fig.\,\ref{fig_disp_eff_all}, using 1SD and 2SD formulations, in one and two dimensions, with different spatial discretization orders.

From Fig.\,\ref{fig_disp_eff_all}, we observe that the $\alpha_{R}$ in 1D is larger than in the 2D case, and the dispersion changes very little with respect to the spatial discretization order.
For the elastic equations, $\alpha_{R}$ decreases even more, and there is only a small difference between the spatial discretization order. Moreover, in agreement with the stability analysis, when increasing the polynomial degree, this leads to larger $\alpha_{R}$.
However, in contrast with the results on stability, there is a stronger impact of the equation formulation, and higher polynomial degrees implies fewer computations. This is in agreement with the expected behavior, since larger polynomial degrees allows larger time-steps, diminishing the number of time-steps needed for the computations and, then, reducing the dispersion error.

\section{Convergence and efficiency}\label{sec_convergence_efficiency}

In this section, we address the numerical convergence and computational efficiency of Faber polynomial approximations for the full equation sets, on the limited area domain with the PML absorbing conditions. Due to the complexity of the equations, the analysis is purely numerical, relying on the test cases shown in Section \ref{sec_test_cases}. To ensure robustness of the results, and as explained in Section \ref{sec_test_cases}, the seven experiments vary in levels of complexity and problem specifications. Then, we assess the convergence on the experiments, by means of the approximation error in $L_2$, which is computed for a wide range of Faber polynomial degrees and time-step sizes $\Delta t$.

To approximate the solutions using Faber polynomials, we define a spatial step size $\Delta x=0.0025$ for the 1D examples as well as $\Delta x=0.02$ for the 2D examples and use a finite difference scheme with 4th and 8th spatial order.
For comparison purposes, we use a nine-stages seventh-order temporal Runge-Kutta scheme RK(9,7) recommended for hyperbolic problems (see \citet{calvo1996explicit})) with a small time-step size ($\Delta t=\Delta x/(8c_\text{max}))$, where $c_\text{max}$ is the maximum velocity.
The equation formulation and spatial discretization size and order used in the RK(9,7) are the same as the one adopted in the Faber approximations, so that the spatial operator is exactly the same as the one used in the Faber approximation scheme. Therefore, only temporal effects are visible in the numerical errors shown.  

\begin{figure}[tbh]\hspace{-0.1cm}
    \subfloat[TC\#2, in 1D, acoustic,$\quad$\linebreak 1SD, spatial order 4.]{\includegraphics[scale=0.25]{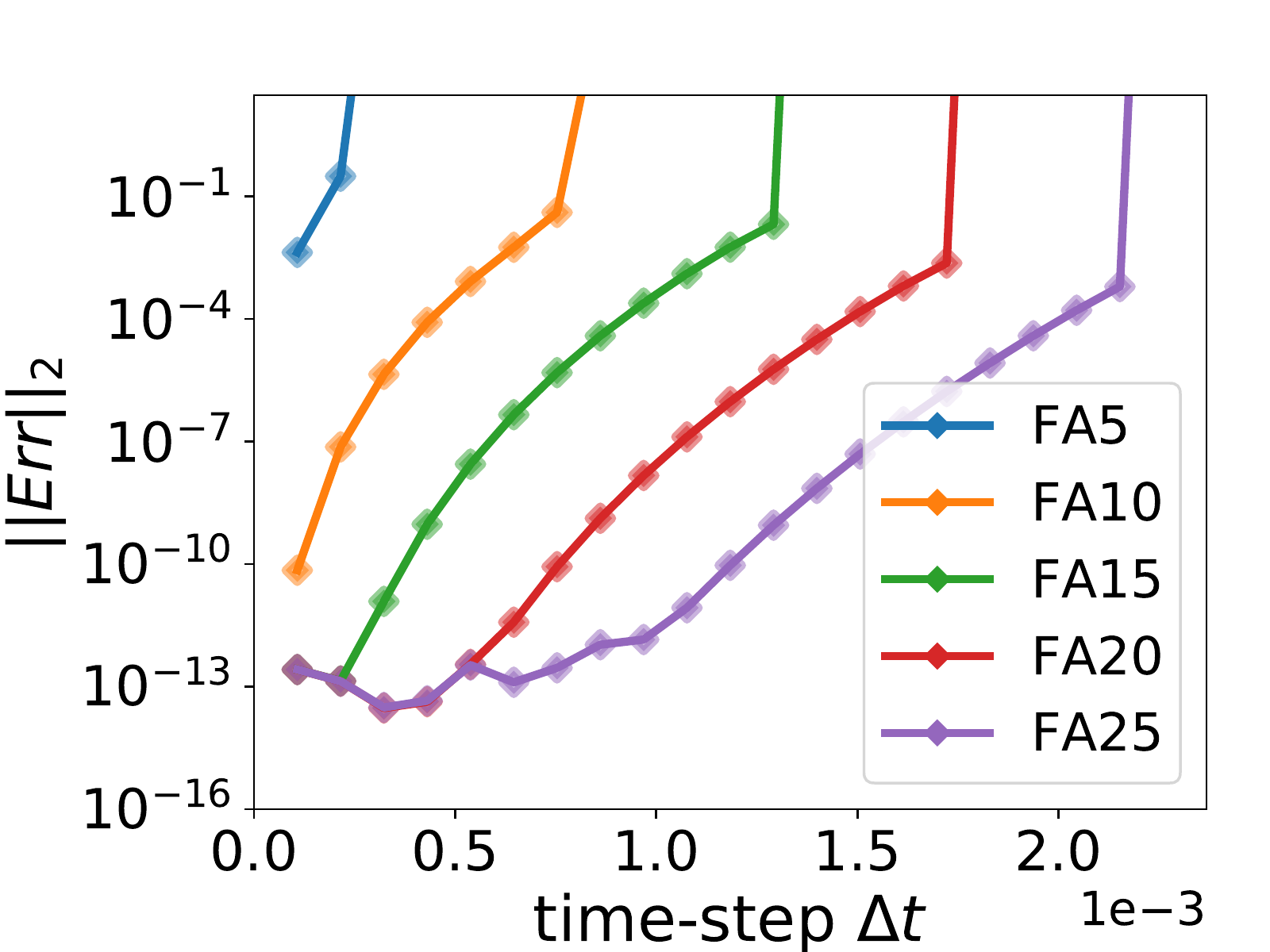}}
    \subfloat[TC\#6, in 2D, acoustic,$\quad$\linebreak acoustic 2SD, spatial order 4.]{\includegraphics[scale=0.25]{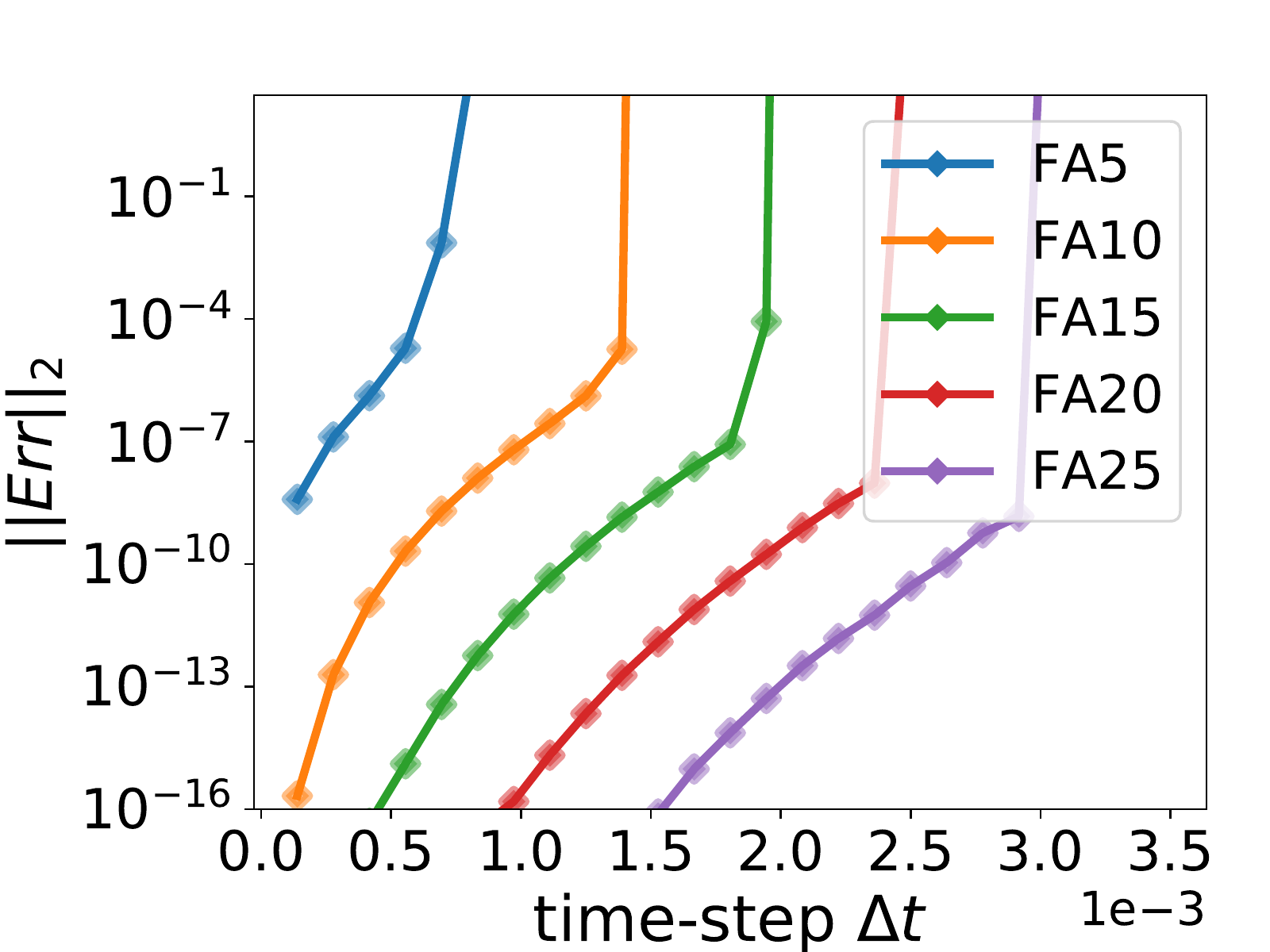}}
    \subfloat[TC\#7, in 2D, elastic,$\quad$\linebreak spatial order 8.]{\includegraphics[scale=0.25]{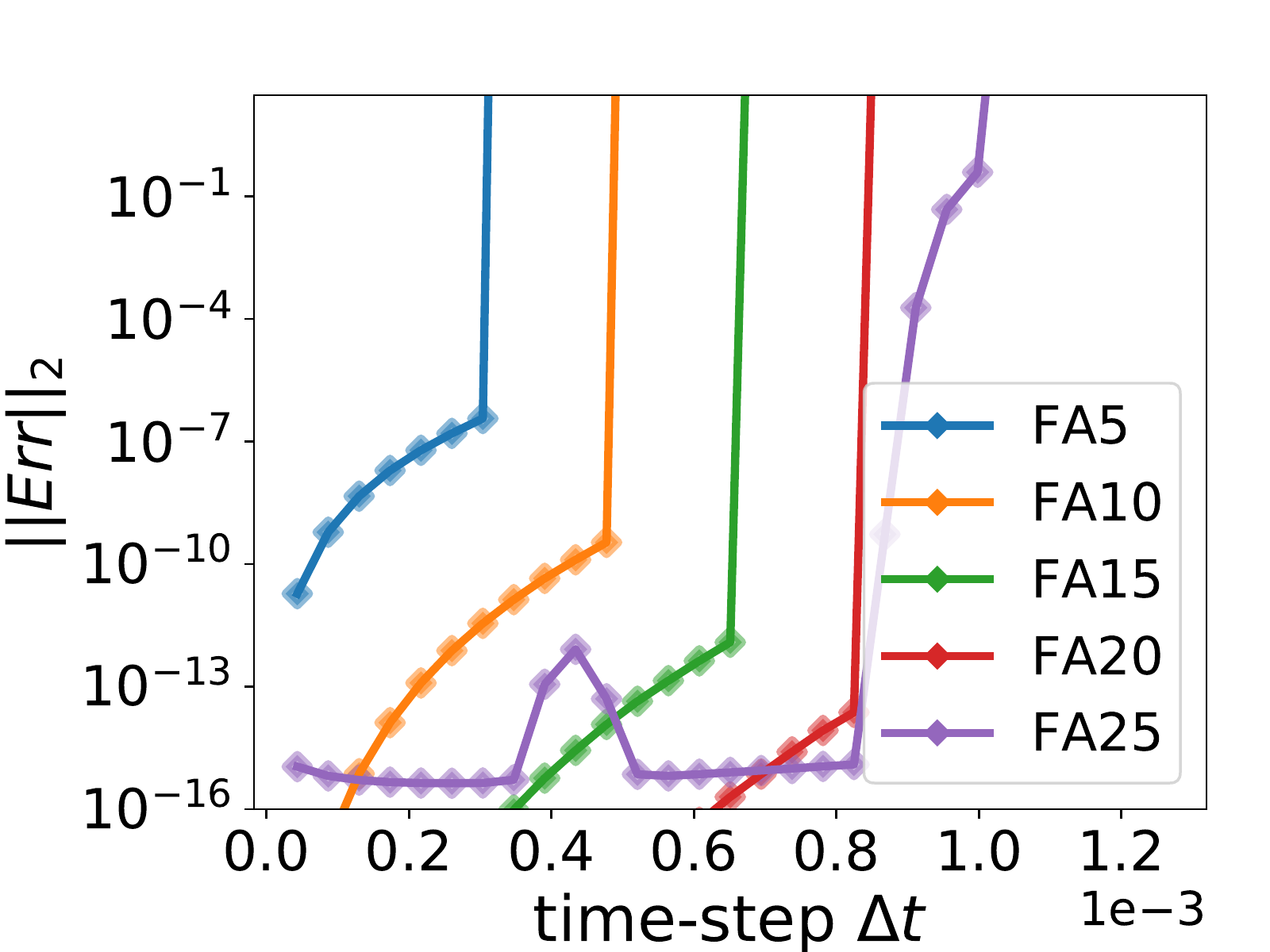}}
\caption{Approximation error of Faber polynomials using 1SD, 2SD and elastic formulations, to solve TC\#1, TC\#6 and TC\#7, respectively. The curves represent the error when using polynomial of degrees $m=\{5,10,\dots,25\}$. Increasing the degree of the polynomial allows for larger steps in time, while keeping the approximation error smaller than a fixed threshold.}\label{fig_approx_test_167}
\end{figure}

Fig.\,\ref{fig_approx_test_167} shows the approximation error for three different examples, with one and two dimensions, using three different formulations and for different spatial approximations order.
We observe in all cases that increasing the polynomial degree of the method allows larger stable time-step sizes.
If we define a maximum $\Delta t$ such that the approximation error is bounded by a fixed threshold, we note that the magnitude of this $\Delta t_\text{max}$ changes depending on the problem specifications.
We also notice that for each polynomial degree the convergence deteriorates before reaching the critical $\Delta t$.

Moreover, we point out that the solution behavior described in Fig.\,\ref{fig_approx_test_167} is sustained for other scenarios, and is independent of the wave formulation, spatial discretization, and numerical examples, considered in this paper.
This is due to the fact that larger polynomial degrees relate to higher approximation orders, hence allowing larger time-steps.

\begin{figure}[tbh]
	\begin{tabular}{l|l|l}
		(a) TC\#1 and TC\#3,
		&
		(b) TC\#4 and TC\#5,
		&
		(c) TC\#7, spatial
		\\
		\hspace{0.5cm}spatial order 8 & \hspace{0.55cm}spatial order 4 & \hspace{0.45cm} order 4 \& 8\\
	    \includegraphics[scale=0.24]{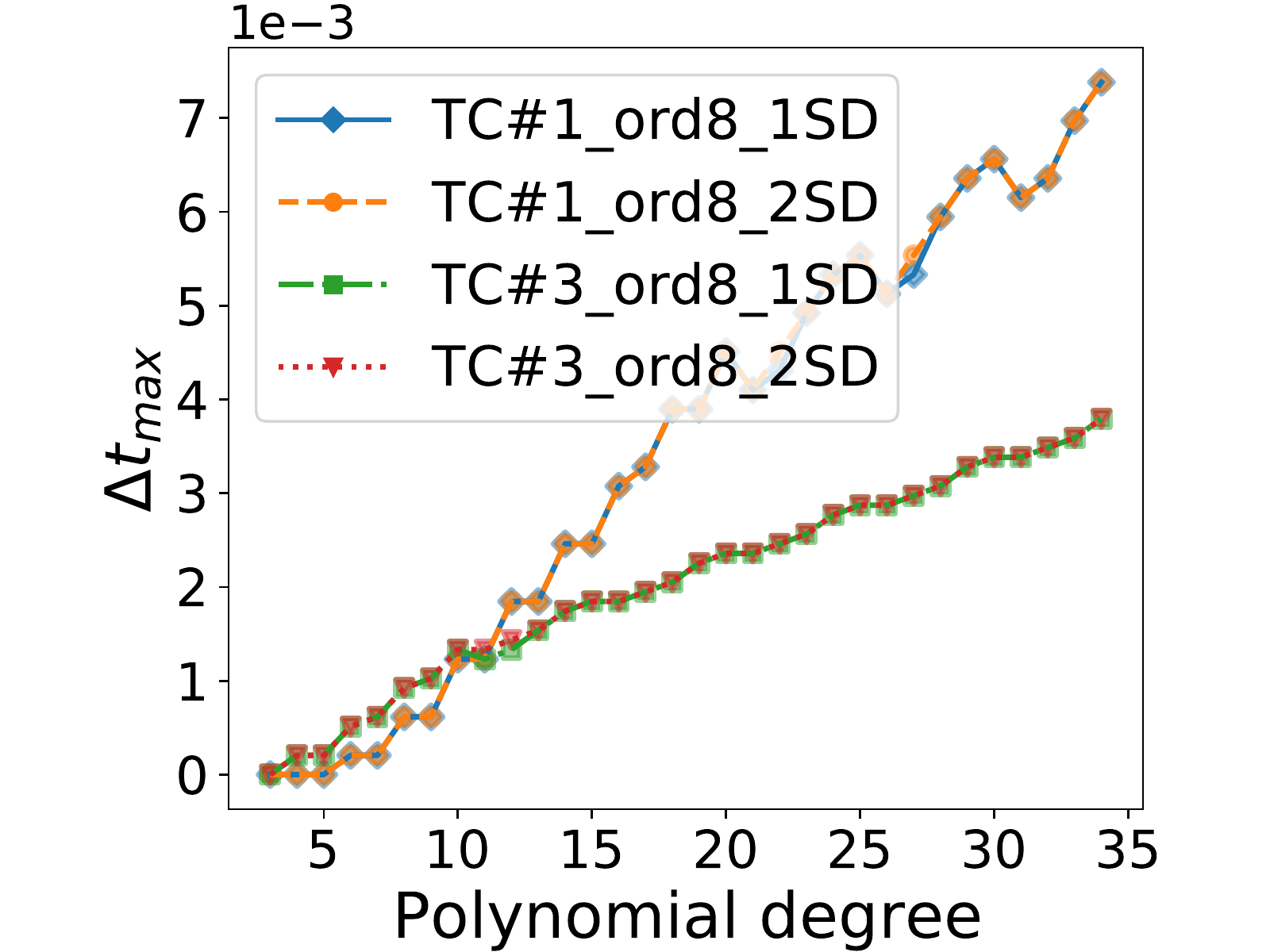}\hspace{-0.3cm}
	 	&
	    \includegraphics[scale=0.24]{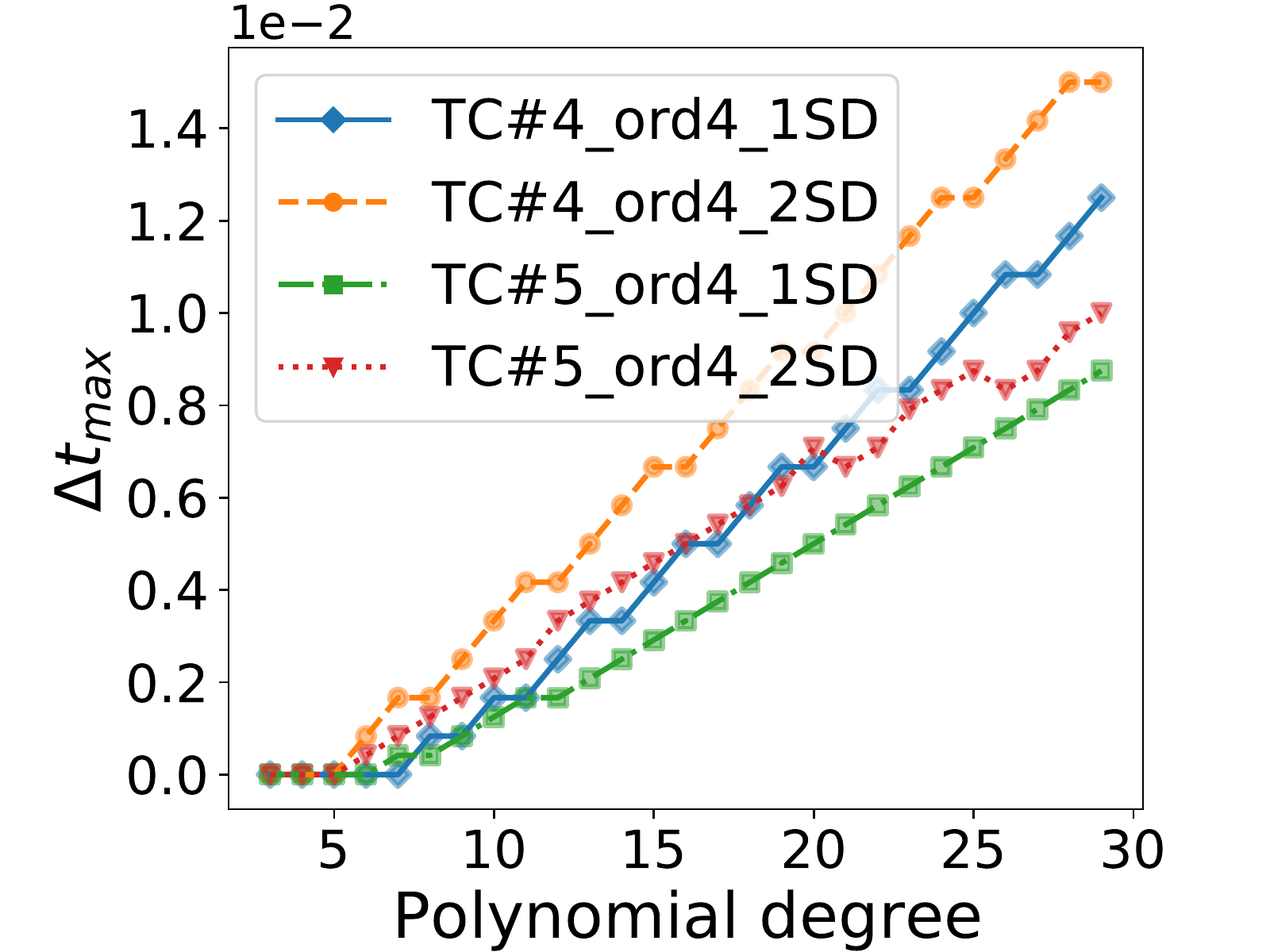}\hspace{-0.3cm}
	    &
	    \includegraphics[scale=0.24]{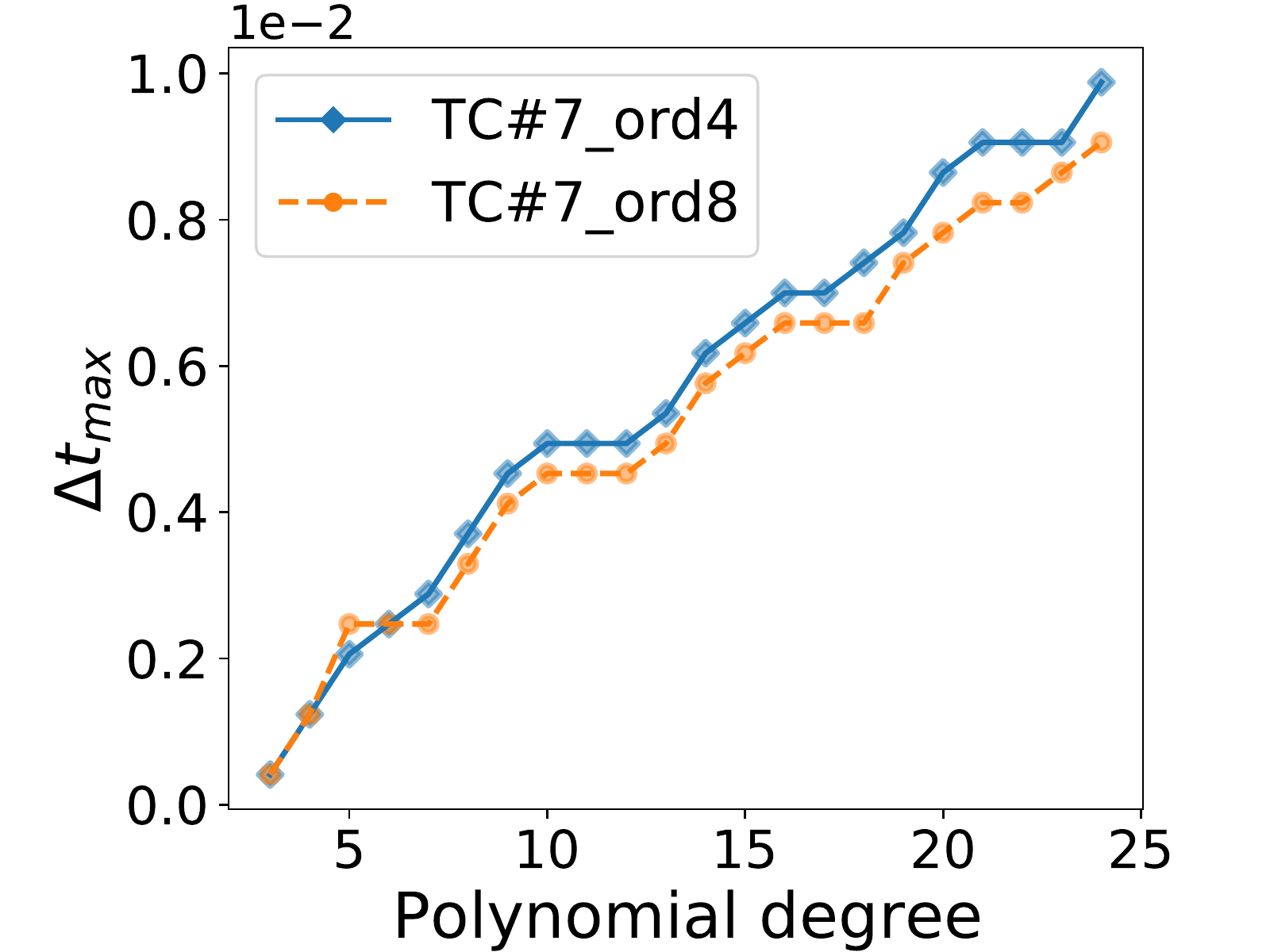}\hspace{-0.3cm}
	    \\
	    \includegraphics[scale=0.24]{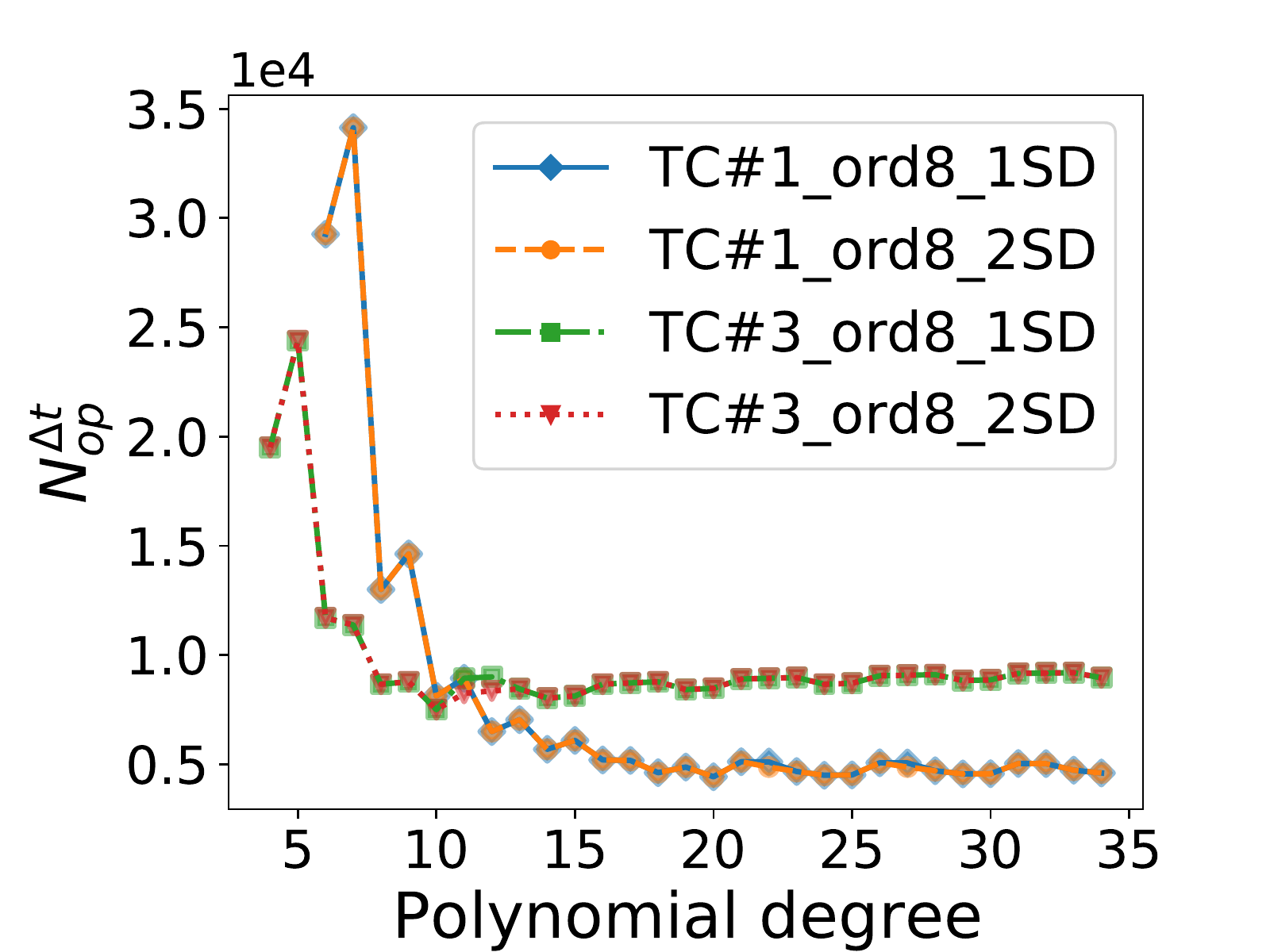}\hspace{-0.3cm}
	    &
	    \includegraphics[scale=0.24]{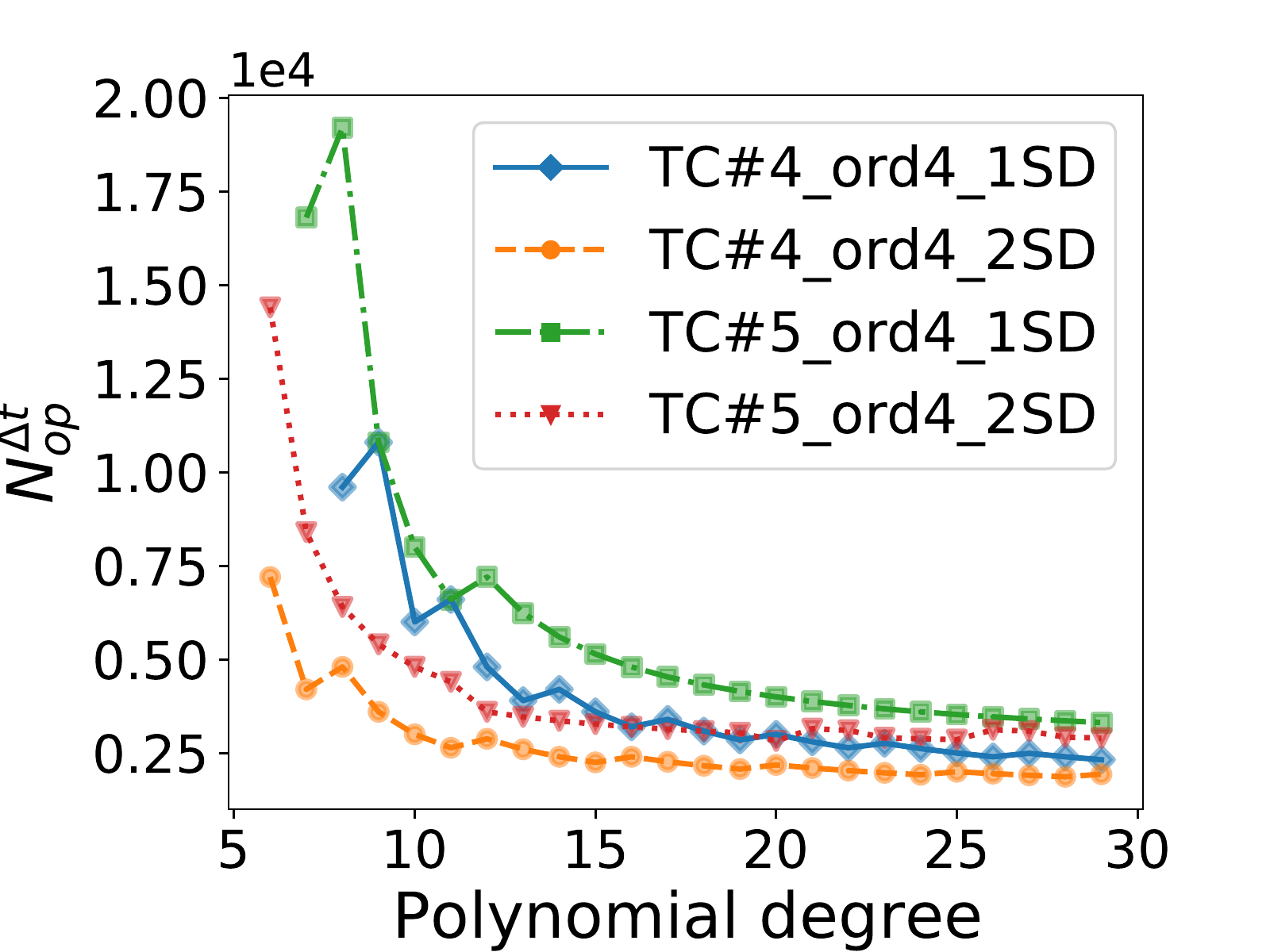}\hspace{-0.3cm}
	    &
	    \includegraphics[scale=0.24]{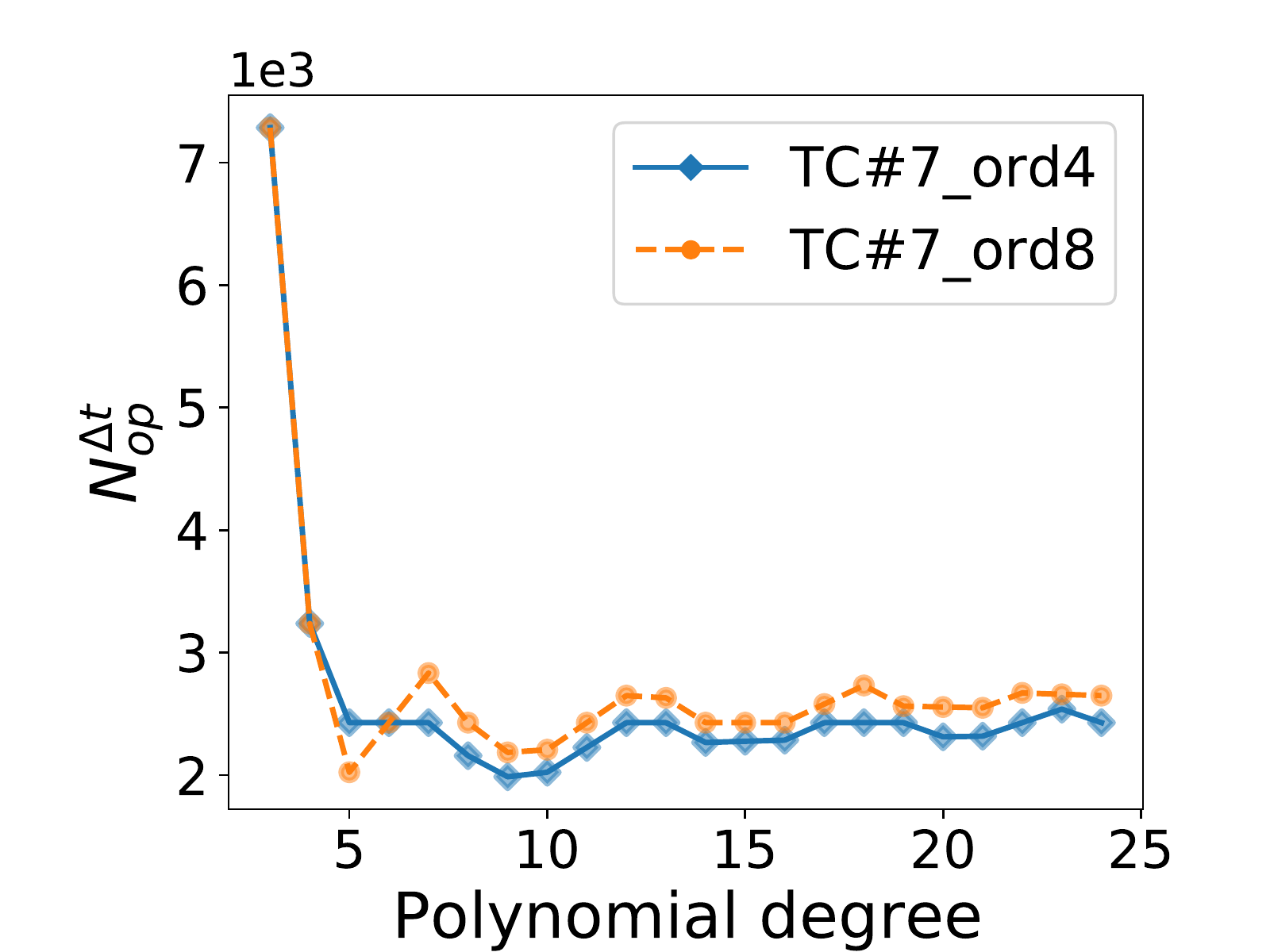}\hspace{-0.3cm}
    \end{tabular}
	\caption{ Convergence in polynomial order for 1SD, 2SD, and elastic formulations, using different experimental sets, spatial discretization orders, and a wide range of polynomial degrees.
	The maximum $\Delta t$ such that the error of Faber approximations is less than $10^{-6}$ is shown in the top row. In the bottom row, we show the number of operations using the values of $\Delta t_\text{max}$ of the upper line. When the polynomial degree increases, the maximum allowed time-step size also increases, together with a decrease of the number of operations.
	\label{fig_conv_eff_all}}
\end{figure}

Next, we further investigate $\Delta t_\text{max}$, the maximum time-step size allowed by the polynomials while maintaining an error lower than $\epsilon_{\Delta t}=10^{-6}$ (see Fig.\,\ref{fig_conv_eff_all}). In (a) we notice that the medium velocity influences $\Delta t_\text{max}$, as expected, since the maximum velocity in TC\#1 is two times lower than the maximum velocity in TC\#3. From the same figure, we note also that even when there is a consistent increase of maximum time-step size, the behavior is not monotonic (TC\#1 lines in (a)). Moreover, there are experiments where the equation formulation seems to have no influence in the convergence (SubFig.\,(a)), but there are others where the 2SD formulation seems to perform better than the 1SD (SubFig.\,(b)).
We account for that by an increased difficulty of the 2D experiments TC\#4 and TC\#5, when compared to the 1D TC\#1 and TC\#3, revealing the differences between formulations 1SD and 2SD, but with further investigation required for confirmation of this hypothesis.
In general, we have observed that 2SD always perform similarly or better than 1SD.
This suggests the continuum formulation of the wave equation with PML to be an important factor to consider when solving the equations. 

In addition, from Fig.\,\ref{fig_conv_eff_all}(c), we see that the spatial discretization order has little influence for the elastic equations.
We again remark that our experiments consider the reference solution to have the same discrete operator as the exponential scheme, therefore only temporal effects are to be noted.
Here, we see that the increase in spatial order has only a minor effect of reducing the maximum time step size. 

We now define a measure of the number of operations depending on the time step sizes (analogously to Eq.\,\eqref{eq_efficiency_s} and Eq.\,\eqref{eq_efficiency_d}) as 
\begin{equation*}
    \text{N}^{\Delta t}_{\text{op}}=\frac{\text{\# MVOs}}{\Delta t_{\text{max}}},
\end{equation*}
where again $\Delta t_{max}$ is the maximum $\Delta t$ such that the approximation error is bounded by a fixed threshold $\gamma=10^{-6}$.

In Sections \ref{sec_stab} and \ref{sec_disp} we used  simplified formulations on periodic domains. Here, we finally consider the full equation sets on the limited area domain with the PML absorbing conditions.
Therefore, $N^{\Delta t}_{\text{op}}$ represents a realistic measure of the amount of computations by unit of time.

A general behavior in the number of operations graphics in Fig.\,\ref{fig_conv_eff_all} is to have a lot of computations for low degree approximations, followed by a pattern of declination, and seems approximating to an equilibrium.

\subsection*{The corner model}

\begin{figure}[H]
\subfloat[High velocity contrast$\quad$\linebreak heterogeneous medium.]{\includegraphics[scale=0.24]{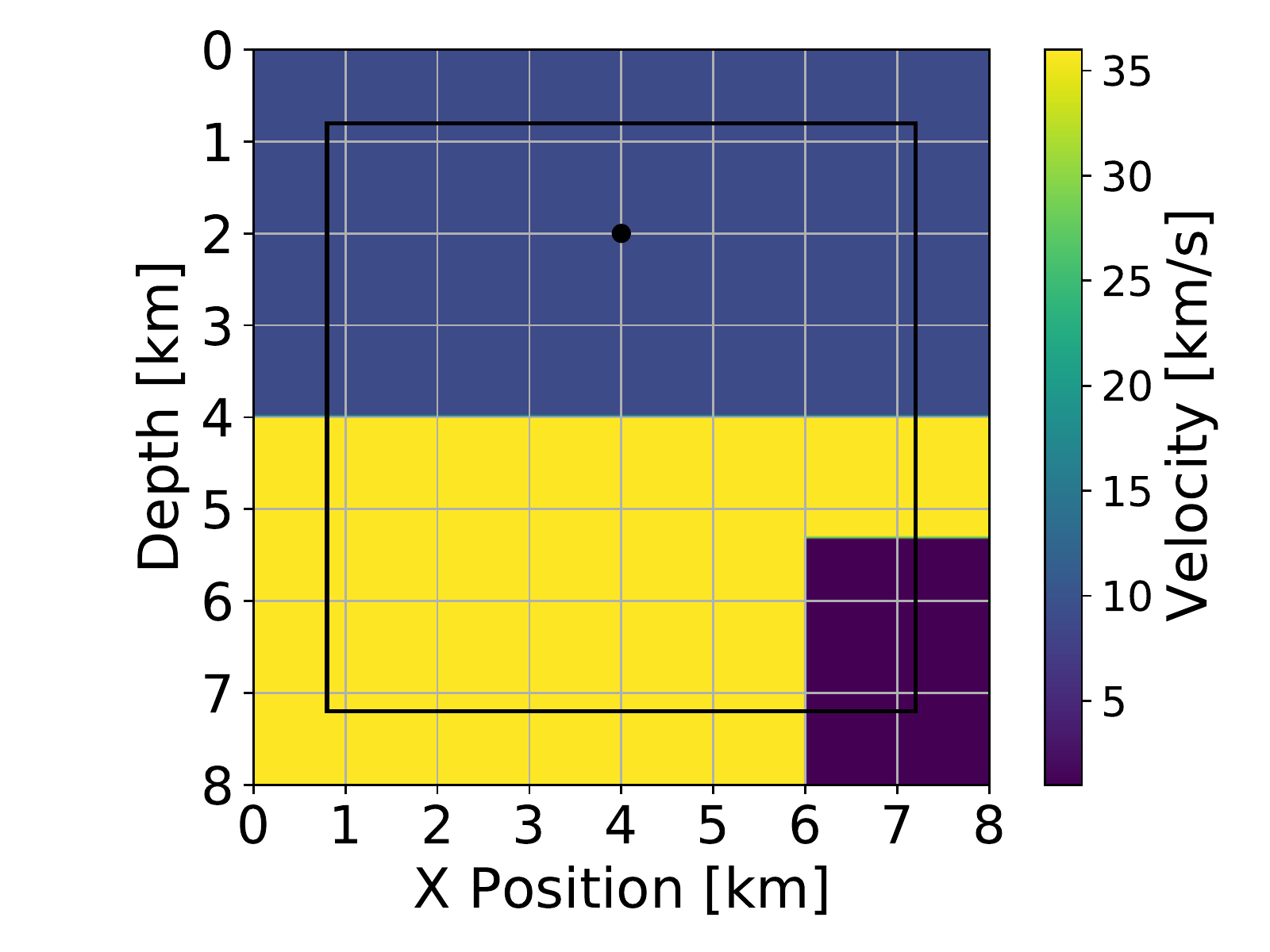}}
    \hfill
    \subfloat[Displacement snapshot$\quad$\linebreak at time $t=0.3$.]{\includegraphics[scale=0.24]{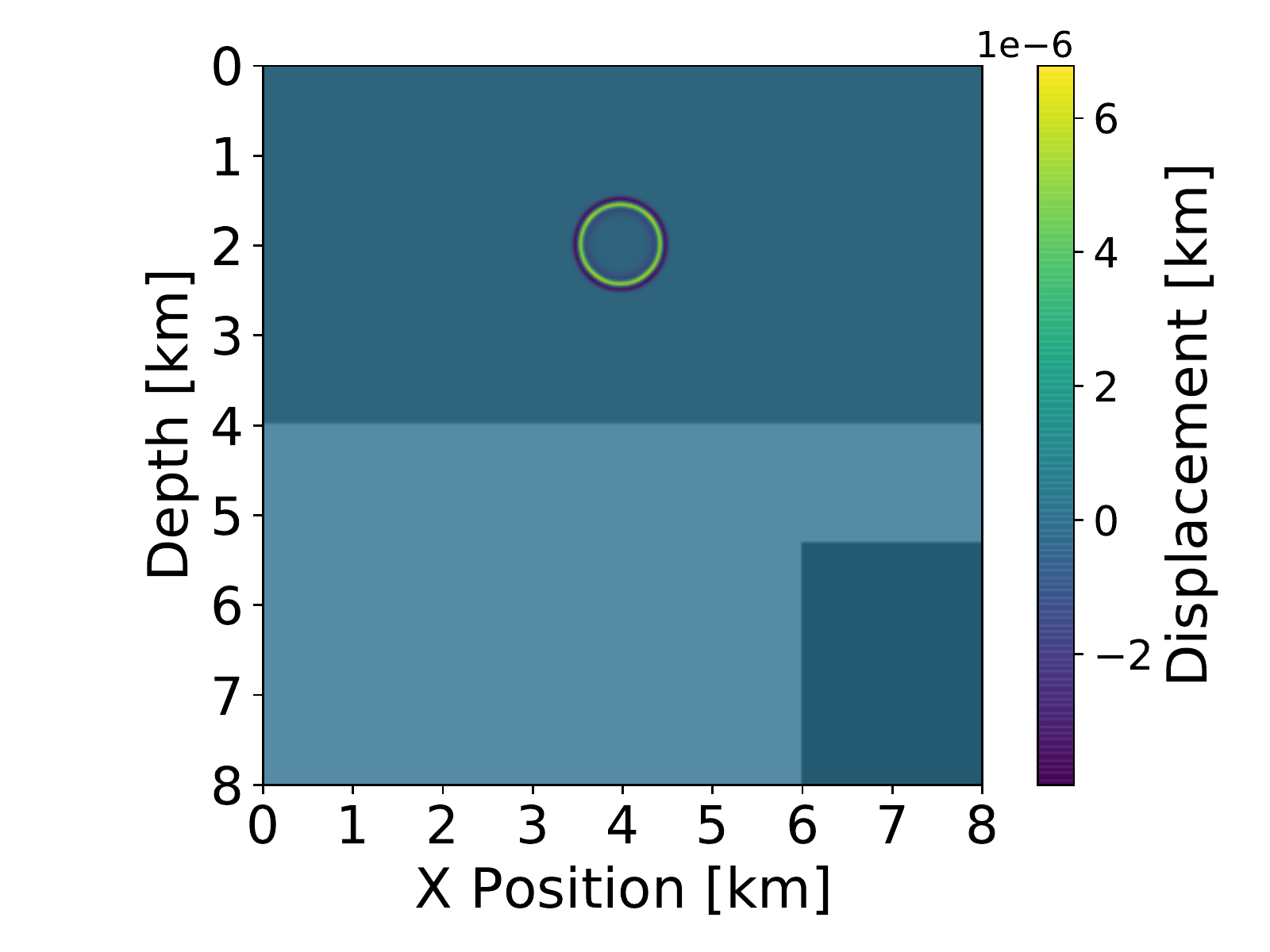}}
    \hfill
    \subfloat[Displacement snapshot$\quad$\linebreak at time $t=0.6$.]{\includegraphics[scale=0.24]{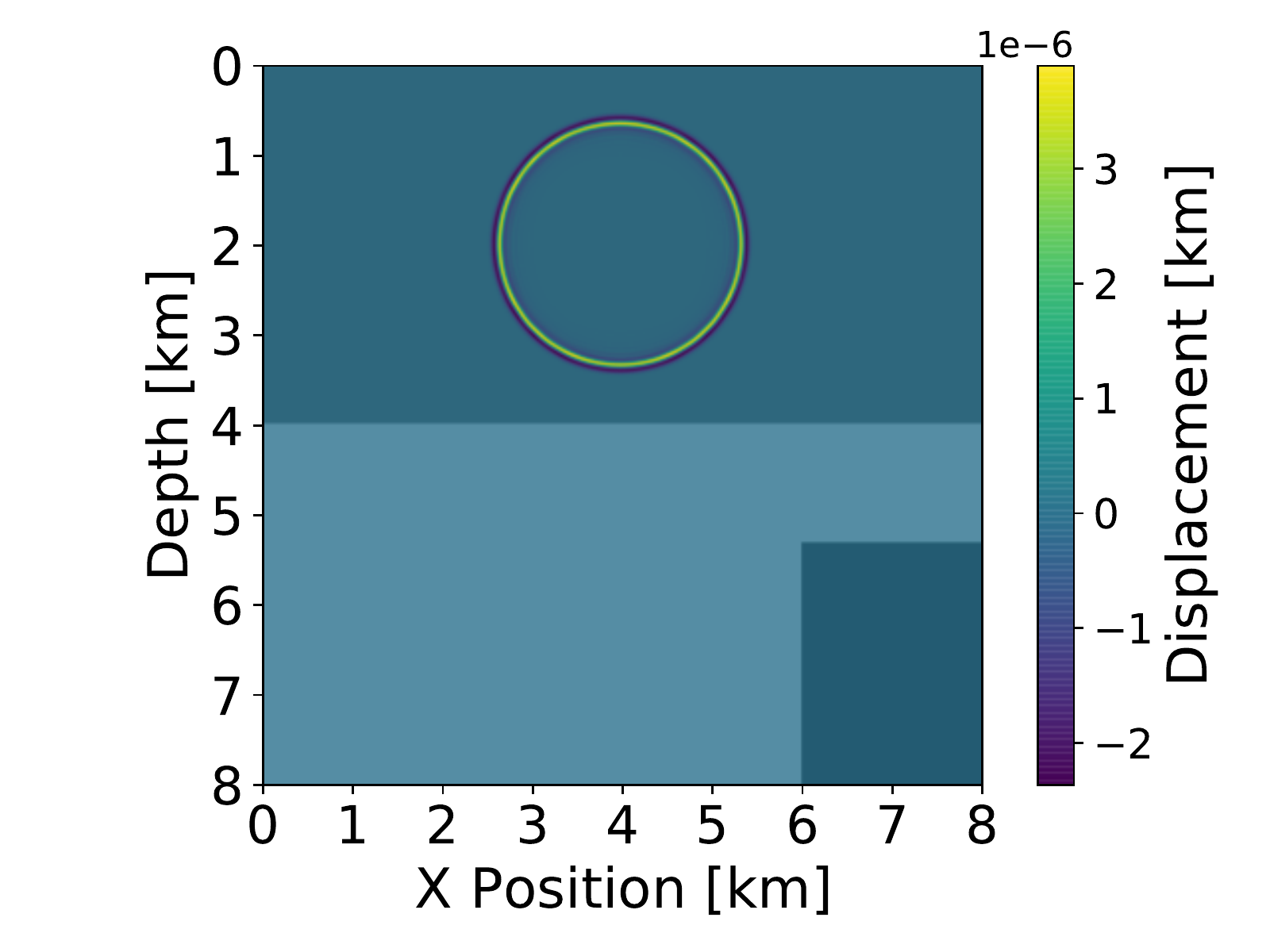}}
    \\
    \subfloat[Displacement snapshot$\quad$\linebreak at time $t=0.9$.]{\includegraphics[scale=0.24]{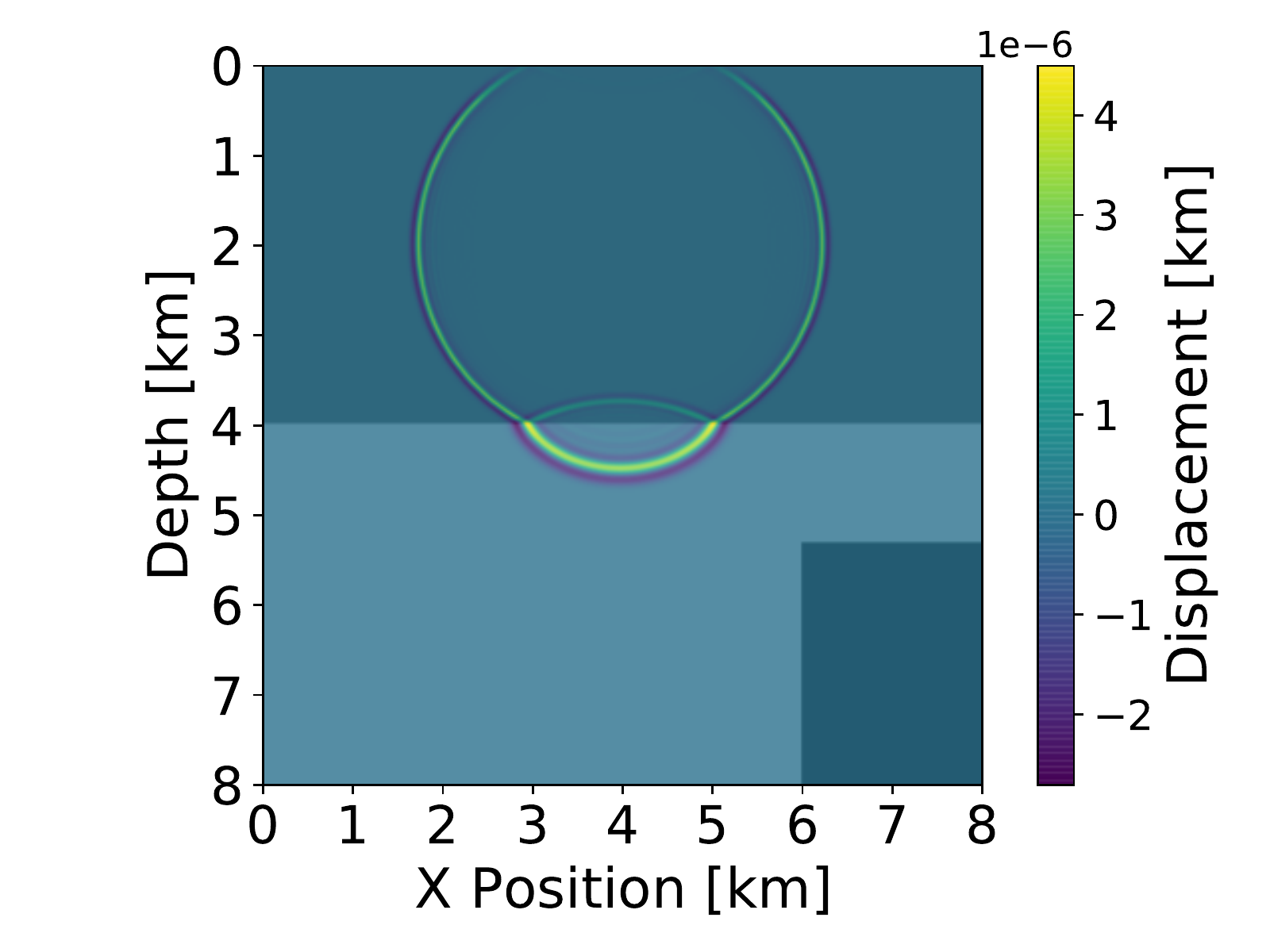}}
    \hfill
    \subfloat[Displacement snapshot$\quad$\linebreak at time $t=1.2$.]{\includegraphics[scale=0.24]{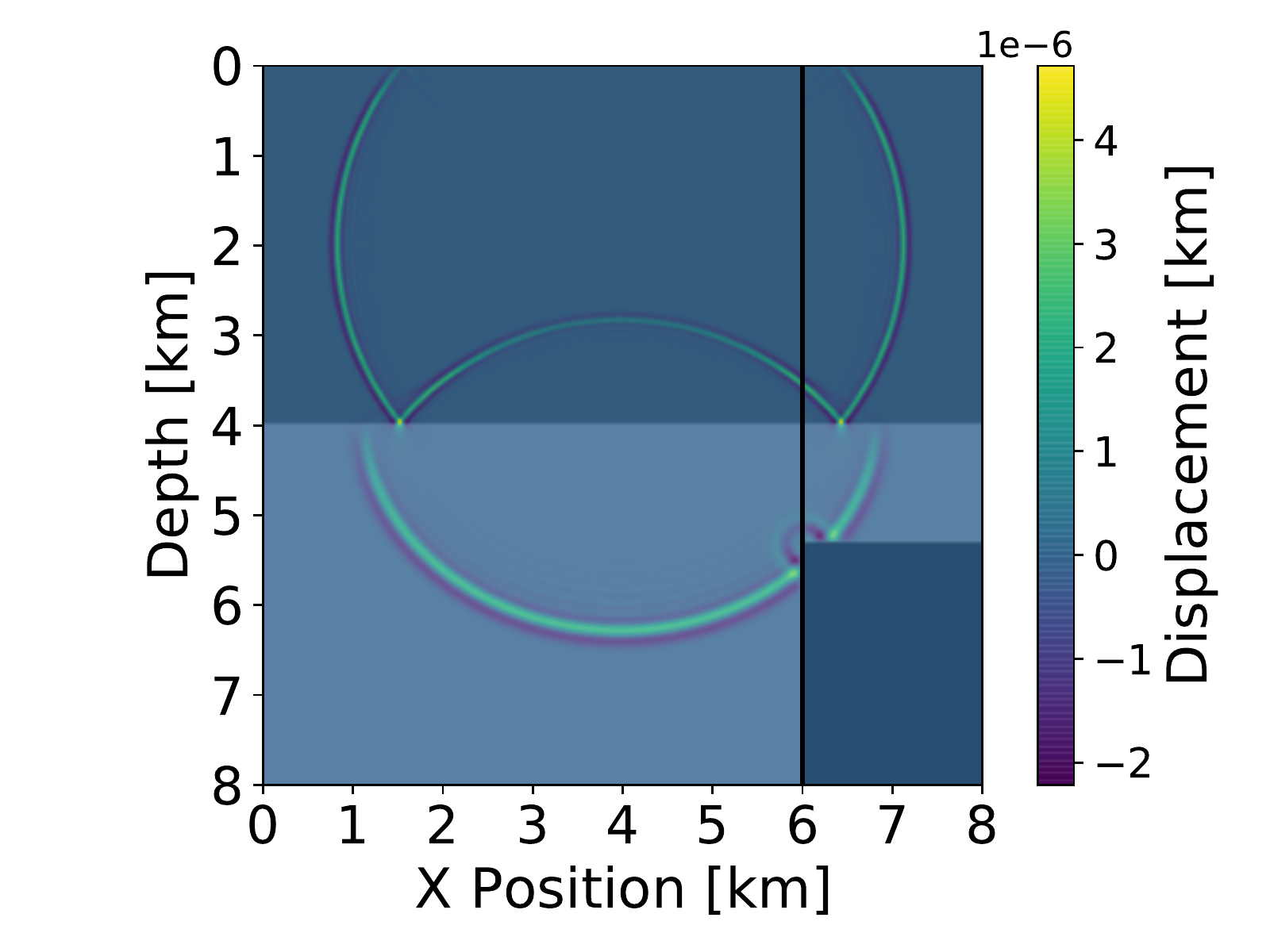}}
    \subfloat[Reference solution at time $t=1.2s.$]{\includegraphics[scale=0.25]{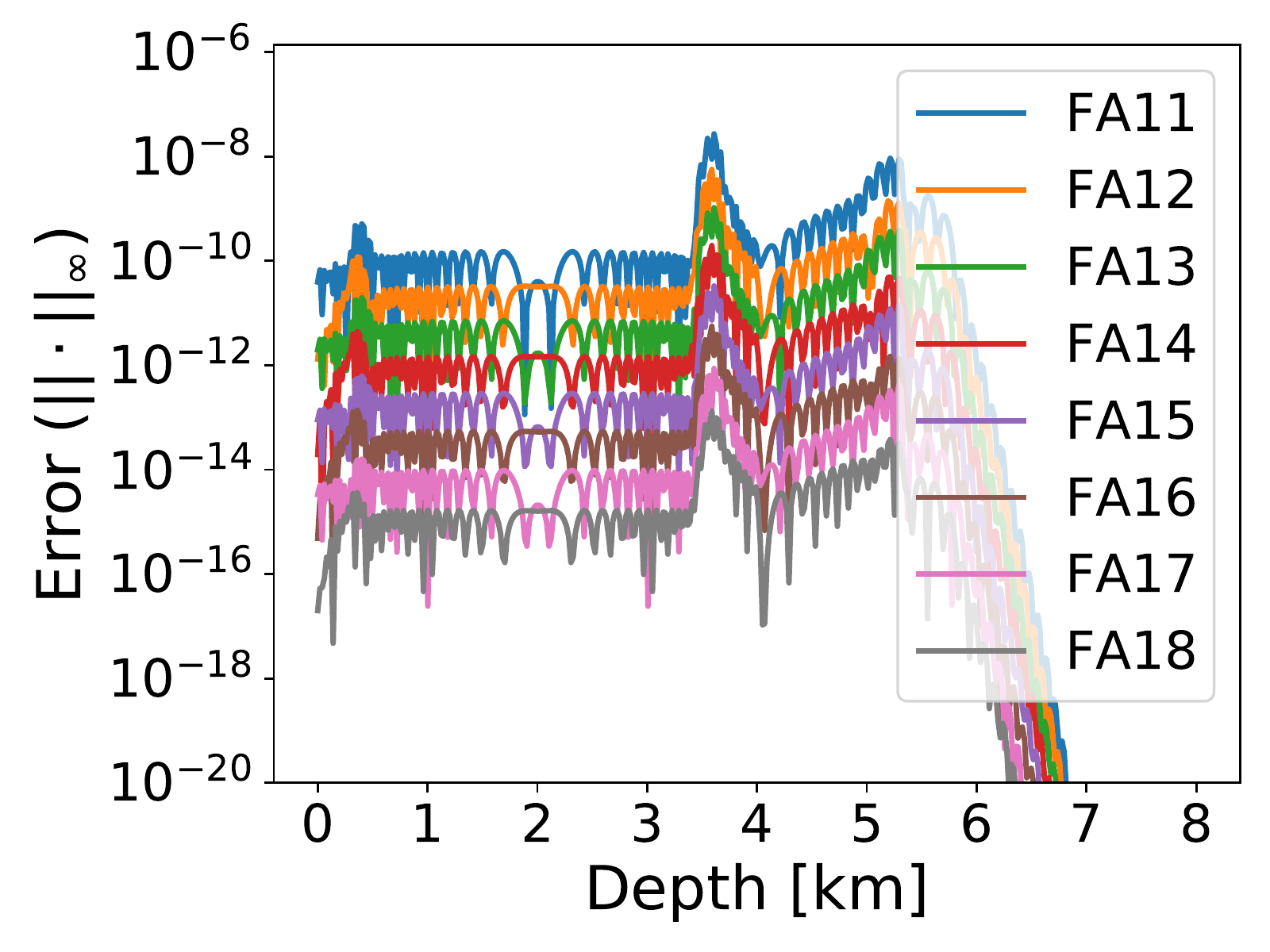}}
    \caption{
    (Subfigure (a))
    Acoustic wave propagation in a heterogeneous medium with a sharp corner.
    The physical domain is the region enclosed by the black squared contour, the region outside it is the PML domain, and the initial wave explosion arise at the dot black point inside the physical region.
    (Subfigures (b)-(e)) The reference solution is calculated at the final time and at intermediary time instants.
    (Subfigure (f)) Approximation errors at the final time $t=1.3$ for several degrees of Faber polynomials $m=\{10,\dots,18\}$ over the black vertical line in Subfigure (e).
    When the degree of the polynomial increases, the error calculated over the black line in Subfigure (e) diminishes by several orders of magnitude (Subfigure (f)).
    }\label{fig_snapshoots}
\end{figure}

Although the previous figures of convergence are useful to understand the approximation error with different setups, they offer little insight on the spatial distribution of the error. Now, we use the TC\#5, which has a high velocity contrast heterogeneous medium, to compute the approximation error along a straight line that cuts vertically the space (as shown in Fig. \ref{fig_snapshoots} (e)). The Faber solution is calculated with a timestep size $\Delta t$ which is $11 \times$ larger than the one used in the reference solution and its error is computed for several polynomial degrees.

Fig.\,\ref{fig_snapshoots} is composed by the velocity field of TC\#5, four wave propagation snap-shoots for different time instants, and Subfigure (f) showing the spatial error for different polynomial degrees over the line drawn in Subfigure (e).
If the reference solution was to be calculated with the same $\Delta t$ as the one used in the Faber approximations, it would have not converged to the solution and the same applies for the Faber polynomials with degree lower than $11$.  However, for the polynomials of degree higher than 10, we observe that the error decreases almost an order of magnitude for each increment of the degree.
This means that for the large fixed $\Delta t$, the addition of an approximation order produces a reduction of one order in the numerical error.
We again highlight that the dominant error show here is the temporal one, since our reference solution uses the same discrete operator as the exponential scheme.

We found that the numerical experiments are in agreement with the theoretical analysis of stability and dispersion. For higher polynomial degrees we have a larger CFL number, which translates into a larger time step sizes, such that high accuracy solution can be achieved with large time step sizes.

\section{Discussion, contributions \& outlook}\label{sec_conclussion}

Faber's polynomials provide a way to generalize exponential integration based on Chebychev polynomials to non-symmetric or non-antisymmetric matrices.
In the present work, we have developed sharper bounds of Faber approximations for normal matrices and a discussion about the importance of the conics used in the construction of the method. We showed that these conics are of utmost importance to ensure a fast convergence of the polynomial approximation.

We provided explicit bounds of the convex hull of the spectrum using different scenarios with respect to continuous model formulation and discretization schemes for the acoustic and elastic wave operators with the popular PML absorbing boundaries. These estimates remove the necessity of computing the eigenvalues of the full operator matrix, requiring only a calibration about the growth of the imaginary part with large $\Delta x$ for only one velocity model.
This way, we can specify a-priori the ranges of the spectrum on the discrete operator based on a simplified velocity model, which implies significant savings in the computational effort on the wave inversion step, where several wave propagation runs with changes in the velocity need to be solved.

Furthermore, we perform a study of stability, dispersion error, numerical convergence, and computational cost, for the Faber approximation of arbitrary order, using different wave equations formulations with PML parameters, spatial discretization order, and dimensions. We observed relatively small differences between distinct formulations of the continuum operator. In particular, we observed that 2SD performs better than the other formulation.
This suggests that the choice of an adequate formulation could improve even more the performance of Faber approximation.
We also found that increasing the order of the approximation implies a larger the CFL number, and the solution calculated for larger time steps maintains high accuracy. Due to the high order expansion of the source term (integral of Eq. \eqref{eq_const_var}), we also note that these results of large time step sizes are not affected by time-frequency sampling issues of the source term. Moreover, from the computational efficiency results, we conclude that, at least for the Faber polynomials applied to seismic waves, the increase in the polynomial degrees is also computationally more efficient than using lower order polynomials. We assume that this could lead to real improvements in seismic imaging, where dispersion errors are extremely important, and is very demanding on computational resources. 

In a future research, we intend to work on the open problem of Section \ref{sec_conics} and explore the conics that grants the fastest convergence for Faber polynomials.
Additionally, further investigation is required to prove the convergence of the spectra of the discrete operators to the continuum operators, as it has been done for the Schr\"{o}dinger equations in \citet{nakamura2021continuum}. This would lead to a theoretical prediction of the asymptotic behavior of the discrete operator's eigenvalues when $\Delta x$ tends to zero. 

We plan a follow-up work on this one on a comparison of the Faber polynomial exponential scheme with other exponential integrators and several classic methods in the context of the wave propagation equations and realistic seismic wave problems.

\begin{acknowledgements}

This research was carried out in association with the ongoing R\&D project registered as ANP20714-2 STMI - Software Technologies for Modelling and Inversion, with applications in seismic imaging (USP/Shell Brasil/ANP). It was funded in part by the Coordenação de Aperfeiçoamento de Pessoal de Nível Superior - Brasil (CAPES) - Finance Code 001, and in part by Conselho Nacional de Desenvolvimento Científico e Tecnológico (CNPq) - Brasil. Fundação de Amparo à Pesquisa do Estado de São Paulo (FAPESP) grant 2021/06176-0 is also acknowledged. It has also partially received funding from the Federal Ministry of Education and Research and the European High-Performance Computing Joint Undertaking (JU) under grant agreement No 955701, Time-X. The JU receives support from the European Union’s Horizon 2020 research and innovation programme and Belgium, France, Germany, Switzerland.
\end{acknowledgements}

\noindent{\small\textbf{Data Avaliability} All the data used in the paper was synthetically generated, following the instructions on the paper. The data and the methods codes are avaliable in the git-hub link
\hypertarget{https://github.com/fernanvr/Faber_1D_2D}{https://github.com/fernanvr/Faber\_1D\_2D}.}\\

\noindent{\small\textbf{Competing Interests} The authors have no relevant financial or non-financial interests to disclose.}

%
%

\bibliographystyle{spbasic}      
\bibliography{bibliography}   



\appendix
\section{Appendix}\label{sec_appendix}

\subsection{Continuous framework}\label{sec_appendix_continuous_framework}

We formally write the equations formulations, with PML conditions, used throughout the work.
For 1D the domain $\Omega=[a_1,a_2]$ is an interval and for 2D $\Omega=[0,a_2]\times[0,b_2]$ is considered a square.
The particular values of $a_1$ and $a_2$ are fixed in the numerical tests.

\noindent
\textbf{One dimensional acoustic waves with PML:}

\begin{enumerate}
    \item Using second order spatial derivatives (2SD)
    \begin{flalign}
            \frac{\partial u}{\partial t}(x,t) &=v(x,t), &&\\
             \frac{\partial v}{\partial t}(x,t)&=-\beta_x(x)v(x,t)+c^2(x)\left(\frac{\partial^2 u}{\partial x^2}(x,t)+\frac{\partial w}{\partial x}(x,t)\right)+f(x,t),&&\\\
             \frac{\partial w}{\partial t}(x,t)&=-\beta_x(x)\left(w(x,t)+\frac{\partial u}{\partial x}(x,t)\right).&&
\end{flalign}
    \item Using only first order spatial derivatives \citep{chern2019reflectionless} (1SD)
    \begin{flalign}
            \frac{\partial u}{\partial t}(x,t) &=c^2(x)\left(\frac{\partial v}{\partial x}(x,t) - w(x,t)\right) +\int\limits_{t_0}^t f(x,s)ds,\\
             \frac{\partial v}{\partial t}(x,t)&=-\beta_x(x)v(x,t)+\frac{\partial u}{\partial x}(x,t),\\
             \frac{\partial w}{\partial t}(x,t)&=\beta_x(x)\left(-w(x,t)+\frac{\partial v}{\partial x}(x,t)\right).&&
\end{flalign}
\end{enumerate}

\noindent
\textbf{Two dimension acoustic waves  with PML:}

\begin{enumerate}
    \item Using second order spatial derivatives (2SD)
    \begin{flalign}
            \frac{\partial u}{\partial t}(x,y,t) &=v(x,y,t),&&\\
             \frac{\partial v}{\partial t}(x,y,t)&=-\big(\beta_x(x)+\beta_y(y)\big)v(x,y,t)-\beta_x(x)\beta_y(y)u(x,y,t)\nonumber\\
             &+c^2(x,y)\left(\frac{\partial^2 u}{\partial x^2}(x,y,t)+\frac{\partial^2 u}{\partial y^2}(x,y,t)+\frac{\partial w_x}{\partial x}(x,y,t)\right.\nonumber\\
             &\left.+\frac{\partial w_y}{\partial y}(x,y,t)\right)+f(x,y,t),\\
             \frac{\partial w_x}{\partial t}(x,y,t)&=-\beta_x(x)w_x(x,y,t)+(\beta_y(y)-\beta_x(x))\frac{\partial u}{\partial x}(x,y,t),\\
             \frac{\partial w_y}{\partial t}(x,y,t)&=-\beta_y(y)w_y(x,y,t)+(\beta_x(x)-\beta_y(y))\frac{\partial u}{\partial y}(x,y,t).
\end{flalign}
    \item Using only first order spatial derivatives \citep{chern2019reflectionless} (1SD)
    \begin{flalign}
            \frac{\partial u}{\partial t}(x,y,t) &=c^2(x,y)\bigg(\frac{\partial v_x}{\partial x}(x,y,t) +\frac{\partial v_y}{\partial y}(x,y,t) -w_x(x,y,t)\nonumber\\
            &-w_y(x,y,t)\bigg)+\int\limits_{t_0}^tf(x,y,s)ds, \\
             \frac{\partial v_x}{\partial t}(x,y,t)&=-\beta_x(x)v_x(x,y,t)+\frac{\partial u}{\partial x}(x,y,t),\\
            \frac{\partial v_y}{\partial t}(x,y,t)&=-\beta_y(y)v_y(x,y,t)+\frac{\partial u}{\partial y}(x,y,t),\\
             \frac{\partial w_x}{\partial t}(x,y,t)&=\beta_x(x)\left(-w_x(x,y,t)+\frac{\partial v_x}{\partial x}(x,y,t)\right),\\
             \frac{\partial w_y}{\partial t}(x,y,t)&=\beta_y(y)\left(-w_y(x,y,t)+\frac{\partial v_y}{\partial y}(x,y,t)\right).&&
\end{flalign}
\end{enumerate}

\noindent
\textbf{Two-dimensional  elastic waves with PML, see \citet{assi2017compact}:}

{
\small
\begin{flalign}
        \frac{\partial u_x}{\partial t}(x,y,t) &=v_x(x,y,t),\\
        \frac{\partial u_y}{\partial t}(x,y,t) &=v_y(x,y,t),\\
             \frac{\partial v_x}{\partial t}(x,y,t)&=-\big(\beta_x(x)+\beta_y(y)\big)v_x(x,y,t)-\beta_x(x)\beta_y(y)u_x(x,y,t)\nonumber\\
             &+\frac{1}{\rho(x,y)}\left[\frac{\partial}{\partial x}\left(T_{xx}(x,y,t)+w_{xx}(x,y,t)\right)\right.\nonumber\\
             &\left.+\frac{\partial}{\partial y}\left(T_{xy}(x,y,t)+w_{xy}(x,y,t)\right)\right]+f_x(x,y,t),\\
             \frac{\partial v_y}{\partial t}(x,y,t)&=-\big(\beta_x(x)+\beta_y(y)\big)v_y(x,y,t)-\beta_x(x)\beta_y(y)u_y(x,y,t)\nonumber\\
             &+\frac{1}{\rho}\left[\frac{\partial}{\partial x}\left(T_{xy}(x,y,t)+w_{yx}(x,y,t)\right)\right.\nonumber\\
             &\left.+\frac{\partial}{\partial y}\left(T_{yy}(x,y,t)+w_{yy}(x,y,t)\right)\right]+f_y(x,y,t),\\
             \frac{\partial T_{xx}}{\partial t}(x,y,t)&=\big(2\mu(x,y)+\lambda(x,y)\big)\frac{\partial v_x}{\partial x}(x,y,t)+\lambda(x,y)\frac{\partial v_y}{\partial y}(x,y,t),\\
             \frac{\partial T_{xy}}{\partial t}(x,y,t)&=\mu(x,y)\left(\frac{\partial v_x}{\partial y}(x,y,t)+\frac{\partial v_y}{\partial x}(x,y,t)\right),\\
             \frac{\partial T_{yy}}{\partial t}(x,y,t)&=\lambda(x,y)\frac{\partial v_x}{\partial x}(x,y,t)+\big(2\mu(x,y)+\lambda(x,y)\big)\frac{\partial v_y}{\partial y}(x,y,t),\\
             \frac{\partial w_{xx}}{\partial t}(x,y,t)&=-\beta_x(x)w_{xx}(x,y,t)+\big(\beta_y(y)-\beta_x(x)\big)\big(2\mu(x,y)+\lambda(x,y)\big)\frac{\partial u_x}{\partial x}(x,y,t),\\
             \frac{\partial w_{xy}}{\partial t}(x,y,t)&=-\beta_y(y)w_{xy}(x,y,t)+\big(\beta_x(x)-\beta_y(y)\big)\mu(x,y)\frac{\partial u_x}{\partial y}(x,y,t),\\
             \frac{\partial w_{yx}}{\partial t}(x,y,t)&=-\beta_x(x)w_{yx}(x,y,t)+\big(\beta_y(y)-\beta_x(x)\big)\mu(x,y)\frac{\partial u_y}{\partial x}(x,y,t),\\
             \frac{\partial w_{yy}}{\partial t}(x,y,t)&=-\beta_y(y)w_{yy}(x,y,t)+\big(\beta_x(x)-\beta_y(y)\big)\big(2\mu(x,y)+\lambda(x,y)\big)\frac{\partial u_y}{\partial y}(x,y,t),
\end{flalign}
}
with the variables and parameters described in Table \ref{tab_const}.

\renewcommand{\arraystretch}{1.2}
\begin{center}
\begin{table}[H]
\begin{tabular}{|p{3cm}|p{8cm}|} 
	\hline
	\textbf{Symbol} & \textbf{Description}\\
	\hline
	$x,\;y$ & Spatial variables\\
	\hline
	$t$ & Time \\
	\hline
	$t_0$ & Initial time instant\\
	\hline
	$u(x,y,t)$ & Displacement of acoustic waves\\
	\hline
	$v(x,y,t)$ & Displacement velocity for 2SD, and material velocity for 1SD\\\hline
	$u_x(x,y,t)$, $u_y(x,y,t)$ & Displacement of elastic waves in $x$ and $y$ directions, respectively\\
	\hline
	$v_x(x,y,t)$, $v_x(x,y,t)$ & Displacement velocity for elastic waves, and material velocity for 1SOD, in $x$ and $y$ directions, respectively\\
	\hline
	$T_{xx}(x,y,t)$, $T_{xy}(x,y,t)$, $T_{yy}(x,y,t)$ & Stress components of elastic waves\\\hline
	$w(x,y,t)$, $w_x(x,y,t)$, $w_y(x,y,t)$, $w_{xx}(x,y,t)$, $w_{xy}(x,y,t)$, $w_{yx}(x,y,t)$, $w_{yy}(x,y,t)$
	& auxiliary variables of the PML boundary condition\\
	\hline
	$c(x,y)$ & wave propagation velocities  in 1D and 2D, respectively\\
	\hline
	$\mu(x,y),\;\lambda(x,y)$ & Lamé parameters\\
	\hline
	$\rho(x,y)$ & density\\
	\hline
	$\beta_x(x,y),\beta_y(x,y)$ & wave damping functions\\
	\hline
\end{tabular}
\caption{Variables used in the equations and their description}\label{tab_const}
\end{table}
\end{center}

The damping functions $\beta_z$, related to the absorption factor are defined as

\begin{equation}
\beta_z(z)=\left\{\begin{array}{ll}
0,&\text{ if } d(z,\partial\Omega)>\delta\\
\beta_0\left(\frac{d(z,\Omega_1)}{\delta}\right)^2,& \text{ if }d(z,\partial\Omega)\leq\delta
\end{array}\right.,\quad z=x,\;y
\end{equation}
where $d(z,\partial \Omega)$ is the distance from $z$ to the boundary of $\Omega$, $\delta$ is the thickness of the PML domain, $\beta_0$ is the magnitude of the absorption factor, and $\Omega_1$ is the numerical domain without the damping layer (physical domain). Thus, $\Omega$ is composed by the union of $\Omega_1$ and a damping layer of thickness $\delta$ extending on the boundary of $\Omega_1$.\\

\subsection{Discrete framework}\label{sec_appendix_numerical_framework}

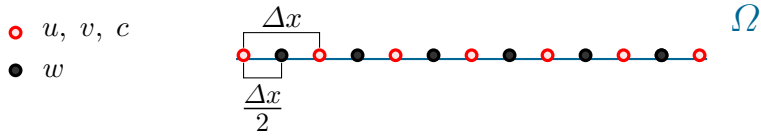
\begin{figure}[H]
\centering
\begin{tikzpicture}
\draw[MidnightBlue, thick] (-0.1,-0.05) -- (6.1,-0.05);
\node at (6.6,0.5) {\Large\textcolor{MidnightBlue}{$\Omega$}};
\filldraw[color=red, fill=red!10, very thick](0,0) circle (0.07);
\filldraw[color=red, fill=red!10, very thick](1,0) circle (0.07);
\filldraw[color=red, fill=red!10, very thick](2,0) circle (0.07);
\filldraw[color=red, fill=red!10, very thick](3,0) circle (0.07);
\filldraw[color=red, fill=red!10, very thick](4,0) circle (0.07);
\filldraw[color=red, fill=red!10, very thick](5,0) circle (0.07);
\filldraw[color=red, fill=red!10, very thick](6,0) circle (0.07);
\filldraw[color=black, fill=black!80, very thick](0.5,0) circle (0.07);
\filldraw[color=black, fill=black!80, very thick](1.5,0) circle (0.07);
\filldraw[color=black, fill=black!80, very thick](2.5,0) circle (0.07);
\filldraw[color=black, fill=black!80, very thick](3.5,0) circle (0.07);
\filldraw[color=black, fill=black!80, very thick](4.5,0) circle (0.07);
\filldraw[color=black, fill=black!80, very thick](5.5,0) circle (0.07);
\filldraw[color=black, fill=black!80, very thick](5.5,0) circle (0.07);
\filldraw[color=red, fill=red!10, very thick](-3,0.3) circle (0.07);
\filldraw[color=black, fill=black!80, very thick](-3,-0.2) circle (0.07);

\node at (-2.08,0.3) {\large $u,\;v,\;c$};
\node at (-2.5,-0.2) {\large $w$};
\draw[black] (0,0.1) -- (0,0.3) -- (1,0.3) -- (1,0.1);
\node at (0.5,0.5) {\large$\Delta x$};
\draw[black] (0,-0.1) -- (0,-0.3) -- (0.5,-0.3) -- (0.5,-0.1);
\node at (0.25,-0.7) {\Large$\frac{\Delta x}{2}$};

\end{tikzpicture}
\caption{Staggered grid in 1D with the relative positions of the (2SD) and (1SD) wave equation variables and parameters. $u,\;v$ and $c$ are collocated (centered) and $w$ is staggered in the grid.}\label{fig_staggered_1D_b}
\end{figure}

\begin{figure}[H]
\centering
\begin{tikzpicture}
\draw[blue, very thick] (0,0) -- (0,3) -- (6,3) -- (6,0) -- (0,0);
\node at (6.8,2.4) {\Huge\textcolor{blue}{$\Omega$}};
\filldraw[color=red, fill=red!10, very thick](0,0) circle (0.07);
\filldraw[color=red, fill=red!10, very thick](1,0) circle (0.07);
\filldraw[color=red, fill=red!10, very thick](2,0) circle (0.07);
\filldraw[color=red, fill=red!10, very thick](3,0) circle (0.07);
\filldraw[color=red, fill=red!10, very thick](4,0) circle (0.07);
\filldraw[color=red, fill=red!10, very thick](5,0) circle (0.07);
\filldraw[color=red, fill=red!10, very thick](6,0) circle (0.07);
\filldraw[color=black, fill=black!80, very thick](0.5,0) circle (0.07);
\filldraw[color=black, fill=black!80, very thick](1.5,0) circle (0.07);
\filldraw[color=black, fill=black!80, very thick](2.5,0) circle (0.07);
\filldraw[color=black, fill=black!80, very thick](3.5,0) circle (0.07);
\filldraw[color=black, fill=black!80, very thick](4.5,0) circle (0.07);
\filldraw[color=black, fill=black!80, very thick](5.5,0) circle (0.07);
\filldraw[color=black, fill=black!80, very thick](5.5,0) circle (0.07);

\filldraw[color=red, fill=red!10, very thick](0,1) circle (0.07);
\filldraw[color=red, fill=red!10, very thick](1,1) circle (0.07);
\filldraw[color=red, fill=red!10, very thick](2,1) circle (0.07);
\filldraw[color=red, fill=red!10, very thick](3,1) circle (0.07);
\filldraw[color=red, fill=red!10, very thick](4,1) circle (0.07);
\filldraw[color=red, fill=red!10, very thick](5,1) circle (0.07);
\filldraw[color=red, fill=red!10, very thick](6,1) circle (0.07);
\filldraw[color=black, fill=black!80, very thick](0.5,1) circle (0.07);
\filldraw[color=black, fill=black!80, very thick](1.5,1) circle (0.07);
\filldraw[color=black, fill=black!80, very thick](2.5,1) circle (0.07);
\filldraw[color=black, fill=black!80, very thick](3.5,1) circle (0.07);
\filldraw[color=black, fill=black!80, very thick](4.5,1) circle (0.07);
\filldraw[color=black, fill=black!80, very thick](5.5,1) circle (0.07);
\filldraw[color=black, fill=black!80, very thick](5.5,1) circle (0.07);

\filldraw[color=red, fill=red!10, very thick](0,2) circle (0.07);
\filldraw[color=red, fill=red!10, very thick](1,2) circle (0.07);
\filldraw[color=red, fill=red!10, very thick](2,2) circle (0.07);
\filldraw[color=red, fill=red!10, very thick](3,2) circle (0.07);
\filldraw[color=red, fill=red!10, very thick](4,2) circle (0.07);
\filldraw[color=red, fill=red!10, very thick](5,2) circle (0.07);
\filldraw[color=red, fill=red!10, very thick](6,2) circle (0.07);
\filldraw[color=black, fill=black!80, very thick](0.5,2) circle (0.07);
\filldraw[color=black, fill=black!80, very thick](1.5,2) circle (0.07);
\filldraw[color=black, fill=black!80, very thick](2.5,2) circle (0.07);
\filldraw[color=black, fill=black!80, very thick](3.5,2) circle (0.07);
\filldraw[color=black, fill=black!80, very thick](4.5,2) circle (0.07);
\filldraw[color=black, fill=black!80, very thick](5.5,2) circle (0.07);
\filldraw[color=black, fill=black!80, very thick](5.5,2) circle (0.07);

\filldraw[color=red, fill=red!10, very thick](0,3) circle (0.07);
\filldraw[color=red, fill=red!10, very thick](1,3) circle (0.07);
\filldraw[color=red, fill=red!10, very thick](2,3) circle (0.07);
\filldraw[color=red, fill=red!10, very thick](3,3) circle (0.07);
\filldraw[color=red, fill=red!10, very thick](4,3) circle (0.07);
\filldraw[color=red, fill=red!10, very thick](5,3) circle (0.07);
\filldraw[color=red, fill=red!10, very thick](6,3) circle (0.07);
\filldraw[color=black, fill=black!80, very thick](0.5,3) circle (0.07);
\filldraw[color=black, fill=black!80, very thick](1.5,3) circle (0.07);
\filldraw[color=black, fill=black!80, very thick](2.5,3) circle (0.07);
\filldraw[color=black, fill=black!80, very thick](3.5,3) circle (0.07);
\filldraw[color=black, fill=black!80, very thick](4.5,3) circle (0.07);
\filldraw[color=black, fill=black!80, very thick](5.5,3) circle (0.07);
\filldraw[color=black, fill=black!80, very thick](5.5,3) circle (0.07);

\filldraw[color=Plum, fill=Plum!10, very thick](0,0.5) circle (0.07);
\filldraw[color=Plum, fill=Plum!10, very thick](1,0.5) circle (0.07);
\filldraw[color=Plum, fill=Plum!10, very thick](2,0.5) circle (0.07);
\filldraw[color=Plum, fill=Plum!10, very thick](3,0.5) circle (0.07);
\filldraw[color=Plum, fill=Plum!10, very thick](4,0.5) circle (0.07);
\filldraw[color=Plum, fill=Plum!10, very thick](5,0.5) circle (0.07);
\filldraw[color=Plum, fill=Plum!10, very thick](6,0.5) circle (0.07);
\filldraw[color=Plum, fill=Plum!10, very thick](0,1.5) circle (0.07);
\filldraw[color=Plum, fill=Plum!10, very thick](1,1.5) circle (0.07);
\filldraw[color=Plum, fill=Plum!10, very thick](2,1.5) circle (0.07);
\filldraw[color=Plum, fill=Plum!10, very thick](3,1.5) circle (0.07);
\filldraw[color=Plum, fill=Plum!10, very thick](4,1.5) circle (0.07);
\filldraw[color=Plum, fill=Plum!10, very thick](5,1.5) circle (0.07);
\filldraw[color=Plum, fill=Plum!10, very thick](6,1.5) circle (0.07);
\filldraw[color=Plum, fill=Plum!10, very thick](0,2.5) circle (0.07);
\filldraw[color=Plum, fill=Plum!10, very thick](1,2.5) circle (0.07);
\filldraw[color=Plum, fill=Plum!10, very thick](2,2.5) circle (0.07);
\filldraw[color=Plum, fill=Plum!10, very thick](3,2.5) circle (0.07);
\filldraw[color=Plum, fill=Plum!10, very thick](4,2.5) circle (0.07);
\filldraw[color=Plum, fill=Plum!10, very thick](5,2.5) circle (0.07);
\filldraw[color=Plum, fill=Plum!10, very thick](6,2.5) circle (0.07);
\node at (0,0.5) {+};
\node at (1,0.5) {+};
\node at (2,0.5) {+};
\node at (3,0.5) {+};
\node at (4,0.5) {+};
\node at (5,0.5) {+};
\node at (6,0.5) {+};
\node at (0,1.5) {+};
\node at (1,1.5) {+};
\node at (2,1.5) {+};
\node at (3,1.5) {+};
\node at (4,1.5) {+};
\node at (5,1.5) {+};
\node at (6,1.5) {+};
\node at (0,2.5) {+};
\node at (1,2.5) {+};
\node at (2,2.5) {+};
\node at (3,2.5) {+};
\node at (4,2.5) {+};
\node at (5,2.5) {+};
\node at (6,2.5) {+};

\filldraw[color=Dandelion, fill=Dandelion!10, very thick](0.5,0.5) circle (0.07);
\filldraw[color=Dandelion, fill=Dandelion!10, very thick](1.5,0.5) circle (0.07);
\filldraw[color=Dandelion, fill=Dandelion!10, very thick](2.5,0.5) circle (0.07);
\filldraw[color=Dandelion, fill=Dandelion!10, very thick](3.5,0.5) circle (0.07);
\filldraw[color=Dandelion, fill=Dandelion!10, very thick](4.5,0.5) circle (0.07);
\filldraw[color=Dandelion, fill=Dandelion!10, very thick](5.5,0.5) circle (0.07);
\filldraw[color=Dandelion, fill=Dandelion!10, very thick](0.5,1.5) circle (0.07);
\filldraw[color=Dandelion, fill=Dandelion!10, very thick](1.5,1.5) circle (0.07);
\filldraw[color=Dandelion, fill=Dandelion!10, very thick](2.5,1.5) circle (0.07);
\filldraw[color=Dandelion, fill=Dandelion!10, very thick](3.5,1.5) circle (0.07);
\filldraw[color=Dandelion, fill=Dandelion!10, very thick](4.5,1.5) circle (0.07);
\filldraw[color=Dandelion, fill=Dandelion!10, very thick](5.5,1.5) circle (0.07);
\filldraw[color=Dandelion, fill=Dandelion!10, very thick](0.5,2.5) circle (0.07);
\filldraw[color=Dandelion, fill=Dandelion!10, very thick](1.5,2.5) circle (0.07);
\filldraw[color=Dandelion, fill=Dandelion!10, very thick](2.5,2.5) circle (0.07);
\filldraw[color=Dandelion, fill=Dandelion!10, very thick](3.5,2.5) circle (0.07);
\filldraw[color=Dandelion, fill=Dandelion!10, very thick](4.5,2.5) circle (0.07);
\filldraw[color=Dandelion, fill=Dandelion!10, very thick](5.5,2.5) circle (0.07);
\node at (0.5,0.5) {x};
\node at (1.5,0.5) {x};
\node at (2.5,0.5) {x};
\node at (3.5,0.5) {x};
\node at (4.5,0.5) {x};
\node at (5.5,0.5) {x};
\node at (0.5,1.5) {x};
\node at (1.5,1.5) {x};
\node at (2.5,1.5) {x};
\node at (3.5,1.5) {x};
\node at (4.5,1.5) {x};
\node at (5.5,1.5) {x};
\node at (0.5,2.5) {x};
\node at (1.5,2.5) {x};
\node at (2.5,2.5) {x};
\node at (3.5,2.5) {x};
\node at (4.5,2.5) {x};
\node at (5.5,2.5) {x};

\filldraw[color=red, fill=red!10, very thick](0,5.5) circle (0.07);
\filldraw[color=black, fill=black!80, very thick](0,5) circle (0.07);
\filldraw[color=Plum, fill=Plum!10, very thick](0,4.5) circle (0.07);
\filldraw[color=Dandelion, fill=Dandelion!10, very thick](0,4) circle (0.07);
\node at (0,4.5) {+};
\node at (0,4) {x};

\node at (1.8,5.5) {\large $u,\;v,\;u_x,\;v_x,\;c,\;\rho$};
\node at (2.8,5) {\large $T_{xx},\;T_{yy},\;w_x,\;w_{xx},\;w_{yy},\;\mu,\;\lambda$};
\node at (2.17,4.5) {\large $T_{xy},\;w_y,\;w_{xy},\;w_{yx},\;\mu$};
\node at (1.13,4) {\large $u_y,\;v_y,\;\rho$};
\draw[black] (0,-0.1) -- (0,-0.3) -- (0.5,-0.3) -- (0.5,-0.1);
\node at (0.25,-0.7) {\Large$\frac{\Delta x}{2}$};
\draw[black] (-0.1,0.5) -- (-0.3,0.5) -- (-0.3,1) -- (-0.1,1);
\node at (-0.7,0.75) {\Large$\frac{\Delta x}{2}$};
\end{tikzpicture}
\caption{Staggered grid in 2D with the relative positions of the (2SD and elastic) wave equations variables and parameters. $u,\;v,\;u_x,\;v_x,\;\rho$ and $c$ are collocated.}\label{fig_staggered_2D_b}
\end{figure}
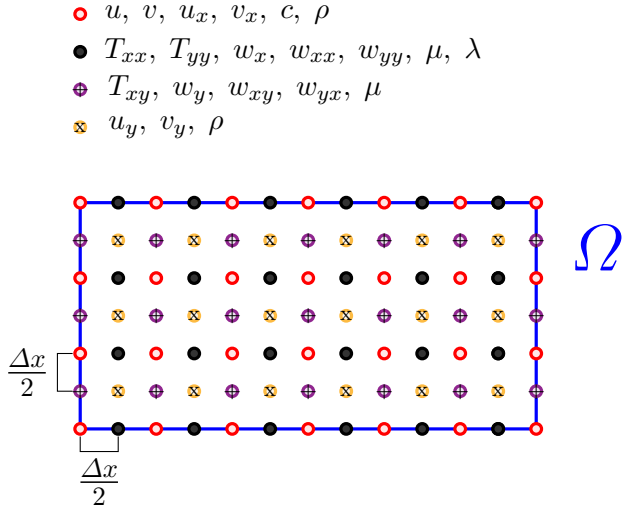

The spatial discretizations are based on a staggered grid using 4th and 8th order approximation of the spatial derivatives defined by equations \eqref{eq_spatial_4th_1} and \eqref{eq_spatial_8th_1}. Figs.\,\ref{fig_staggered_1D_b} and \ref{fig_staggered_2D_b} describe the discrete space for the 1SD and 2SD formulation in 1D, and the 2SD and elastic formulations in 2D, respectively. 

\subsection{Numerical benchmarks}\label{sec_appendix_experiments_formulation}

We define the numerical experiments, called ``Test Case'' used through the paper.
In all the tests we use a zero Dirichlet condition on the boundary of the domain $\Omega$, a PML layer thickness of $\delta=0.8\text{km}$, a damping parameter $\beta_0=30$, and Ricker peak frequency of $f_0=25\text{Hz}$.
If not otherwise stated, the initial condition for all the variables is zero.
The particular benchmarks are then defined as follows:

~\\

\noindent
\textbf{Test Case \#1:} $\Omega=[0,10.5\text{km}]$
    \begin{equation*}
        c\equiv 1.524\text{km}/\text{s},\quad u_0(x)=((1-10(x-5.25)^2)e^{-10(x-5.25)^2},\quad f\equiv 0
    \end{equation*}

\noindent
\textbf{Test Case \#2:} $\Omega=[0,10.5\text{km}]$
    \begin{align*}
        c(x)&=\left\{\begin{array}{ll}
             1.524\;\text{km}/\text{s}, & \text{if } x<5.25 \\
             3.048\;\text{km}/\text{s} &   \text{if } 5.25\leq x < 7\\
             0.1524\;\text{km}/\text{s} &   \text{if } 7\leq x
        \end{array}\right.,\\
        u_0&=\left\{\begin{array}{ll}
            0, & \text{if } |x-2.6|\geq 0.01\\
            e^{\frac{(x-2.6)^2}{(x-2.6)^2-0.01^2}}, & \text{if } |x-2.6|<0.01
        \end{array}\right.,\quad f\equiv 0
    \end{align*}
    
\noindent
\textbf{Test Case \#3:} $\Omega=[0,10.5\text{km}]$
    \begin{align*}
        c(x)&=\left\{\begin{array}{ll}
             1.524\;\text{km}/\text{s}, & \text{if } x<5.25 \\
             3.048\;\text{km}/\text{s} &   \text{if } 5.25\leq\\
        \end{array}\right.,\quad u_0\equiv0,\\
        f(x,t)&=\left\{\begin{array}{ll}
            0, & \text{if } |x-2.6|\geq 0.01\\
            e^{\frac{(x-2.6)^2}{(x-2.6)^2-0.01^2}}(1-f_0^2\pi^2(t-t_0)^2)e^{-f_0^2\pi^2(t-t_0)^2}, & \text{if } |x-2.6|<0.01
        \end{array}\right.
    \end{align*}

\noindent
\textbf{Test Case \#4:} $\Omega=[0,8\;\text{km}]\times[0,8\;\text{km}]$
    \begin{align*}
        c&\equiv 3\;\text{km}/\text{s},\quad u_0(x,y)=\left\{\begin{array}{ll}
             0,& \text{if } \|(x,y)-(4,2)\|\geq 0.01  \\
             e^{\frac{\|(x,y)-(4,2)\|^2}{\|(x,y)-(4,2)\|^2-0.01^2}},& \text{if } \|(x,y)-(4,2)\|< 0.01 
        \end{array}\right.,\\
        f&\equiv 0
    \end{align*}

\noindent
\textbf{Test Case \#5:} $\Omega=[0,8\;\text{km}]\times[0,8\;\text{km}]$
    \begin{align*}
        c(x,y)&=\left\{\begin{array}{ll}
             3\;\text{km}/\text{s}, & \text{if } y\geq4 \\
             6\;\text{km}/\text{s}, &   \text{if } (y<4 \text{ and } x\leq6) \text{ or } (16/3< y<4 \text{ and } x>6)\\
             1\;\text{km}/\text{s}, &   \text{if } 6\leq x \text{ and } y\leq16/3
        \end{array}\right.,\\
        u_0&=\left\{\begin{array}{ll}
             0,& \text{if } \|(x,y)-(4,2)\|\geq 0.01  \\
             e^{\frac{\|(x,y)-(4,2)\|^2}{\|(x,y)-(4,2)\|^2-0.01^2}},& \text{if } \|(x,y)-(4,2)\|< 0.01 
        \end{array}\right.,\quad f\equiv 0
    \end{align*}
    
\noindent
\textbf{Test Case \#6:} $\Omega=[0,8\;\text{km}]\times[0,8\;\text{km}]$
    \begin{align*}
        c(x,y)&=\left\{\begin{array}{ll}
             3\;\text{km}/\text{s}, & \text{if } y\geq4 \\
             6\;\text{km}/\text{s}, &   \text{if } y<4\\
        \end{array}\right.,\quad u_0\equiv0,\\
        f(x,y,t)&=\left\{\begin{array}{ll}
            0, & \text{if } \|(x,y)-(4,2)\|\geq 0.01\\
            e^{\frac{\|(x,y)-(4,2)\|^2}{\|(x,y)-(4,2)\|^2-0.01^2}}(1-f_0^2\pi(t-t_0)^2)e^{-f_0^2\pi^2(t-t_0)^2}, & \text{if } \|(x,y)-(4,2)\|<0.01
        \end{array}\right.
    \end{align*}
    
\noindent
\textbf{Test Case \#7:} $\Omega=[0,8\;\text{km}]\times[0,8\;\text{km}]$ for elastic waves
    \begin{align*}
        \rho&\equiv 0.25,\quad \mu(x,y)=\left\{\begin{array}{ll}
             1\;\text{km}/\text{s}, & \text{if } y\geq4 \\
             1.5\;\text{km}/\text{s}, &   \text{if } (y<4 \text{ and } x\leq6) \text{ or } (16/3< y<4 \text{ and } x>6)\\
             2.25\;\text{km}/\text{s}, &   \text{if } 6\leq x \text{ and } y\leq16/3
        \end{array}\right.\\
        \lambda(x,y)&=\left\{\begin{array}{ll}
             8\;km/s, & \text{if } y\geq4 \\
             12\;km/s, &   \text{if } (y<4 \text{ and } x\leq6) \text{ or } (16/3< y<4 \text{ and } x>6)\\
             18\;km/s, &   \text{if } 6\leq x \text{ and } y\leq16/3
        \end{array}\right.,\quad u_0\equiv0,\\
        f(x,y,t)&=\left\{\begin{array}{ll}
            0, & \text{if } \|(x,y)-(4,2)\|\geq 0.01\\
            e^{\frac{\|(x,y)-(4,2)\|^2}{\|(x,y)-(4,2)\|^2-0.01^2}}(1-f_0^2\pi(t-t_0)^2)e^{-f_0^2\pi^2(t-t_0)^2}, & \text{if } \|(x,y)-(4,2)\|<0.01
        \end{array}\right.
    \end{align*}

\subsection{Stability and dispersion}\label{sec_appendix_stability_dispersion}

Here we present the operators and results of stability and dispersion for all the systems of equations considered in Section \ref{sec_appendix_numerical_framework} (assuming no PML and no source term), with a spatial discretization of fourth and eighth orders.\\

\begin{enumerate}

\item 1SD and 2SD in one dimension
\begin{equation*}
\Delta t\boldsymbol{H}=\frac{c\Delta t}{\Delta x}\begin{pmatrix}
0&g_{11}\\
g_{21}&0
\end{pmatrix},\qquad
\Delta t\boldsymbol{H}=\frac{c\Delta t}{\Delta x}\begin{pmatrix}
0&1\\
h_{11}&0
\end{pmatrix}.
\end{equation*}

    \item 1SD and 2SD in two dimension
    \begin{equation*}
        \Delta t\boldsymbol{H}=\frac{c\Delta t}{\Delta x}\begin{pmatrix}
        0&g_{11}&g_{12}\\
        g_{21}&0&0\\
        g_{22}&0&0
        \end{pmatrix},\qquad \Delta t\boldsymbol{H}=\frac{c\Delta t}{\Delta x}\begin{pmatrix}
        0&1\\
        h_{11}+h_{22}&0
        \end{pmatrix}.
    \end{equation*}
    \item elastic in two dimension (without considering the decoupled two first equations)
    \begin{equation*}
        \Delta t\boldsymbol{H}=\frac{\Delta t}{\Delta x}\frac{2\mu+\lambda}{\rho}\begin{pmatrix}
        0&0&\frac{1}{2\mu+\lambda}g_{11}&\frac{1}{2\mu+\lambda}g_{12}&0\\
        0&0&0&\frac{1}{2\mu+\lambda}g_{21}&\frac{1}{2\mu+\lambda}g_{22}\\
        \rho g_{21}&\frac{\rho\lambda}{2\mu+\lambda} g_{12}&0&0&0\\
        \frac{\rho\mu}{2\mu+\lambda} g_{22}&\frac{\rho\mu}{2\mu+\lambda} g_{11}&0&0&0\\
        \frac{\rho\lambda}{2\mu+\lambda} g_{21}&\rho g_{12}&0&0&0
        \end{pmatrix}.
    \end{equation*}
\end{enumerate}
Where
\begin{enumerate}
    \item For 4th order
    \begin{flalign*}
        g_{11}&=\frac{1}{24}\left(27(1-e^{-ik_x\Delta x})+e^{-2k_x\Delta x}-e^{k_x\Delta x}\right)\\
        g_{12}&=\frac{1}{24}\left(27(1-e^{-ik_y\Delta x})+e^{-2k_y\Delta x}-e^{k_y\Delta x}\right)\\
        g_{21}&=\frac{1}{24}\left(27(e^{ik_x\Delta x}-1)+e^{-k_x\Delta x}-e^{2k_x\Delta x}\right)\\
        g_{22}&=\frac{1}{24}\left(27(e^{ik_x\Delta x}-1)+e^{-k_x\Delta x}-e^{2k_x\Delta x}\right)\\
        h_{11}&=-\frac{1}{6}\cos(2\theta_x)+\frac{8}{3}\cos(\theta_x)-\frac{5}{2}\\
        h_{22}&=-\frac{1}{6}\cos(2\theta_y)+\frac{8}{3}\cos(\theta_y)-\frac{5}{2}
    \end{flalign*}

    \item For 8th order
    \begin{flalign*}
        g_{11}&=\frac{1225}{1024}\left(1-e^{-ik_x\Delta x}+\frac{1}{15}(e^{-2k_x\Delta x}-e^{k_x\Delta x})+\frac{1}{125}(e^{2k_x\Delta x}-e^{-3k_x\Delta x})\right.\\
        &\left.+\frac{1}{1715}(e^{-4k_x\Delta x}-e^{3k_x\Delta x})\right)\\
        g_{12}&=\frac{1225}{1024}\left(1-e^{-ik_y\Delta x}+\frac{1}{15}(e^{-2k_y\Delta x}-e^{k_y\Delta x})+\frac{1}{125}(e^{2k_y\Delta x}-e^{-3k_y\Delta x})\right.\\
        &\left.+\frac{1}{1715}(e^{-4k_y\Delta x}-e^{3k_y\Delta x})\right)\\
        g_{21}&=\frac{1225}{1024}\left(e^{ik_x\Delta x}-1+\frac{1}{15}(e^{-k_x\Delta x}-e^{2k_x\Delta x})+\frac{1}{125}(e^{3k_x\Delta x}-e^{-2k_x\Delta x})\right.\\
        &\left.+\frac{1}{1715}(e^{-3k_x\Delta x}-e^{4k_x\Delta x})\right)\\
        g_{22}&=\frac{1225}{1024}\left(e^{ik_y\Delta x}-1+\frac{1}{15}(e^{-k_y\Delta x}-e^{2k_y\Delta x})+\frac{1}{125}(e^{3k_y\Delta x}-e^{-2k_y\Delta x})\right.\\
        &\left.+\frac{1}{1715}(e^{-3k_y\Delta x}-e^{4k_y\Delta x})\right)\\
        h_{11}&=-\frac{1}{560}(e^{-4k_x\Delta x}+e^{4k_x\Delta x})+\frac{8}{315}(e^{-3k_x\Delta x}+e^{3k_x\Delta x})-\frac{1}{5}(e^{-2k_x\Delta x}+e^{2k_x\Delta x})\\
        &+\frac{8}{5}(e^{-k_x\Delta x}+e^{k_x\Delta x})-\frac{205}{72}\\
        h_{22}&=-\frac{1}{560}(e^{-4k_y\Delta x}+e^{4k_y\Delta x})+\frac{8}{315}(e^{-3k_y\Delta x}+e^{3k_y\Delta x})-\frac{1}{5}(e^{-2k_y\Delta x}+e^{2k_y\Delta x})\\
        &+\frac{8}{5}(e^{-k_y\Delta x}+e^{k_y\Delta x})-\frac{205}{72}
\end{flalign*}
\end{enumerate}

\end{document}